\documentclass{amsart}
\usepackage{latexsym}
\usepackage{amsmath}
\usepackage{amssymb}
\usepackage{amsfonts}
\usepackage{amsthm}
\usepackage[all]{xy}

\newtheorem{thm}{Theorem}[section]
\newtheorem*{thm*}{Theorem}
\newtheorem{Definition}[thm]{Definition}
\newtheorem{Example}[thm]{Example}
\newtheorem{Remark}[thm]{Remark}
\newtheorem{Proposition}[thm]{Proposition}
\newtheorem{Corollary}[thm]{Corollary}
\newtheorem{Theorem}[thm]{Theorem}
\newtheorem{Hypothesis}[thm]{Hypothesis}
\newtheorem{Lemma}[thm]{Lemma}

\newcommand {\La}    {\ensuremath{\mathcal{L}}}
\newcommand {\Ca}    {\ensuremath{\mbox{$\mathcal{C}$}}}
\newcommand {\Na}    {\ensuremath{\mbox{$\mathcal{N}$}}}
\newcommand {\proj}  {\ensuremath{\mathbb{P}}}
\newcommand {\real}  {\ensuremath{\mathbb{R}}}
\newcommand {\Z}  {\ensuremath{\mathbb{Z}}}
\newcommand {\Hom}   {\ensuremath{\operatorname{Hom}}}

\newcommand {\sym}   {\ensuremath{\operatorname{\text{Sym}}}}
\newcommand {\qad}   {\ensuremath{\operatorname{\text{Quad}}}}
\newcommand {\unil}  {\hbox{\rm UNil}}

\begin{document}


\title{Generalized Arf invariants in algebraic $L$-theory}

\author{Markus Banagl}

\address{Department of Mathematical Sciences\newline
\indent University of Cincinnati\newline
\indent Cincinnati, OH 45221-0025\newline
\indent U.S.A.}
\email{banagl@math.uc.edu}

\thanks{The first author was in part supported by the
 European Union TMR network project on ``K-theory, linear
 algebraic groups and related structures.''}

\author{Andrew Ranicki}

\address{School of Mathematics\newline
\indent University of Edinburgh\newline
\indent King's Buildings\newline
\indent Edinburgh EH9 3JZ\newline
\indent SCOTLAND, UK}

\email{aar@maths.ed.ac.uk}



\begin{abstract}
The difference between the quadratic $L$-groups $L_*(A)$ and the
symmetric $L$-groups $L^*(A)$ of a ring with involution $A$ is detected
by generalized Arf invariants.  The special case $A=\Z[x]$\ gives a
complete set of invariants for the Cappell UNil-groups ${\rm
UNil}_*(\Z;\Z,\Z)$ for the infinite dihedral group
$D_{\infty}=\Z_2*\Z_2$, extending the results of Connolly and
Ranicki \cite{connran}, Connolly and Davis \cite{conndav}.
\end{abstract}


\maketitle


\tableofcontents


\section*{Introduction}

The invariant of Arf \cite{arf} is a basic ingredient in the
isomorphism classification of quadratic forms over a field of
characteristic 2.  The algebraic $L$-groups of a ring with involution
$A$ are Witt groups of quadratic structures on $A$-modules and
$A$-module chain complexes, or equivalently the cobordism groups of
algebraic Poincar\'e complexes over $A$.  The cobordism formulation of
algebraic $L$-theory is used here to obtain generalized Arf invariants
detecting the difference between the quadratic and symmetric $L$-groups
of an arbitrary ring with involution $A$, with applications to the
computation of the Cappell $\unil$-groups.

The (projective) {\it quadratic $L$-groups} of Wall \cite{wall}
are 4-periodic groups
$$L_n(A)~=~L_{n+4}(A)~.$$
The $2k$-dimensional $L$-group $L_{2k}(A)$ is the Witt group
of nonsingular $(-1)^k$-quadratic forms $(K,\psi)$ over $A$,
where $K$ is a f.g. projective $A$-module
and $\psi$ is an equivalence class of $A$-module morphisms
$$\psi~:~K \to K^*~=~{\rm Hom}_A(K,A)$$
such that $\psi+(-1)^k\psi^*:K \to K^*$
is an isomorphism, with $\psi \sim \psi+\chi+(-1)^{k+1}\chi^*$
for $\chi \in {\rm Hom}_A(K,K^*)$. A lagrangian $L$ for $(K,\psi)$
is a direct summand $L \subset K$ such that
$$\begin{array}{l}
L^{\perp}~=~L, ~{\rm where}~
L^{\perp}~=~\{x \in K\,\vert\,(\psi+(-1)^k\psi^*)(x)(y)=0~{\rm for~all}~
y \in L\}~,\\[1ex]
\psi(x)(x) \in \{a+(-1)^{k+1}\overline{a}\,\vert\, a \in A\}~{\rm for~all}~
x \in L~.
\end{array}$$
A form $(K,\psi)$ admits a lagrangian $L$ if and only if it is isomorphic
to the hyperbolic form $H_{(-1)^k}(L)=(L\oplus L^*,
\begin{pmatrix} 0 & 1 \\ 0 & 0 \end{pmatrix})$, in which case
$$(K,\psi)~=~H_{(-1)^k}(L)~=~0 \in L_{2k}(A)~.$$
The $(2k+1)$-dimensional $L$-group
$L_{2k+1}(A)$ is the Witt group of $(-1)^k$-quadratic
formations $(K,\psi;L,L')$ over $A$, with $L,L' \subset K$
lagrangians for $(K,\psi)$.

The {\it symmetric $L$-groups}
$L^n(A)$ of Mishchenko \cite{mish} are the cobordism
groups of {\it $n$-dimensional symmetric Poincar\'e complexes}
$(C,\phi)$ over $A$, with $C$ an $n$-dimensional f.g. projective
$A$-module chain complex
$$C~:~\dots \to 0 \to C_n \to C_{n-1} \to \dots \to C_1 \to C_0 \to 0 \to \dots$$
and $\phi \in Q^n(C)$ an element of the $n$-dimensional symmetric $Q$-group of $C$
(about which more in \S1 below) such that $\phi_0:C^{n-*} \to C$ is a chain
equivalence.  In particular, $L^0(A)$ is the Witt group of nonsingular
symmetric forms $(K,\phi)$ over $A$, with
$$\phi~=~\phi^*~:~K \to K^*$$
an isomorphism, and  $L^1(A)$ is the Witt group
of symmetric formations $(K,\phi;L,L')$ over $A$.
An analogous cobordism formulation of the quadratic $L$-groups was
obtained in Ranicki \cite{ranicki1}, expressing $L_n(A)$ as the
cobordism group of {\it $n$-dimensional quadratic Poincar\'e complexes}
$(C,\psi)$, with $\psi \in Q_n(C)$ an element of the $n$-dimensional
quadratic $Q$-group
of $C$ such that $(1+T)\psi_0:C^{n-*} \to C$ is a chain equivalence.
The {\it hyperquadratic $L$-groups} $\widehat{L}^n(A)$ of \cite{ranicki1} are the
cobordism groups of $n$-dimensional (symmetric, quadratic) Poincar\'e
pairs $(f:C \to D,(\delta\phi,\psi))$ over $A$  such that
$$(\delta\phi_0,(1+T)\psi_0)~:~D^{n-*} \to \Ca(f)$$
is a chain equivalence, with $\Ca(f)$ the algebraic mapping cone of $f$.
The various $L$-groups are related by an exact sequence
    $$\xymatrix@C-5pt{ \dots \ar[r] & L_n(A) \ar[r]^{1+T} & L^n(A)
    \ar[r] & \widehat{L}^n(A)\ar[r]^-{\partial} &  L_{n-1}(A) \ar[r] & \dots~.}$$
The symmetrization maps $1+T:L_*(A) \to L^*(A)$ are isomorphisms modulo
8-torsion, so that the hyperquadratic $L$-groups $\widehat{L}^*(A)$ are
of exponent 8.  The symmetric and hyperquadratic $L$-groups are not
4-periodic in general. However, there are defined natural maps
$$L^n(A) \to L^{n+4}(A)~,~\widehat{L}^n(A) \to \widehat{L}^{n+4}(A)$$
(which are isomorphisms modulo 8-torsion), and there are 4-periodic
versions of the $L$-groups
$$L^{n+4*}(A)~=~\lim\limits_{k\to\infty}L^{n+4k}(A)~,~
\widehat{L}^{n+4*}(A)~=~\lim\limits_{k\to\infty}\widehat{L}^{n+4k}(A)~.$$
The 4-periodic symmetric $L$-group $L^{n+4*}(A)$ is the cobordism group
of $n$-dimensional symmetric Poincar\'e complexes $(C,\phi)$ over $A$
with $C$ a finite (but not necessarily $n$-dimensional) f.g. projective
$A$-module chain complex, and similarly for $\widehat{L}^{n+4*}(A)$.

The Tate $\Z_2$-cohomology groups of a ring with involution $A$
$$\widehat{H}^n(\Z_2;A)~=~\dfrac{\{x \in A\,\vert\,\overline{x}=(-1)^nx\}}
{\{y+(-1)^n\overline{y}\,\vert\,y \in A\}}~~(n(\bmod\,2))$$
are $A$-modules via
$$A \times \widehat{H}^n(\Z_2;A) \to \widehat{H}^n(\Z_2;A)~;~
(a,x) \mapsto ax \overline{a}~.$$
The Tate $\Z_2$-cohomology $A$-modules give an indication of the
difference between the quadratic and symmetric $L$-groups of $A$.
If $\widehat{H}^*(\Z_2;A)=0$ (e.g. if $1/2 \in A$) then
the symmetrization maps $1+T:L_*(A) \to L^*(A)$ are isomorphisms
and $\widehat{L}^*(A)=0$. If $A$ is such that $\widehat{H}^0(\Z_2;A)$
and $\widehat{H}^1(\Z_2;A)$ have 1-dimensional f.g.  projective
$A$-module resolutions then the symmetric and hyperquadratic $L$-groups
of $A$ are 4-periodic (Proposition \ref{4period}).

We shall say that a ring with the involution $A$ is {\it $r$-even} for some
$r \geqslant 1$ if
\begin{itemize}
\item[(i)] $A$ is commutative with the identity involution,
so that $\widehat{H}^0(\Z_2;A)=A_2$ as an additive group with
$$A \times \widehat{H}^0(\Z_2;A) \to \widehat{H}^0(\Z_2;A)~;~(a,x)
\mapsto a^2x~,$$
and
$$\widehat{H}^1(\Z_2;A)~=~\{a \in A\,\vert\,2a=0\}~,$$
\item[(ii)] $2 \in A$ is a non-zero divisor, so that $\widehat{H}^1(\Z_2;A)=0$,
\item[(iii)] $\widehat{H}^0(\Z_2;A)$ is a f.g. free $A_2$-module of rank $r$
with a basis $\{x_1=1,x_2,\dots,x_r\}$.
\end{itemize}
If $A$ is $r$-even then $\widehat{H}^0(\Z_2;A)$
has a 1-dimensional f.g. free $A$-module resolution
$$0 \to A^r \xymatrix{\ar[r]^-{2}&} A^r \xymatrix{\ar[r]^-{x}&}
\widehat{H}^0(\Z_2;A) \to 0$$
so that the symmetric and hyperquadratic $L$-groups of $A$ are
4-periodic (\ref{4period}).\\
A ring with involution $A$ is 1-even if and only if it
satisfies (i), (ii) and also
$$a-a^2 \in 2A~\hbox{for all}~a \in A~.$$

\begin{Theorem} \label{main1}
The hyperquadratic $L$-groups of a 1-even ring with involution $A$ are given by~:
$$\widehat{L}^n(A)~=~\begin{cases}
A_8&\hbox{\it if $n \equiv 0(\bmod\, 4)$}\\
A_2&\hbox{\it if $n \equiv 1,3(\bmod\, 4)$}\\
0&\hbox{\it if $n \equiv 2(\bmod\, 4)$}~.
\end{cases}$$
The boundary maps $\partial:\widehat{L}^n(A) \to L_{n-1}(A)$  are given by~:
$$\begin{array}{l}
\partial~:~\widehat{L}^0(A)=A_8 \to L_{-1}(A)~;~a
\mapsto (A\oplus A,\begin{pmatrix} 0 & 1 \\ 0 & 0 \end{pmatrix}
;A,{\rm im}(\begin{pmatrix} 1-a \\ a \end{pmatrix}:
A \to A \oplus A))~,\\[2ex]
\partial~:~\widehat{L}^1(A)=A_2 \to L_0(A)~;~a
\mapsto (A\oplus A,\begin{pmatrix} (a-a^2)/2 & 1-2a  \\ 0 & -2 \end{pmatrix})~,\\[2ex]
\partial~:~\widehat{L}^3(A)=A_2 \to L_2(A)~;~a \mapsto
(A\oplus A , \begin{pmatrix} a & 1 \\ 0 & 1 \end{pmatrix})~.
\end{array}$$
The map
$$L^0(A) \to \widehat{L}^0(A)~=~A_8~;~(K,\phi) \mapsto \phi(v,v)$$
is defined using any element $v \in K$ such that
$$\phi(u,u)~=~\phi(u,v) \in A_2~~(u \in K)~.$$
\hfill$\qed$
\end{Theorem}

Theorem \ref{main1} is proved in \S2 (Corollary \ref{cor1}).  In
particular, $A=\Z$ is 1-even, and in this case Theorem \ref{main1}
recovers the computation of $\widehat{L}^*(\Z)$ obtained in
\cite{ranicki1} -- the algebraic $L$-theory of $\Z$ is recalled further
below in the Introduction.

\begin{Theorem} \label{main2}
If $A$ is 1-even then the polynomial ring $A[x]$ is 2-even,
with $A[x]_2$-module basis $\{1,x\}$ for $\widehat{H}^0(\Z_2;A[x])$.
The hyperquadratic $L$-groups of $A[x]$ are given by~:
$$\widehat{L}^n(A[x])~=~\begin{cases}
A_8\oplus A_4[x] \oplus A_2[x]^3&\hbox{\it if $n \equiv 0(\bmod\, 4)$}\\
A_2&\hbox{\it if $n \equiv 1(\bmod\, 4)$}\\
0&\hbox{\it if $n \equiv 2(\bmod\, 4)$}\\
A_2[x]&\hbox{\it if $n \equiv 3(\bmod\, 4)$}~.
\end{cases}$$
\hfill$\qed$
\end{Theorem}

Theorems \ref{main1} and \ref{main2} are special cases of the following
computation~:

\begin{Theorem} \label{main3}
The hyperquadratic $L$-groups of an $r$-even ring with involution
$A$ are given by~:
$$\begin{array}{l}
\widehat{L}^n(A)~=\\[3ex]
\begin{cases}
\dfrac{\{M \in {\rm Sym}_r(A)\,\vert\, M-MXM \in {\rm Quad}_r(A)\}}
{4{\rm Quad}_r(A) + \{2(N+N^t)-N^tXN\,\vert\,N \in M_r(A)\}}&
\hbox{\it if $n=0$}\\[3ex]
\dfrac{\{N \in M_r(A)\,\vert\,N+N^t \in 2{\rm Sym}_r(A),
\dfrac{1}{2}(N+N^t)-N^tXN \in {\rm Quad}_r(A)\}}{2M_r(A)}&\text{if $n=1$}\\[2ex]
0&\hbox{\it if $n=2$}\\[2ex]
\dfrac{{\rm Sym}_r(A)}
{{\rm Quad}_r(A) + \{L-LXL\,\vert\,L \in {\rm Sym}_r(A)\}}&
\hbox{\it if $n=3$}
\end{cases}
\end{array}$$
with ${\rm Sym}_r(A)$ the additive group of symmetric $r \times r$ matrices
$(a_{ij})=(a_{ji})$ in $A$, ${\rm Quad}_r(A) \subset {\rm Sym}_r(A)$ the subgroup of
the matrices such that $a_{ii} \in 2A$, and
$$X~=~\begin{pmatrix} x_1 & 0 & \dots & 0 \\ 0 & x_2 & \dots & 0 \\
\vdots & \vdots & \ddots & 0 \\
0 & 0 & \dots & x_r \end{pmatrix} \in {\rm Sym}_r(A)$$
for an $A_2$-module basis
$\{x_1=1,x_2,\dots,x_r\}$ of $\widehat{H}^0(\Z_2;A)$.
The boundary maps $\partial:\widehat{L}^n(A) \to L_{n-1}(A)$  are given by~:
$$\begin{array}{l}
\partial~:~\widehat{L}^0(A) \to L_{-1}(A)~;~
M \mapsto (H_{-}(A^r);A^r,{\rm im}(\begin{pmatrix}
1-XM \\ M \end{pmatrix}:A^r \to A^r\oplus (A^r)^*))~,\\[2ex]
\partial~:~\widehat{L}^1(A) \to L_0(A)~;~N \mapsto
(A^r \oplus A^r, \begin{pmatrix}
\dfrac{1}{4}(N+N^t-2N^tXN) & 1-2NX \\
0 & -2X \end{pmatrix})~,\\[1ex]
\partial~:~\widehat{L}^3(A) \to L_2(A)~;~M \mapsto
(A^r \oplus (A^r)^*,\begin{pmatrix} M & 1 \\ 0 & X \end{pmatrix})~.
\end{array}$$
\hfill$\qed$
\end{Theorem}

In \S\S1,2 we recall and extend the $Q$-groups and algebraic chain
bundles of Ranicki \cite{ranicki1}, \cite{ranicki4} and Weiss
\cite{weiss1}, including a proof of Theorem \ref{main3} (Theorem
\ref{equ.ba-ba0}).

We shall be dealing with two types of generalized Arf invariant: for
forms on f.g.  projective modules, and for linking forms on homological
dimension 1 torsion modules, which we shall be considering separately.

In \S3 we define the {\it generalized Arf invariant} of a nonsingular
$(-1)^k$-quadratic form $(K,\psi)$ over an arbitrary ring with
involution $A$ with a
lagrangian $L \subset K$ for $(K,\psi+(-)^k\psi^*)$ to be the element
$$(K,\psi;L) \in \widehat{L}^{4*+2k+1}(A)$$
with image
$$(K,\psi) \in {\rm im}(\partial:\widehat{L}^{4*+2k+1}(A) \to L_{2k}(A))
~=~{\rm ker}(1+T:L_{2k}(A) \to L^{4*+2k}(A))~.$$
Theorem \ref{genArfthm} gives an explicit formula for the generalized
Arf invariant $(K,\psi;L) \in \widehat{L}^3(A)$ for an $r$-even $A$.
Generalizations of the Arf invariants in $L$-theory have been previously
studied by Clauwens \cite{clauwens} and Bak \cite{bak}.

In \S4 we consider a ring with involution $A$ with a localization
$S^{-1}A$ inverting a multiplicative subset $S \subset A$ of
central non-zero divisors such that $\widehat{H}^*(\Z_2;S^{-1}A)=0$
(e.g. if $2 \in S$).
The relative $L$-group $L_{2k}(A,S)$ in the localization exact sequence
$$\dots \to L_{2k}(A) \to L_{2k}(S^{-1}A) \to L_{2k}(A,S)
\to L_{2k-1}(A) \to L_{2k-1}(S^{-1}A) \to \dots$$
is the Witt group of nonsingular $(-1)^k$-quadratic linking forms
$(T,\lambda,\mu)$ over $(A,S)$, with $T$ a homological dimension 1
$S$-torsion $A$-module, $\lambda$ an $A$-module isomorphism
$$\lambda~=~(-1)^k\lambda\widehat{~}~:~T \to
T\widehat{~}~=~{\rm Ext}^1_A(T,A)~=~{\rm Hom}_A(T,S^{-1}A/A)$$
and
$$\mu~:~T \to Q_{(-1)^k}(A,S)~=~\dfrac{
\{b \in S^{-1}A \,\vert\, \overline{b}=(-1)^kb\}}
{\{a+(-1)^k\overline{a}\,\vert\,a\in A\}}$$
a $(-1)^k$-quadratic function for $\lambda$. The {\it linking Arf invariant} of a
nonsingular $(-1)^k$-quadratic linking form $(T,\lambda,\mu)$ over
$(A,S)$ with a lagrangian $U \subset T$ for $(T,\lambda)$ is
defined to be an element
$$(T,\lambda,\mu;U) \in \widehat{L}^{4*+2k}(A)~,$$
with properties analogous to the generalized Arf invariant defined
for forms in \S3. Theorem \ref{linkArfthm} gives an explicit formula for
the linking Arf invariant
$(T,\lambda,\mu;U) \in \widehat{L}^{2k}(A)$ for an $r$-even $A$, using
$$S~=~(2)^{\infty}~=~\{2^i\,\vert\,i \geqslant 0\} \subset A~~,~~
S^{-1}A~=~A[1/2]~.$$

In \S5 we apply the generalized and linking Arf invariants to the
algebraic $L$-groups of a polynomial extension $A[x]$
($\overline{x}=x$) of a ring with involution $A$, using the exact sequence
   $$\xymatrix@C-5pt{ \dots \ar[r] & L_n(A[x]) \ar[r]^{1+T} & L^n(A[x])
    \ar[r] & \widehat{L}^n(A[x])\ar[r] &  L_{n-1}(A[x]) \ar[r] & \dots~.}$$
For a Dedekind ring $A$ the quadratic $L$-groups of $A[x]$ are related to the
UNil-groups ${\rm UNil}_*(A)$ of Cappell \cite{cappell2} by
the splitting formula of Connolly and Ranicki \cite{connran}
$$L_n(A[x])~=~L_n(A) \oplus {\rm UNil}_n(A)~,$$
and the symmetric and hyperquadratic $L$-groups of $A[x]$ are
4-periodic, and such that
$$L^n(A[x])~=~L^n(A)~,~
\widehat{L}^{n+1}(A[x])~=~\widehat{L}^{n+1}(A) \oplus {\rm UNil}_n(A)~.$$
Any computation of $\widehat{L}^*(A)$ and $\widehat{L}^*(A[x])$ thus
gives a computation of ${\rm UNil}_*(A)$.  Combining the splitting
formula with Theorems \ref{main1}, \ref{main2} gives~:

\begin{Theorem} \label{main4} If $A$ is a 1-even Dedekind ring then
$$\begin{array}{ll}
{\rm UNil}_n(A)&=~\widehat{L}^{n+1}(A[x])/\widehat{L}^{n+1}(A)\\[1ex]
&=~\begin{cases}
0&\hbox{\rm if $n \equiv 0,1(\bmod\, 4)$}\\
xA_2[x]&\hbox{\rm if $n \equiv 2(\bmod\, 4)$}\\
A_4[x] \oplus A_2[x]^3&\hbox{\rm if $n \equiv 3(\bmod\, 4)$}~.
\end{cases}
\end{array}$$
\hfill$\qed$
\end{Theorem}

In particular, Theorem \ref{main4} applies to $A=\Z$.
The twisted quadratic $Q$-groups were first used
in the partial computation of
$${\rm UNil}_n(\Z)~=~\widehat{L}^{n+1}(\Z[x])/\widehat{L}^{n+1}(\Z)$$
by Connolly and Ranicki \cite{connran}. The calculation
in \cite{connran} was almost complete, except that ${\rm UNil}_3(\Z)$ was
only obtained up to extensions.  The calculation was first completed by
Connolly and Davis \cite{conndav}, using linking forms.  We are
grateful to them for sending us a preliminary version of their paper.
The calculation of ${\rm UNil}_3(\Z)$ in \cite{conndav} uses the
results of \cite{connran} and the classifications of quadratic and symmetric
linking forms over $(\Z[x],(2)^{\infty})$.  The calculation of ${\rm
UNil}_3(\Z)$ here uses the linking Arf invariant measuring the
difference between the Witt groups of  quadratic and symmetric
linking forms over $(\Z[x],(2)^{\infty})$, developing further
the $Q$-group strategy of \cite{connran}.

The algebraic $L$-groups of $A=\Z_2$ are given by~:
$$\begin{array}{l}
L^n(\Z_2)~=~\begin{cases}
\Z_2~\hbox{(rank ($\bmod\, 2$))}&\hbox{if $n \equiv 0(\bmod\, 2)$}\\
0&\hbox{if $n \equiv 1(\bmod\, 2)$}~,
\end{cases}\\[4ex]
L_n(\Z_2)~=~\begin{cases}
\Z_2~\hbox{(Arf invariant)}&\hbox{if $n \equiv 0(\bmod\, 2)$}\\
0&\hbox{if $n \equiv 1(\bmod\, 2)$}~,
\end{cases}\\[4ex]
\widehat{L}^n(\Z_2)~=~\Z_2
\end{array}$$
with $1+T=0:L_n(\Z_2) \to L^n(\Z_2)$.
The classical Arf invariant is defined for a nonsingular quadratic
form $(K,\psi)$ over $\Z_2$ and a lagrangian $L\subset K$ for the
symmetric form $(K,\psi+\psi^*)$ to be
$$(K,\psi;L)~=~\sum\limits^{\ell}_{i=1}\psi(e_i,e_i).\psi(e^*_i,e^*_i) \in
\widehat{L}^1(\Z_2)~=~L_0(\Z_2)~=~\Z_2$$
with $e_1,e_2,\dots,e_{\ell}$ any basis for $L\subset K$, and
$e^*_1,e^*_2,\dots,e^*_{\ell}$ a basis for a direct summand
$L^* \subset K$ such that
$$(\psi+\psi^*)(e^*_i,e^*_j)~=~0~,~(\psi+\psi^*)(e^*_i,e_j)~=~\begin{cases}
1&\hbox{if $i=j$}\\
0&\hbox{if $i\neq j$~.}
\end{cases}$$
The Arf invariant is independent of the choices of $L$ and $L^*$.

The algebraic $L$-groups of $A=\Z$ are given by~:
$$\begin{array}{l}
L^n(\Z)~=~\begin{cases}
\Z~\hbox{(signature)}&\hbox{if $n \equiv 0(\bmod\, 4)$}\\
\Z_2~\hbox{(de Rham invariant)}&\hbox{if $n \equiv 1(\bmod\, 4)$}\\
0&\hbox{otherwise~,}
\end{cases}\\[6ex]
L_n(\Z)~=~\begin{cases}
\Z~\hbox{(signature/8)}&\hbox{if $n \equiv 0(\bmod\, 4)$}\\
\Z_2~\hbox{(Arf invariant)}&\hbox{if $n \equiv 2(\bmod\, 4)$}\\
0&\hbox{otherwise~,}
\end{cases}\\[6ex]
\widehat{L}^n(\Z)~=~\begin{cases}
\Z_8~\hbox{(signature ($\bmod\, 8$))}&\hbox{if $n \equiv 0(\bmod\, 4)$}\\
\Z_2~\hbox{(de Rham invariant)}&\hbox{if $n \equiv 1(\bmod\, 4)$}\\
0&\hbox{if $n \equiv 2(\bmod\, 4)$}\\
\Z_2~\hbox{(Arf invariant)}&\hbox{if $n \equiv 3(\bmod\, 4)$~.}
\end{cases}
\end{array}$$
Given a nonsingular symmetric form $(K,\phi)$ over $\Z$
there is a congruence (Hirzebruch, Neumann and Koh \cite[Theorem 3.10]{hirz})
$${\rm signature}(K,\phi)~\equiv~\phi(v,v)~(\bmod\, 8)$$
with $v \in K$ any element such that
$$\phi(u,v)~\equiv~\phi(u,u) ~(\bmod\, 2)~~(u \in K)~,$$
so that
$$\begin{array}{l}
(K,\phi)~=~{\rm signature}(K,\phi)~=~\phi(v,v)\\[1ex]
\hphantom{(K,\phi)~=~} \in
{\rm coker}(1+T:L_0(\Z) \to L^0(\Z))~=~\widehat{L}^0(\Z)~=~
{\rm coker}(8:\Z \to \Z)~=~\Z_8~.
\end{array}$$
The projection $\Z \to \Z_2$ induces an isomorphism
$L_2(\Z) \cong L_2(\Z_2)$, so that the Witt class
of a nonsingular $(-1)$-quadratic
form $(K,\psi)$ over $\Z$ is given by the Arf invariant of the mod 2 reduction
$$(K,\psi;L)~=~\Z_2\otimes_\Z(K,\psi;L) \in L_2(\Z)~=~L_2(\Z_2)~=~\Z_2$$
with $L \subset K$ a lagrangian for the $(-1)$-symmetric form $(K,\psi-\psi^*)$.
Again, this is independent of the choice of $L$.

The $Q$-groups are defined for an $A$-module chain complex $C$ to be
$\Z_2$-hyperhomology invariants of the $\Z[\Z_2]$-module chain complex
$C \otimes_AC$.  The involution on $A$ is used to define the tensor
product over $A$ of left $A$-module chain complexes $C,D$, the abelian
group chain complex
$$C\otimes_AD~=~\displaystyle{\frac{C\otimes_{\Z}D}{
\{ax \otimes y - x\otimes \overline{a}y \,\vert\, a \in A,x \in C,y \in D\}}}~.$$
Let $C\otimes_AC$ denote the
$\Z[\Z_2]$-module chain complex defined by $C\otimes_AC$ via the
transposition involution
$$T~:~C_p \otimes_AC_q \to C_q \otimes_AC_p~;~
x \otimes y \mapsto (-1)^{pq} y \otimes x~.$$
The $\begin{cases} \hbox{\it symmetric} \cr
\hbox{\it quadratic} \cr
\hbox{\it hyperquadratic} \end{cases}$
{\it $Q$-groups} of $C$ are defined by
$$\begin{cases} Q^n(C)~=~H^n(\Z_2;C\otimes_AC) \cr
Q_n(C)~=~H_n(\Z_2;C\otimes_AC) \cr
\widehat{Q}^n(C)~=~\widehat{H}^n(\Z_2;C\otimes_AC)~.
\end{cases}$$
The $Q$-groups are covariant in $C$, and are chain homotopy invariant.
The $Q$-groups are related by an exact sequence
    $$\xymatrix@C-5pt{ \dots \ar[r] & Q_n(C) \ar[r]^-{1+T} & Q^n(C)
    \ar[r]^-{J} & \widehat{Q}^n(C)\ar[r]^-{H} &  Q_{n-1}(C) \ar[r] & \dots~.}$$
\indent
A {\it chain bundle} $(C,\gamma)$ over $A$ is an $A$-module chain complex $C$
together with an element $\gamma \in \widehat{Q}^0(C^{-*})$.
The {\it twisted quadratic $Q$-groups} $Q_*(C,\gamma)$ were defined
in Weiss \cite{weiss1} using simplicial abelian groups, to fit into
an exact sequence
    $$\xymatrix@C-5pt{ \dots \ar[r] & Q_n(C,\gamma) \ar[r]^-{N_{\gamma}} &
    Q^n(C) \ar[r]^-{J_{\gamma}} & \widehat{Q}^n(C)\ar[r]^-{H_\gamma} &
    Q_{n-1}(C,\gamma) \ar[r] & \dots}$$
with
$$J_{\gamma} ~:~Q^n(C) \to \widehat{Q}^n(C)~;~\phi \mapsto J(\phi)
-(\widehat{\phi}_0)^{\%}(\gamma)~.$$

An {\it $n$-dimensional algebraic normal complex} $(C,\phi,\gamma,\theta)$
over $A$ is an $n$-dimensional symmetric complex $(C,\phi)$ together with
a chain bundle $\gamma \in \widehat{Q}^0(C^{-*})$ and
an element $(\phi,\theta) \in Q_n(C,\gamma)$ with image $\phi \in Q^n(C)$.
Every $n$-dimensional symmetric Poincar\'e complex $(C,\phi)$
has the structure of an algebraic normal complex $(C,\phi,\gamma,\theta)$~:
the {\it Spivak normal chain bundle} $(C,\gamma)$ is characterized by
$$(\widehat{\phi}_0)^{\%}(\gamma)~=~J(\phi) \in Q^n(C)~,$$
with
$$(\widehat{\phi}_0)^{\%}~:~\widehat{Q}^0(C^{-*})~=~\widehat{Q}^n(C^{n-*})
\to \widehat{Q}^n(C)$$
the isomorphism induced by the Poincar\'e duality
chain equivalence $\phi_0:C^{n-*} \to C$, and the
{\it algebraic normal invariant} $(\phi,\theta) \in Q_n(C,\gamma)$ is such that
$$N_{\gamma}(\phi,\theta)~=~\phi \in Q^n(C)~.$$
See Ranicki \cite[\S7]{ranicki4} for the one-one correspondence between
the homotopy equivalence classes of $n$-dimensional (symmetric, quadratic)
Poincar\'e pairs and $n$-dimensional algebraic normal complexes.
Specifically, an $n$-dimensional algebraic normal complex $(C,\phi,\gamma,\theta)$
determines an $n$-dimensional (symmetric, quadratic) Poincar\'e pair
$(\partial C \to C^{n-*},(\delta\phi,\psi))$ with
$$\partial C~=~\Ca(\phi_0:C^{n-*} \to C)_{*+1}~.$$
Conversely, an $n$-dimensional (symmetric, quadratic) Poincar\'e pair
$(f:C \to D,(\delta\phi,\psi))$ determines an $n$-dimensional algebraic
normal complex $(\Ca(f),\gamma,\phi,\theta)$, with
$\gamma \in \widehat{Q}^0(\Ca(f)^{-*})$ the Spivak normal chain bundle
and $\phi=\delta\phi/(1+T)\psi$~;
the class $(\phi,\theta) \in Q_n(\Ca(f),\gamma)$ is the algebraic
normal invariant of $(f:C \to D,(\delta\phi,\psi))$.
Thus $\widehat{L}^n(A)$ is the cobordism group of $n$-dimensional
normal complexes over $A$.

Weiss \cite{weiss1} established that for any ring with involution $A$
there exists a universal chain bundle $(B^A,\beta^A)$ over $A$, such
that every chain bundle $(C,\gamma)$ is classified by a chain bundle map
$$(g,\chi)~:~(C,\gamma) \to (B^A,\beta^A)~,$$
with
$$H_*(B^A)~=~\widehat{H}^*(\Z_2;A)~.$$
The function
$$\widehat{L}^{n+4*}(A)\to Q_n(B^A,\beta^A)~;~
(C,\phi,\gamma,\theta) \mapsto (g,\chi)_{\%}(\phi,\theta)$$
was shown in \cite{weiss1} to be an isomorphism.  Since the $Q$-groups
are homological in nature (rather than of the Witt type)
they are in principle effectively computable.  The algebraic
normal invariant defines the isomorphism
$$\begin{array}{l}
{\rm ker}(1+T:L_n(A) \to L^{n+4*}(A)) \xymatrix{\ar[r]^-{\displaystyle{\cong}}&}
{\rm coker}(L^{n+4*+1}(A)\to Q_{n+1}(B^A,\beta^A))~;\\[1ex]
\hphantom{{\rm ker}(1+T:L_n(A) \to L^{n+4*}(A))}
(C,\psi) \mapsto (g,\chi)_{\%}(\phi,\theta)
\end{array}$$
with $(\phi,\theta) \in Q_{n+1}(\Ca(f),\gamma)$ the algebraic
normal invariant of any $(n+1)$-dimensional
(symmetric, quadratic) Poincar\'e pair
$(f:C \to D,(\delta\phi,\psi))$, with
classifying chain bundle map $(g,\chi):(\Ca(f),\gamma) \to
(B^A,\beta^A)$. For $n=2k$ such a pair with $H_i(C)=H_i(D)=0$ for $i \neq k$
is just a nonsingular $(-1)^k$-quadratic form $(K=H^k(C),\psi)$ with
a lagrangian
$$L~=~{\rm im}(f^*:H^k(D) \to H^k(C)) \subset K~=~H^k(C)$$
for $(K,\psi+(-1)^k\psi^*)$,
such that the generalized Arf invariant is the image of the algebraic
normal invariant
$$(K,\psi;L)~=~(g,\chi)_{\%}(\phi,\theta) \in
\widehat{L}^{4*+2k+1}(A)~=~Q_{2k+1}(B^A,\beta^A)~.$$
For $A=\Z_2$ and $n=0$ this is just the classical Arf invariant isomorphism
$$\begin{array}{l}
L_0(\Z_2)~=~{\rm ker}(1+T=0:L_0(\Z_2) \to L^0(\Z_2))\\[1ex]
\hskip100pt \xymatrix{\ar[r]^-{\displaystyle{\cong}}&}
{\rm coker}(L^1(\Z_2)=0 \to Q_1(B^{\Z_2},\beta^{\Z_2}))~=~\Z_2~;\\[1ex]
\hskip150pt (K,\psi) \mapsto (K,\psi;L)
\end{array}$$
with $L\subset K$ an arbitrary lagrangian of $(K,\psi+\psi^*)$. The isomorphism
$${\rm coker}(1+T:L_n(A) \to L^{n+4*}(A))
\xymatrix{\ar[r]^-{\displaystyle{\cong}}&}
{\rm ker}(\partial:Q_n(B^A,\beta^A)\to L_{n-1}(A))$$
is a generalization from $A=\Z$, $n=0$ to arbitrary $A$, $n$
of the identity ${\rm signature}(K,\phi)\equiv \phi(v,v)~(\bmod\, 8)$
described above.

[Here is some of the geometric background.  Chain bundles are algebraic
analogues of vector bundles and spherical fibrations, and the twisted
$Q$-groups are the analogues of the homotopy groups of the Thom spaces.
A $(k-1)$-spherical fibration $\nu:X \to BG(k)$ over a connected $CW$
complex $X$ determines a chain bundle $(C(\widetilde{X}),\gamma)$ over
$\Z[\pi_1(X)]$, with $C(\widetilde{X})$ the cellular
$\Z[\pi_1(X)]$-module chain complex of the universal cover
$\widetilde{X}$, and there are defined Hurewicz-type morphisms
$$\pi_{n+k}(T(\nu)) \to Q_n(C(\widetilde{X}),\gamma)$$
with $T(\nu)$ the Thom space. An $n$-dimensional normal space
$(X,\nu:X \to BG(k),\rho:S^{n+k} \to T(\nu))$ (Quinn \cite{quinn})
determines an $n$-dimensional algebraic normal complex
$(C(\widetilde{X}),\phi,\gamma,\theta)$ over $\Z[\pi_1(X)]$.
An $n$-dimensional geometric Poincar\'e complex $X$ has a Spivak normal
structure $(\nu,\rho)$ such that the composite of the Hurewicz map and
the Thom isomorphism
$$\pi_{n+k}(T(\nu)) \to \widetilde{H}_{n+k}(T(\nu))~\cong~H_n(X)$$
sends $\rho$ to the fundamental class $[X] \in H_n(X)$, and there is
defined an $n$-dimensional symmetric Poincar\'e complex
$(C(\widetilde{X}),\phi)$ over $\Z[\pi_1(X)]$, with
$$\phi_0~=~[X] \cap -~:~C(\widetilde{X})^{n-*} \to C(\widetilde{X})~.$$
The symmetric signature of $X$ is the symmetric Poincar\'e cobordism class
$$\sigma^*(X)~=~(C(\widetilde{X}),\phi) \in L^n(\Z[\pi_1(X)])$$
which is both a homotopy and $K(\pi_1(X),1)$-bordism invariant. The algebraic normal invariant
of a normal space $(X,\nu,\rho)$
$$[\rho]~=~(\phi,\theta) \in Q_n(C(\widetilde{X}),\gamma)$$
is a homotopy invariant. The classifying chain bundle map
$$(g,\chi)~:~(C(\widetilde{X}),\gamma) \to (B^{\Z[\pi_1(X)]},\beta^{\Z[\pi_1(X)]})$$
sends $[\rho]$ to the hyperquadratic signature of $X$
$$\widehat{\sigma}^*(X)~=~[\phi,\theta] \in Q_n(B^{\Z[\pi_1(X)]},\beta^{\Z[\pi_1(X)]})~=~
\widehat{L}^{n+4*}(\Z[\pi_1(X)])~,$$
which is both a homotopy and $K(\pi_1(X),1)$-bordism invariant.  The
(simply-connected) symmetric signature of a $4k$-dimensional geometric
Poincar\'e complex $X$ is just the signature
$$\sigma^*(X)~=~{\rm signature}(X) \in L^{4k}(\Z)~=~\Z$$
and the hyperquadratic signature is the mod 8 reduction of the signature
$$\widehat{\sigma}^*(X)~=~{\rm signature}(X) \in \widehat{L}^{4k}(\Z)~=~
\Z_8~.$$
See Ranicki \cite{ranicki4} for a more extended discussion of the
connections between chain bundles and their geometric models.]

\section{The $Q$- and $L$-groups}

\subsection{Duality}

Let $T \in \Z_2$ be the generator. The {\it Tate $\Z_2$-cohomology} groups
of a $\Z[\Z_2]$-module $M$ are given by
$$\widehat{H}^n(\Z_2;M)~=~
\frac{\{x \in M \,\vert\, T(x)=(-1)^nx\}}
{\{y+(-1)^nT(y)\,\vert\, y \in M\}}~,$$
and the
$\begin{cases} \hbox{\it $\Z_2$-cohomology}\\
\hbox{\it $\Z_2$-homology} \end{cases}$
groups are given by
$$\begin{array}{l}
H^n(\Z_2;M)~=~
\begin{cases} \{x \in M \,\vert\, T(x)=x\}&\hbox{if $n=0$}\\
\widehat{H}^n(\Z_2;M)&\hbox{if $n>0$}\\
0&\hbox{if $n<0$}~,\end{cases}\\[6ex]
H_n(\Z_2;M)~=~
\begin{cases} M/\{y-T(y)\,\vert\,y \in M \}&\hbox{if $n=0$}\\
\widehat{H}^{n+1}(\Z_2;M)&\hbox{if $n>0$}\\
0&\hbox{if $n<0$}~.\end{cases}
\end{array}$$
\indent We recall some standard properties of $\Z_2$-(co)homology~:

\begin{Proposition}\label{Tate}
Let $M$ be a $\Z[\Z_2]$-module.\\
{\rm (i)} There is defined an exact sequence
$$\dots \to H_n(\Z_2;M) \xymatrix{\ar[r]^N&}
H^{-n}(\Z_2;M) \to \widehat{H}^n(\Z_2;M) \to H_{n-1}(\Z_2;M) \to \dots$$
with
$$N~=~1+T~:~H_0(\Z_2;M) \to H^0(\Z_2;M)~;~x \mapsto x+T(x)~.$$
{\rm (ii)} The Tate $\Z_2$-cohomology groups are 2-periodic and of exponent 2
$$\widehat{H}^*(\Z_2;M)~=~\widehat{H}^{*+2}(\Z_2;M)~,~
2 \widehat{H}^*(\Z_2;M)~=~0~.$$
{\rm (iii)} $\widehat{H}^*(\Z_2;M)=0$ if $M$ is a free $\Z[\Z_2]$-module.
\hfill$\qed$
\end{Proposition}

Let $A$ be an associative ring with $1$, and with an involution
$$\bar{ }~:~A \to A~;~a\mapsto \overline{a}~,$$
such that
$$\overline{a+b}~=~\overline{a}+\overline{b}~,~
\overline{ab}~=~\overline{b}.\overline{a}~,~
\overline{1}~=~1~,~\overline{\overline{a}}~=~a~.$$
When a ring $A$ is declared to be commutative it is given the identity
involution.

\begin{Definition} {\rm For a ring with involution $A$ and $\epsilon=\pm 1$
let $(A,\epsilon)$ denote the $\Z[\Z_2]$-module given by $A$ with $T \in \Z_2$
acting by
$$T_{\epsilon}~:~A \to A~;~a \mapsto \epsilon \overline{a}~.$$
\hfill$\qed$}
\end{Definition}

For $\epsilon=1$ we shall write
$$\widehat{H}^*(\Z_2;A,1)~=~\widehat{H}^*(\Z_2;A)~,~
H^*(\Z_2;A,1)~=~H^*(\Z_2;A)~,~H_*(\Z_2;A,1)~=~H_*(\Z_2;A)~.$$

The {\it dual} of a f.g.  projective (left) $A$-module $P$ is the f.g.
projective $A$-module
$$P^*~=~{\rm Hom}_A(P,A)~,~A \times P^* \to P^*~;~(a,f) \mapsto
(x \mapsto f(x)\overline{a})~.$$
The natural $A$-module isomorphism
$$P \to P^{**}~;~x \mapsto (f \mapsto \overline{f(x)})$$
is used to identify
$$P^{**}~=~P~.$$
For any f.g. projective $A$-modules $P,Q$ there is defined an
isomorphism
$$P\otimes_AQ \to \Hom_A (P^*,Q)~;~x \otimes y \mapsto
(f \mapsto \overline{f(x)}y)$$
regarding $Q$ as a right $A$-module by
$$Q \times A \to Q~;~(y,a) \mapsto \overline{a}y~.$$
There is also a duality isomorphism
$$T~:~{\rm Hom}_A(P,Q) \to {\rm Hom}_A(Q^*,P^*)~;~ f \mapsto f^*$$
with
$$f^*~:~Q^* \to P^*~;~ g \mapsto (x \mapsto g(f(x)))~.$$

\begin{Definition} {\rm For any f.g.  projective $A$-module $P$ and
$\epsilon = \pm 1$ let $(S(P),T_{\epsilon})$ denote the
$\Z[\Z_2]$-module given by the abelian group
$$S(P)~=~\Hom_A(P,P^*)$$
with $\Z_2$-action by the $\epsilon$-duality involution
$$T_{\epsilon}~:~S(P) \to S(P)~;~\phi \mapsto \epsilon \phi^*~.$$
Furthermore, let
$$\begin{array}{l}
{\rm Sym}(P,\epsilon)~=~H^0(\Z_2;S(P),T_{\epsilon})~=~
\{ \phi \in S(P)\,\vert\, T_{\epsilon} (\phi)=\phi\}~,\\[1ex]
{\rm Quad}(P,\epsilon)~=~H_0(\Z_2;S(P),T_{\epsilon})~=~
\dfrac{S(P)}{\{ \theta \in S(P)\,\vert\, \theta -T_{\epsilon} (\theta)\}}~.
\end{array}$$
\hfill$\qed$}
\end{Definition}

An element $\phi \in S(P)$ can be regarded as a sesquilinear form
$$\phi ~:~ P \times P \to A~;~(x,y) \mapsto \langle x,y\rangle_{\phi}~=~\phi(x)(y)$$
such that
$$\langle ax,by\rangle_{\phi}~=~b\langle x,y\rangle_{\phi}\overline{a} \in A~~
(x,y \in P,a,b \in A)~,$$
with
$$\langle x,y\rangle_{T_{\epsilon}(\phi)}~=~\epsilon
\overline{\langle y,x\rangle}_{\phi} \in A~.$$
An $A$-module morphism $f:P \to Q$ induces contravariantly a
$\Z[\Z_2]$-module morphism
$$S(f)~:~(S(Q),T_{\epsilon}) \to (S(P),T_{\epsilon})~;~\theta \mapsto f^* \theta f~.$$
\indent For a f.g. free $A$-module $P=A^r$ we shall use the $A$-module
isomorphism
$$A^r \to (A^r)^*~;~(a_1,a_2,\dots,a_r) \mapsto
((b_1,b_2,\dots,b_r) \mapsto \sum\limits^r_{i=1}b_i\overline{a}_i)$$
to identify
$$(A^r)^*~=~A^r~,~\Hom_A (A^r,(A^r)^*)~=~M_r (A)$$
noting that the duality involution $T$ corresponds to the conjugate
transposition of a matrix. We can thus identify
$$ \begin{array}{l}
M_r(A)~=~S(A^r)~=~\hbox{\rm additive group of $r\times r$ matrices $(a_{ij})$
with $a_{ij} \in A$}~,\\[1ex]
T~:~M_r(A) \to M_r(A)~;~M~=~(a_{ij}) \mapsto M^t~=~(\overline{a}_{ji})~,\\[1ex]
{\rm Sym}_r(A,\epsilon)~=~{\rm Sym}(A^r,\epsilon)~=~\{(a_{ij}) \in
M_r(A)\,\vert\,a_{ij}=\epsilon\overline{a}_{ji}\}~,\\[1ex]
{\rm Quad}_r(A,\epsilon)~=~{\rm Quad}(A^r,\epsilon)~=~
\dfrac{M_r(A)}{\{
(a_{ij}-\epsilon\overline{a}_{ji})\,\vert\,(a_{ij})\in M_r(A)\}}~,\\[2ex]
1+ T_{\epsilon}~:~\qad_r(A,\epsilon) \to \sym_r(A,\epsilon)~;~
M \mapsto M+\epsilon M^t~.
\end{array} $$
The homology of the chain complex
$$\xymatrix{\dots \ar[r] & M_r(A) \ar[r]^{1-T} & M_r(A) \ar[r]^{1+T}
& M_r(A) \ar[r]^{1-T} & M_r(A) \ar[r]& \dots}$$
is given by
$$\frac{{\rm ker}(1-(-1)^nT:M_r(A) \to M_r(A))}
{{\rm im}(1+(-1)^nT:M_r(A) \to M_r(A))}~=~
\widehat{H}^n(\Z_2;M_r(A))~=~
\bigoplus\limits_r \widehat{H}^n(\Z_2;A)~.$$
The $(-1)^n$-symmetrization map $1+(-1)^nT:\sym_r(A) \to \qad_r(A)$
fits into an exact sequence
$$\begin{array}{ll}
0 \to \bigoplus\limits_r \widehat{H}^{n+1}(\Z_2;A) &\to
\qad_r(A,(-1)^n) \\
&\xymatrix{\ar[r]^{1+(-1)^nT}&}
\sym_r(A,(-1)^n) \to \bigoplus\limits_r \widehat{H}^n(\Z_2;A) \to 0~.
\end{array}$$
For $\epsilon=1$ we abbreviate
$$\begin{array}{l}
\sym(P,1)~=~\sym(P)~,~\qad(P,1)~=~\qad(P)~,\\[1ex]
\sym_r(A,1)~=~\sym_r(A)~,~\qad_r(A,1)~=~\qad_r(A)~.
\end{array}$$

\begin{Definition} \label{even}
{\rm An involution on a ring $A$ is {\it even} if
$$\widehat{H}^1(\Z_2;A)~=0~,$$
that is if
$$\{a \in A \,\vert\, a+\overline{a}=0\}~=~\{b-\overline{b}\,\vert\,
b \in A\}~.$$
\hfill$\qed$}
\end{Definition}

\begin{Proposition} \label{even2}
{\rm (i)} For any f.g. projective $A$-module $P$
there is defined an exact sequence
$$0 \to \widehat{H}^1(\Z_2;S(P),T) \to {\rm Quad}(P)
\xymatrix{\ar[r]^{1+T}&} {\rm Sym}(P)~,$$
with
$$1+T~:~{\rm Quad}(P) \to {\rm Sym}(P)~;~\psi \mapsto \psi+\psi^*~.$$
{\rm (ii)} If the involution on $A$ is even the symmetrization
$1+T:{\rm Quad}(P) \to {\rm Sym}(P)$ is injective, and
$$\widehat{H}^n(\Z_2;S(P),T)~=~
\begin{cases}
\dfrac{{\rm Sym}(P)}{{\rm Quad}(P)}&\text{if $n$ is even}\\
0&\text{if $n$ is odd}~,
\end{cases}$$
identifying ${\rm Quad}(P)$ with ${\rm im}(1+T) \subseteq {\rm Sym}(P)$.
\end{Proposition}
\begin{proof} (i) This is a special case of \ref{Tate} (i).\\
(ii) If $Q$ is a f.g. projective $A$-module such that $P\oplus Q=A^r$
is f.g. free then
$$\begin{array}{ll}
\widehat{H}^1(\Z_2;S(P),T) \oplus \widehat{H}^1(\Z_2;S(Q),T)&=~
\widehat{H}^1(\Z_2;S(P\oplus Q),T)\\[1ex]
&=~\bigoplus\limits_r\widehat{H}^1(\Z_2;A,-T)~=~0
\end{array}$$
and so $\widehat{H}^1(\Z_2;S(P),T)=0$.
\end{proof}

In particular, if the involution on $A$ is even
there is defined an exact sequence
$$0 \to {\rm Quad}_r(A) \xymatrix{\ar[r]^{1+T}&}
{\rm Sym}_r(A) \to \bigoplus\limits_r \widehat{H}^0(\Z_2;A) \to 0$$
with
$${\rm Sym}_r(A) \to \bigoplus\limits_r \widehat{H}^0(\Z_2;A)~;~
(a_{ij}) \mapsto (a_{ii})~.$$

For any involution on $A$,
$\sym_r(A)$ is the additive group of symmetric $r\times r$ matrices
$(a_{ij})=(\overline{a}_{ji})$ with $a_{ij} \in A$.
For an even involution
$\qad_r(A) \subseteq \sym_r(A)$ is the subgroup of the matrices
such that the diagonal terms are of the form $a_{ii}=b_i+\overline{b}_i$
for some $b_i \in A$, with
$$\dfrac{\sym_r(A)}{\qad_r(A)}~=~\bigoplus\limits_r \widehat{H}^0(\Z_2;A)~.$$

\begin{Definition} {\rm A ring $A$ is {\it even} if $2\in A$ is
a non-zero divisor, i.e. $2:A \to A$ is injective.\hfill$\qed$}
\end{Definition}

\begin{Example} \label{expl}
{\rm (i) An integral domain $A$ is even if and only if it has characteristic
$\neq 2$.\\
(ii) The identity involution on a commutative ring $A$ is even (\ref{even})
if and only if $A$ is even, in which case
$$\widehat{H}^n(\Z_2;A)~=~
\begin{cases} A_2&\hbox{if $n \equiv 0(\bmod\, 2)$}\\
0&\hbox{if $n \equiv 1(\bmod\, 2)$}
\end{cases}$$
and
$$\qad_r(A)~=~\{(a_{ij}) \in \sym_r(A)\,\vert\, a_{ii} \in 2A\}~.$$
\hfill$\qed$}
\end{Example}

\begin{Example} {\rm For any group $\pi$ there is defined an involution on the
group ring $\Z[\pi]$
$$\overline{~}~:~\Z[\pi] \to \Z[\pi]~;~\sum\limits_{g \in \pi}n_g g \mapsto
\sum\limits_{g \in \pi}n_g g^{-1}~.$$
If $\pi$ has no 2-torsion this involution is even.\hfill$\qed$}
\end{Example}

\subsection{The Hyperquadratic $Q$-Groups}

Let $C$ be a finite (left) f.g. projective $A$-module chain complex.
The dual of the f.g. projective $A$-module $C_p$ is written
$$C^p~=~(C_p)^*~=~\Hom_A(C_p,A)~.$$
The dual $A$-module chain complex $C^{-*}$ is defined by
$$d_{C^{-*}}~=~(d_C)^*~:~(C^{-*})_r~=~C^{-r} \to
(C^{-*})_{r-1}~=~C^{-r+1}~.$$
The $n$-dual $A$-module chain complex $C^{n-*}$ is defined by
$$d_{C^{n-*}}~=~(-1)^r(d_C)^*~:~
(C^{n-*})_r~=~C^{n-r} \to (C^{n-*})_{r-1}~=~C^{n-r+1}~.$$
Identify
$$C\otimes_AC~=~\Hom_A (C^{-\ast},C)~,$$
noting that a cycle $\phi \in (C\otimes_AC)_n$ is a chain map
$\phi:C^{n-*} \to C$. For $\epsilon= \pm 1$ the
$\epsilon$-transposition involution $T_{\epsilon}$ on $C\otimes_AC$
corresponds to the $\epsilon$-duality involution on $\Hom_A (C^{-\ast},C)$
$$T_{\epsilon}~:~\Hom_A (C^p,C_q) \to  \Hom_A (C^q,C_p)~;~
   \phi  \mapsto  (-1)^{pq}\epsilon \phi^*~.$$

Let $\widehat{W}$ be the complete resolution of the $\Z[\Z_2]$-module $\Z$
$$ \widehat{W}~:~\dots \longrightarrow
\widehat{W}_1=
\Z[\Z_2] \stackrel{1-T}{\longrightarrow}
\widehat{W}_0=\Z[\Z_2] \stackrel{1+T}{\longrightarrow}
\widehat{W}_{-1}=\Z[\Z_2] \stackrel{1-T}{\longrightarrow}
\widehat{W}_{-2}=\Z[\Z_2] \longrightarrow       \dots~. $$
If we set
$$\widehat{W}^\% C~=~\Hom_{\Z[\Z_2]} (\widehat{W}, \Hom_A (C^{-\ast},C))~,$$
then an {\it $n$-dimensional $\epsilon$-hyperquadratic structure on $C$} is a cycle
$\theta \in (\widehat{W}^\% C)_n$, which is just a collection
$\{\theta_s \in {\rm Hom}_A(C^r,C_{n-r+s})\,\vert\,r,s \in \Z\}$
such that
$$d\theta_s+(-1)^r\theta_sd^*+(-1)^{n+s-1}(\theta_{s-1}+(-1)^sT_{\epsilon}
\theta_{s-1})~=~0~:~C^r \to C_{n-r+s-1}~.$$

\begin{Definition}
{\rm The {\it $n$-dimensional $\epsilon$-hyperquadratic $Q$-group}
$\widehat{Q}^n (C,\epsilon)$ is the abelian group of equivalence classes of
$n$-dimensional $\epsilon$-hyperquadratic structures on $C$, that is,
$$\widehat{Q}^n (C,\epsilon)~=~H_n (\widehat{W}^\% C)~.$$
\hfill$\qed$}
\end{Definition}

The $\epsilon$-hyperquadratic $Q$-groups are 2-periodic and of exponent 2
$$\widehat{Q}^*(C,\epsilon)~\cong~ \widehat{Q}^{*+2}(C,\epsilon)~,~2\widehat{Q}^*(C,\epsilon)~=~0~.$$
More precisely, there are defined isomorphisms
$$\widehat{Q}^n(C,\epsilon)
\xymatrix{\ar[r]^-{\displaystyle{\cong}}&} \widehat{Q}^{n+2}(C,\epsilon)~;~
\{\theta_s\} \mapsto \{\theta_{s+2}\}~,$$
and for any $n$-dimensional $\epsilon$-hyperquadratic structure $\{\theta_s\}$
$$2\theta_s~=~d\chi_s+(-1)^r\chi_sd^*+(-1)^{n+s}(\chi_{s-1}+(-1)^sT_{\epsilon}\chi_{s-1})~:~
C^r \to C_{n-r+s}$$
with $\chi_s=(-1)^{n+s-1}\theta_{s+1}$. There are also defined
suspension isomorphisms
$$S~:~\widehat{Q}^n(C,\epsilon)
\xymatrix{\ar[r]^-{\displaystyle{\cong}}&}
\widehat{Q}^{n+1}(C_{*-1},\epsilon)~;~\{\theta_s\} \mapsto \{\theta_{s-1}\}$$
and skew-suspension isomorphisms
$$\overline{S}~:~\widehat{Q}^n(C,\epsilon)
\xymatrix{\ar[r]^-{\displaystyle{\cong}}&}
\widehat{Q}^{n+2}(C_{*-1},-\epsilon)~;~\{\theta_s\} \mapsto \{\theta_s\}~.$$

\begin{Proposition}
Let $C$ be  a f.g. projective $A$-module chain complex
which is concentrated in degree $k$
$$C~:~\dots \to 0 \to C_k \to 0 \to \dots~.$$
The $\epsilon$-hyperquadratic $Q$-groups of $C$ are given by
$$\widehat{Q}^n(C,\epsilon)~=~\widehat{H}^{n-2k}(\Z_2;S(C^k),(-1)^kT_{\epsilon})$$
{\rm (}with $S(C^k)={\rm Hom}_A(C^k,C_k)${\rm )}.
\end{Proposition}
\begin{proof} The $\Z[\Z_2]$-module chain complex
$V=\Hom_A(C^{-\ast},C)$ is given by
$$V~:~\dots \to V_{2k+1}~=~0 \to
V_{2k}~=~S(C^k) \to V_{2k-1}~=~0\to \dots$$
and
$$(\widehat{W}^\% C)_j~=~\Hom_{\Z [\Z_2]} (\widehat{W}_{2k-j}, V_{2k})~=~
\Hom_{\Z [\Z_2]} (\widehat{W}_{2k-j}, S(C^k))~.$$
Thus the chain complex $\widehat{W}^\% C$ is of the form
$$\begin{array}{lcl}
(\widehat{W}^\% C)_{2k+1} &=&\Hom_{\Z [\Z_2]} (\widehat{W}_{-1}, V_{2k})~=~S(C^k)\\
\xymatrix{\ar[d]^-{d_{2k+1}=1+(-1)^kT_{\epsilon}}&\\&} & & \\
(\widehat{W}^\% C)_{2k} &=&\Hom_{\Z [\Z_2]} (\widehat{W}_0, V_{2k})~=~S(C^k) \\
\xymatrix{\ar[d]^-{d_{2k}=1+(-1)^{k+1}T_{\epsilon}}&\\&}  & & \\
(\widehat{W}^\% C)_{2k-1} &=&\Hom_{\Z [\Z_2]} (\widehat{W}_1, V_{2k})~=~S(C^k) \\
\xymatrix{\ar[d]^-{d_{2k-1}=1+(-1)^kT_{\epsilon}}&\\&} & & \\
(\widehat{W}^\% C)_{2k-2} &=&\Hom_{\Z [\Z_2]} (\widehat{W}_2, V_{2k})~=~S(C^k) \\
\xymatrix{\ar[d]&\\&}& &
\end{array}$$
and
$$\widehat{Q}^n(C,\epsilon)~=~H_n(\widehat{W}^{\%}C)~=~
\widehat{H}^{n-2k}(\Z_2;S(C^k),(-1)^kT_{\epsilon})~.$$
\end{proof}

\begin{Example} \label{expl.abshyperQ}
{\rm The $\epsilon$-hyperquadratic $Q$-groups of a 0-dimensional f.g. free $A$-module
chain complex
$$C~:~\dots \to 0 \to C_0~=~A^r \to 0 \to \dots$$
are given by
$$\widehat{Q}^n(C,\epsilon)~=~\bigoplus\limits_r \widehat{H}^n(\Z_2;A,\epsilon)~.$$
\hfill$\qed$}
\end{Example}

The {\it algebraic mapping cone} $\Ca(f)$ of a chain map
$f:C\to D$ is the chain complex defined as usual by
$$d_{\Ca(f)}~=~\begin{pmatrix} d_D& (-1)^{r-1}f \\ 0& d_C\end{pmatrix}~ :~
\Ca(f)_r~=~D_r\oplus C_{r-1} \to \Ca(f)_{r-1}~=~D_{r-1}\oplus C_{r-2}~.$$
The relative homology groups
$$H_n(f)~=~H_n(\Ca(f))$$
fit into an exact sequence
$$\dots \to H_n(C)~\raise4pt\hbox{$f_* \atop \to$}~H_n(D) \to
H_n(f) \to H_{n-1}(C) \to \dots~ .$$
\indent An $A$-module chain map $f:C\rightarrow D$ induces a chain map
$$ \widehat{f}^\%~=~\Hom_{\Z[\Z_2]} (1_{\widehat{W}}, \Hom_A (f^*, f))~:~
 \widehat{W}^\% C \longrightarrow \widehat{W}^\% D $$
which induces
$$ \widehat{f}^\% ~:~\widehat{Q}^n (C,\epsilon)
\longrightarrow \widehat{Q}^n (D,\epsilon) $$
on homology. The {\it relative $\epsilon$-hyperquadratic $Q$-group}
$$\widehat{Q}^n (f,\epsilon)~=~H_n(\widehat{f}^{\%}:\widehat{W}^\% C \to \widehat{W}^\% D)$$
fits into a long exact sequence
$$ \xymatrix{\dots \ar[r] & \widehat{Q}^n (C,\epsilon)
\ar[r]^-{\widehat{f}^\%} &    \widehat{Q}^n (D,\epsilon) \ar[r] &
\widehat{Q}^n (f,\epsilon) \ar[r] & \widehat{Q}^{n-1} (C,\epsilon) \ar[r] & \dots~.}$$

\begin{Proposition} \label{hyperThom}
{\rm (i)} The relative $\epsilon$-hyperquadratic $Q$-groups of an
$A$-module chain map $f:C \to D$ are isomorphic to the absolute
$\epsilon$-hyperquadratic $Q$-groups of the algebraic mapping cone
$\Ca(f)$
$$\widehat{Q}^*(f,\epsilon)~\cong~\widehat{Q}^*(\Ca(f),\epsilon)~.$$
{\rm (ii)} If $f:C \to D$ is a chain equivalence the morphisms
$\widehat{f}^{\%}:\widehat{Q}^*(C,\epsilon) \to \widehat{Q}^*(D,\epsilon)$ are isomorphisms,
and
$$\widehat{Q}^*(f,\epsilon)~=~0~.$$
{\rm (iii)} The $\epsilon$-hyperquadratic $Q$-groups are additive:
for any collection $\{C(i)\,\vert\,i \in \Z\}$ of
f.g. projective $A$-module chain complexes $C(i)$
$$\widehat{Q}^n(\sum\limits_iC(i),\epsilon)~=~
\bigoplus\limits_i\widehat{Q}^n(C(i),\epsilon)~.$$
\end{Proposition}
\begin{proof}
(i) See \cite[\S1,\S3]{ranicki1} for the definition of the
$\Z_2$-isovariant chain map
$t: \Ca (f\otimes f)\rightarrow \Ca (f)\otimes_A \Ca (f)$
inducing the algebraic Thom construction
$$t~:~ \widehat{Q}^n (f,\epsilon) \to  \widehat{Q}^n (\Ca (f),\epsilon)~;~
 (\theta,\partial \theta)  \mapsto  \theta / \partial \theta$$
with
$$ \begin{array}{l}
(\theta/\partial \theta)_s~=~
\begin{pmatrix} \theta_s & 0 \\ \partial \theta_s f^* &
T_{\epsilon}\partial\theta_{s-1}
\end{pmatrix}~:\\[3ex]
\Ca (f)^{n-r+s}~=~D^{n-r+s} \oplus C^{n-r+s-1} \to
\Ca (f)_r~=~D_r \oplus C_{r-1}~~(r,s \in \Z)~.
\end{array} $$
Define a free $\Z[\Z_2]$-module chain complex
$$ E~=~(C_{\ast -1} \otimes_A \Ca (f)) \oplus
       (\Ca (f) \otimes_AC_{\ast -1})$$
with
$$T~:~E \to E~;~(a\otimes b, x\otimes y) \mapsto(y\otimes x, b\otimes a)~,$$
such that
$$ H_\ast (\widehat{W} \otimes_{\Z[\Z_2]} E) ~=~ H_\ast (\Hom_{\Z[\Z_2]} (\widehat{W},E)) ~=~0~.$$
Let $p:\Ca (f) \rightarrow C_{\ast -1}$ be the projection.
It now follows from the chain homotopy cofibration
$$ \xymatrix@C+15pt{\Ca (f\otimes f) \ar[r]^-{\displaystyle{t}}&
   \Ca (f) \otimes_A \Ca (f)
   \ar[r]^-{\begin{pmatrix} p \otimes 1 \\ 1 \otimes p \end{pmatrix}}& E}$$
that $t$ induces isomorphisms
$$\widehat{t}~:~\widehat{Q}^*(f,\epsilon)~ \cong~ \widehat{Q}^*(\Ca(f),\epsilon)~.$$
(ii)+(iii) See \cite[Propositions 1.1,1.4]{ranicki1}.
\end{proof}

\begin{Proposition} \label{hypercor}
Let $C$ be a f.g. projective $A$-module chain complex which is
concentrated in degrees $k,k+1$
$$C~:~\dots \to 0 \to C_{k+1} \xymatrix{\ar[r]^d&} C_k \to 0 \to \dots~.$$
{\rm (i)} The $\epsilon$-hyperquadratic $Q$-groups of $C$ are the relative Tate $\Z_2$-cohomology groups in the exact sequence
$$\begin{array}{l}
\dots \to \widehat{H}^{n-2k}(\Z_2;S(C^{k+1}),(-1)^kT_{\epsilon})
\xymatrix{\ar[r]^{\widehat{d}^{\%}}&}
\widehat{H}^{n-2k}(\Z_2;S(C^k),(-1)^kT_{\epsilon}) \\[2ex]
\hskip100pt \to
\widehat{Q}^n(C,\epsilon) \to \widehat{H}^{n-2k-1}(\Z_2;S(C^{k+1}),(-1)^kT_{\epsilon}) \to \dots
\end{array} $$
that is
$$\widehat{Q}^n(C,\epsilon)~=~
\dfrac{\{(\phi,\theta) \in S(C^{k+1})\oplus S(C^k)
\,\vert\, \phi^*=(-1)^{n+k-1}\epsilon\phi, d\phi d^*=\theta+(-1)^{n+k-1}\epsilon\theta^*\}}
{\{(\sigma+(-1)^{n+k-1}\epsilon\sigma^*, d\sigma d^* + \tau +(-1)^{n+k}\epsilon\tau^*)\,\vert\,
(\sigma,\tau) \in S(C^{k+1}) \oplus S(C^k)\}}~,$$
with $(\phi,\theta)$ corresponding to the cycle $\beta \in \widehat{W}^{\%}C)_n$
given by
$$\begin{array}{l}
\beta_{n-2k-2}~=~\theta~:~C^{k+1} \to C_{k+1}~,~
\beta_{n-2k}~=~\phi~:~C^k \to C_k~,\\[1ex]
\beta_{n-2k-1}~=~\begin{cases} d \phi~:~C^{k+1} \to C_k&\\
0~:~C^k \to C_{k+1}~.\end{cases}
\end{array}$$
{\rm (ii)} If the involution on $A$ is even then
$$\widehat{Q}^n(C)~=~
\begin{cases}
{\rm coker}(\widehat{d}^{\%}:
\dfrac{{\rm Sym}(C^{k+1})}{{\rm Quad}(C^{k+1})}
\to \dfrac{{\rm Sym}(C^k)}{{\rm Quad}(C^k)})
&\text{if $n-k$ is even}~,\\[2ex]
{\rm ker}(\widehat{d}^{\%}:
\dfrac{{\rm Sym}(C^{k+1})}{{\rm Quad}(C^{k+1})}
\to \dfrac{{\rm Sym}(C^k)}{{\rm Quad}(C^k)})&\text{if $n-k$ is odd}~.
\end{cases}$$
\end{Proposition}
\begin{proof} (i) Immediate from Proposition \ref{hyperThom}.\\
(ii) Combine (i) and the vanishing  $\widehat{H}^1(\Z_2;
S(P),T)=0$ given by Proposition \ref{even2} (ii).
\end{proof}

For $\epsilon=1$ we write
$$T_{\epsilon}~=~T~,~\widehat{Q}^n(C,\epsilon)~=~\widehat{Q}^n(C)~,~
\hbox{$\epsilon$-hyperquadratic}~=~
\hbox{hyperquadratic}~.$$

\begin{Example} \label{expl.relhyperQ}
{\rm Let $A$ be a ring with an involution which is even (\ref{even}).\\
(i) The hyperquadratic $Q$-groups of a 1-dimensional f.g. free $A$-module
chain complex
$$C~:~\dots \to 0 \to C_1=A^q \xymatrix{\ar[r]^-{d}&} C_0=A^r \to 0 \to \dots$$
are given by
$$\widehat{Q}^n(C)~=~\dfrac{
\{(\phi,\theta) \in M_q(A) \oplus M_r(A)
\,\vert\, \phi^*=(-1)^{n-1}\phi, d\phi d^*=
\theta+(-1)^{n-1}\theta^*\}}
{\{(\sigma+(-1)^{n-1} \sigma^*, d\sigma d^* +
\tau +(-1)^n\tau^*\,\vert\,(\sigma,\tau) \in M_q(A) \oplus M_r(A)\}}~.$$
Example \ref{expl.abshyperQ} and Proposition \ref{hypercor}
give an exact sequence
$$\begin{array}{l}
\widehat{H}^{1}(\Z_2;S(C^1),T)=0\to \widehat{Q}^1(C)\\[2ex]
\hskip50pt \xymatrix{\ar[r]&}
\widehat{H}^0(\Z_2;S(C^1),T)=\bigoplus\limits_q\widehat{H}^0(\Z_2;A)\\[2ex]
\hskip50pt \xymatrix{\ar[r]^{\widehat{d}^{\%}}&}
\widehat{H}^0(\Z_2;S(C^0),T)=\bigoplus\limits_r\widehat{H}^0(\Z_2;A)\\[2ex]
\hskip150pt \xymatrix{\ar[r]&} \widehat{Q}^0(C) \to
\widehat{H}^{-1}(\Z_2;S(C^1),T)=0~.
\end{array}$$
(ii) If $A$ is an even commutative ring and
$$d~=~2~:~C_1~=~A^r \to C_0~=~A^r$$
then $\widehat{d}^{\%}=0$ and there are defined isomorphisms
$$\begin{array}{l}
\widehat{Q}^0(C) \xymatrix@C-10pt{\ar[r]^-{\cong}&} \dfrac{{\rm Sym}_r(A)}{{\rm Quad}_r(A)}~=~
\bigoplus\limits_rA_2~;~(\phi,\theta) \mapsto
\theta~=~(\theta_{ii})_{1 \leqslant i \leqslant r}\\[2ex]
\widehat{Q}^1(C) \xymatrix@C-10pt{\ar[r]^-{\cong}&} \dfrac{{\rm Sym}_r(A)}{{\rm Quad}_r(A)}~=~\bigoplus\limits_rA_2~;~(\phi,\theta) \mapsto
\phi~=~(\phi_{ii})_{1 \leqslant i \leqslant r}~.
\end{array}$$
\hfill$\qed$}
\end{Example}

\subsection{The Symmetric $Q$-Groups}

Let $W$ be the standard free $\Z[\Z_2]$-module resolution of $\Z$
$$ W~:~ \dots \longrightarrow W_3=\Z[\Z_2] \stackrel{1-T}{\longrightarrow}
       W_2=\Z[\Z_2] \stackrel{1+T}{\longrightarrow}
       W_1=\Z[\Z_2] \stackrel{1-T}{\longrightarrow} W_0=\Z[\Z_2] \longrightarrow 0. $$
Given a f.g. projective $A$-module chain complex $C$ we set
$$ W^\% C~=~\Hom_{\Z[\Z_2]} (W, \Hom_A (C^{-\ast},C))~,$$
with $T \in \Z_2$ acting on $C\otimes_AC=\Hom_A (C^{-\ast},C)$
by the $\epsilon$-duality involution $T_{\epsilon}$.
An {\it $n$-dimensional $\epsilon$-symmetric structure on $C$} is a cycle
$\phi \in (W^\% C)_n$, which is just a collection
$\{\phi_s \in {\rm Hom}_A(C^r,C_{n-r+s})\,\vert\,r \in \Z,s\geqslant 0\}$
such that
$$\begin{array}{c}
d\phi_s+(-1)^r\phi_sd^*+
(-1)^{n+s-1}(\phi_{s-1}+(-1)^sT_{\epsilon}\phi_{s-1})~=~0~:~C^r \to C_{n-r+s-1}\\[1ex]
(r \in \Z,s \geqslant 0,\phi_{-1}=0)~.
\end{array}$$

\begin{Definition} {\rm
The {\it $n$-dimensional $\epsilon$-symmetric $Q$-group} $Q^n(C,\epsilon)$
is the abelian group of equivalence classes of
$n$-dimensional $\epsilon$-symmetric structures on $C$, that is,
$$Q^n (C,\epsilon)~=~H_n (W^\% C)~.$$
\hfill$\qed$}
\end{Definition}

Note that there are defined skew-suspension isomorphisms
$$\overline{S}~:~Q^n(C,\epsilon)
\xymatrix{\ar[r]^-{\displaystyle{\cong}}&}
Q^{n+2}(C_{*-1},-\epsilon)~;~\{\phi_s\} \mapsto \{\phi_s\}~.$$

\begin{Proposition} \label{symmQ}
The $\epsilon$-symmetric $Q$-groups of a f.g. projective $A$-module chain complex
concentrated in degree $k$
$$C~:~\dots \to 0 \to C_k \to 0 \to \dots$$
are given by
$$\begin{array}{ll}
Q^n(C,\epsilon)&=~H^{2k-n}(\Z_2;S(C^k),(-1)^kT_{\epsilon})\\[2ex]
&=~\begin{cases}
\widehat{H}^{2k-n}(\Z_2;S(C^k),(-1)^kT_{\epsilon})&\hbox{\it if $n\leqslant 2k-1$}\\
H^0(\Z_2;S(C^k),(-1)^kT_{\epsilon})&\hbox{\it if $n=2k$}\\
0&\hbox{\it if $n\geqslant 2k+1$.}
\end{cases}
\end{array}$$
\end{Proposition}
\begin{proof} The $\Z[\Z_2]$-module chain complex
$V=\Hom_A(C^{-\ast},C)$ is given by
$$V~:~\dots \to V_{2k+1}~=~0 \to
V_{2k}~=~S(C^k) \to V_{2k-1}~=~0\to \dots$$
and
$$(W^\% C)_j~=~\Hom_{\Z[\Z_2]}(W_{2k-j},V_{2k})~=~
\Hom_{\Z[\Z_2]}(W_{2k-j},S(C^k))$$
which vanishes for $j>2k$. Thus the chain complex $W^\% C$ is of the form
$$\begin{array}{lcl}
(W^\% C)_{2k+1} &=& 0 \\
\xymatrix{\ar[d]^-{d_{2k+1}}&\\&} & & \\
(W^\% C)_{2k} &=&\Hom_{\Z [\Z_2]} (W_0, V_{2k})~=~S(C^k) \\
\xymatrix{\ar[d]^-{d_{2k}=1+(-1)^{k+1}T_{\epsilon}}&\\&}  & & \\
(W^\% C)_{2k-1} &=&\Hom_{\Z [\Z_2]} (W_1, V_{2k})~=~S(C^k) \\
\xymatrix{\ar[d]^-{d_{2k-1}=1+(-1)^kT_{\epsilon}}&\\&} & & \\
(W^\% C)_{2k-2} &=&\Hom_{\Z [\Z_2]} (W_2, V_{2k})~=~S(C^k) \\
\xymatrix{\ar[d]&\\&}& &
\end{array}$$
and
$$Q^n(C,\epsilon)~=~H_n(W^{\%}C)~=~H^{2k-n}(\Z_2;S(C^k),(-1)^kT_{\epsilon})~.$$
\end{proof}

For $\epsilon=1$ we write
$$T_{\epsilon}~=~T~,~Q^n(C,\epsilon)~=~Q^n(C)~,~\hbox{$\epsilon$-symmetric}~=~
\hbox{symmetric}~.$$

\begin{Example} \label{expl.abssymmQ}
{\rm The symmetric $Q$-groups of a 0-dimensional f.g. free $A$-module
chain complex
$$C~:~\dots \to 0 \to C_0~=~A^r \to 0 \to \dots$$
are given by
$$Q^n(C)~=~\begin{cases}
\bigoplus\limits_r \widehat{H}^n(\Z_2;A)&\hbox{\rm if $n<0$}\\
\sym_r(A)&\hbox{\rm if $n=0$}\\
0&\hbox{\rm otherwise.}
\end{cases}$$
\hfill$\qed$}
\end{Example}

An $A$-module chain map $f:C\rightarrow D$ induces a chain map
$$\Hom_A (f^*, f)~:~\Hom_A (C^{-\ast},C) \to  \Hom_A (D^{-\ast},D)~;~
\phi  \mapsto  f \phi f^*$$
and thus a chain map
$$ f^\%~=~\Hom_{\Z[\Z_2]} (1_W, \Hom_A (f^*, f))~:~
 W^\% C \longrightarrow W^\% D $$
which induces
$$ f^\% ~:~ Q^n (C,\epsilon) \longrightarrow Q^n (D,\epsilon) $$
on homology. The {\it relative $\epsilon$-symmetric $Q$-group}
$$Q^n (f,\epsilon)~=~H_n(f^{\%}:W^\% C \to W^\% D)$$
fits into a long exact sequence
$$ \xymatrix{\dots \ar[r] & Q^n (C,\epsilon) \ar[r]^-{f^\%} &    Q^n (D,\epsilon) \ar[r] &
Q^n (f,\epsilon) \ar[r] & Q^{n-1} (C,\epsilon) \ar[r] & \dots~.}$$

\begin{Proposition} \label{symmThom}
{\rm (i)} The relative $\epsilon$-symmetric $Q$-groups of an
$A$-module chain map $f:C \to D$ are related to the
absolute $\epsilon$-symmetric $Q$-groups of the algebraic mapping cone
$\Ca(f)$ by a long exact sequence
$$\dots \to H_n(\Ca(f)\otimes_A C) \xymatrix{\ar[r]^-{F}&}
Q^n(f,\epsilon) \xymatrix{\ar[r]^-{t}&} Q^n(\Ca(f),\epsilon) \to
H_{n-1}(\Ca(f)\otimes_AC ) \to \dots$$
with
$$t~:~Q^n (f,\epsilon) \to  Q^n(\Ca(f),\epsilon)~;~
 (\phi,\partial \phi)  \mapsto  \phi / \partial \phi$$
the algebraic Thom construction
$$ \begin{array}{l}
(\phi/\partial \phi)_s~=~
\begin{pmatrix} \phi_s & 0 \\ \partial \phi_s f^* & T_{\epsilon}\partial\phi_{s-1}
\end{pmatrix}~:\\[3ex]
\Ca (f)^{n-r+s}~=~D^{n-r+s} \oplus C^{n-r+s-1} \to
\Ca (f)_r~=~D_r \oplus C_{r-1}~~(r \in \Z,s \geqslant 0, \phi_{-1}=0)~.
\end{array} $$
An element $(g,h) \in H_n(\Ca(f) \otimes_A C)$ is represented by a
chain map $g:C^{n-1-*} \to C$ together with a chain homotopy
$h:fg \simeq 0:C^{n-1-*} \to D$, and
$$F~:~H_n(\Ca(f) \otimes_A C) \to Q^n(f,\epsilon)~;~(g,h) \mapsto (\phi,\partial \phi)$$
with
$$\partial\phi_s~=~\begin{cases} (1+T_{\epsilon})g&\hbox{if $s=0$}\\
0&\hbox{if $s\geqslant 1$}\end{cases}~,~
\phi_s~=~\begin{cases} (1+T_{\epsilon})hf^*&\hbox{if $s=0$}\\
0&\hbox{if $s\geqslant 1$}~.\end{cases}$$
The map
$$Q^n(\Ca(f),\epsilon) \to H_{n-1}(\Ca(f)\otimes_A C)~;~\phi \mapsto p \phi_0$$
is defined using $p={\rm projection}:\Ca(f) \to C_{*-1}$.\\
{\rm (ii)} If $f:C \to D$ is a chain equivalence the morphisms
$f^{\%}:Q^*(C,\epsilon) \to Q^*(D,\epsilon)$ are isomorphisms, and
$$Q^*(\Ca(f),\epsilon)~=~Q^*(f,\epsilon)~=~0~.$$
{\rm (iii)} For any collection $\{C(i)\,\vert\,i \in \Z\}$ of
f.g. projective $A$-module chain complexes $C(i)$
$$Q^n(\sum\limits_iC(i),\epsilon)~=~\bigoplus\limits_iQ^n(C(i),\epsilon) \oplus
  \bigoplus\limits_{i<j} H_n (C(i)\otimes_AC(j))~.$$
\end{Proposition}
\begin{proof} (i) The long exact sequence is induced by the
chain homotopy cofibration of Proposition \ref{hyperThom}
$$ \xymatrix@C+15pt{\Ca (f\otimes f) \ar[r]^-{\displaystyle{t}}&
   \Ca (f) \otimes_A \Ca (f)
   \ar[r]^-{\begin{pmatrix} p \otimes 1 \\ 1 \otimes p \end{pmatrix}}& E}$$
with
$$\begin{array}{l}
E~=~(C_{*-1}\otimes_A\Ca(f)) \oplus (\Ca(f)\otimes_AC_{*-1})~,\\[1ex]
H_*(W^{\%}E)~=~H_*({\rm Hom}_{\Z[\Z_2]}(W,E))~=~
H_{*-1}(C \otimes_A \Ca(f))~.
\end{array}$$
(ii)+(iii) See \cite[Propositions 1.1,1.4]{ranicki1}.
\end{proof}

\begin{Proposition} \label{symmcor}
Let $C$ be a f.g. projective $A$-module chain complex which is
concentrated in degrees $k,k+1$
$$C~:~\dots \to 0 \to C_{k+1} \xymatrix{\ar[r]^d&} C_k \to 0 \to \dots~.$$
The absolute $\epsilon$-symmetric $Q$-groups $Q^*(C,\epsilon)$ and the
relative $\epsilon$-symmetric $Q$-groups $Q^*(d,\epsilon)$ of $d:C_{k+1} \to C_k$
regarded as a morphism of chain complexes concentrated in degree $k$
are given as follows.\\
{\rm (i)} For $n \neq 2k,2k+1,2k+2$
$$Q^n(C,\epsilon)~=~Q^n(d,\epsilon)~=~\begin{cases}\widehat{Q}^n(d,\epsilon)~=~\widehat{Q}^n(C,\epsilon)&
\hbox{if $n\leqslant 2k-1$}\\
0&\hbox{if $n\geqslant 2k+3$}\end{cases}$$
with $\widehat{Q}^n(C,\epsilon)$ as given by Proposition \ref{hypercor}.\\
{\rm (ii)} For $n=2k,2k+1,2k+2$ there are exact sequences
$$\begin{array}{l}
0 \to Q^{2k+1}(d,\epsilon) \xymatrix@C-7pt{\ar[r]&}
Q^{2k}(C_{k+1},\epsilon)=H^0(\Z_2;S(C^{k+1}),(-1)^kT_{\epsilon})\\[1ex]
\hskip25pt \xymatrix@C-7pt{\ar[r]^-{d^{\%}}&}
Q^{2k}(C_k,\epsilon)=H^0(\Z_2;S(C^k),(-1)^kT_{\epsilon}) \xymatrix@C-7pt{\ar[r]&} Q^{2k}(d,\epsilon)
\\[1ex]
\hskip25pt
\xymatrix@C-7pt{\ar[r]&}Q^{2k-1}(C_{k+1},\epsilon)=
\widehat{H}^1(\Z_2;S(C^{k+1}),(-1)^kT_{\epsilon})\\[1ex]
\hskip25pt \xymatrix@C-7pt{\ar[r]^-{d^{\%}}&}
Q^{2k-1}(C_k,\epsilon)=\widehat{H}^1(\Z_2;S(C^k),(-1)^kT_{\epsilon})~, \\[2ex]
Q^{2k+2}(d,\epsilon)=0 \to Q^{2k+2}(C,\epsilon) \to C_{k+1}\otimes_A H_{k+1}(C)
\xymatrix@C-10pt{\ar[r]^-{F}&} Q^{2k+1}(d,\epsilon)\\[1ex]
\hskip25pt  \xymatrix@C-10pt{\ar[r]^-{t}&} Q^{2k+1}(C,\epsilon) \to
C_{k+1}\otimes_A H_k(C) \xymatrix@C-10pt{\ar[r]^-{F}&}
Q^{2k}(d,\epsilon) \xymatrix@C-10pt{\ar[r]^-{t}&} Q^{2k}(C,\epsilon) \to 0~.
\end{array}$$
\end{Proposition}
\begin{proof} The $\Z[\Z_2]$-module chain complex
$V=\Hom_A(C^{-\ast},C)$ is such that
$$V_n~=~\begin{cases}
S(C^k)&\hbox{if $n=2k$}\\
{\rm Hom}_A(C^k,C_{k+1}) \oplus {\rm Hom}_A(C^{k+1},C_k)
&\hbox{if $n=2k+1$}\\
S(C^{k+1}) &\hbox{if $n=2k+2$}\\
0&\hbox{otherwise}
\end{cases}$$
and
$$(W^{\%}C)_n~=~ \sum\limits_{s=0}^{\infty}{\rm Hom}_A(W_s,V_{n+s})~=~
0~{\rm for}~n \geqslant 2k+3~.$$
\end{proof}

\begin{Example} \label{expl.relsymmQ}
{\rm Let $C$ be a 1-dimensional f.g. free $A$-module chain complex
$$C~:~\dots \to 0 \to C_1=A^q \xymatrix{\ar[r]^-{d}&} C_0=A^r \to 0 \to \dots~,$$
so that $C=\Ca(d)$ is the algebraic mapping cone of the
chain map $d:C_1 \rightarrow C_0$ of 0-dimensional complexes, with
$$d^{\%}~ :~ \Hom_A (C^1, C_1)~=~M_q(A) \to \Hom_A (C^0, C_0)~=~M_r(A)~;~
\phi \mapsto d\phi d^*~.$$
Example \ref{expl.abssymmQ} and Proposition \ref{symmcor} give exact sequences
$$\begin{array}{l}
Q^1(C_0)=0 \to Q^1(d) \xymatrix{\ar[r]&}Q^0(C_1)=\sym_q(A)
\xymatrix{\ar[r]^-{d^{\%}}&} Q^0(C_0)=\sym_r(A)\\[1ex]
\hskip10pt
\xymatrix{\ar[r]&} Q^0(d)
\xymatrix@C-10pt{\ar[r]&}  Q^{-1}(C_1)=\bigoplus\limits_q \widehat{H}^1(\Z_2;A)
\xymatrix{\ar[r]^-{d^{\%}}&} Q^{-1}(C_0)=\bigoplus\limits_r \widehat{H}^1(\Z_2;A)\\[2ex]
H_1(C)\otimes_AC_1  \xymatrix@C-10pt{\ar[r]^-{F}&}
Q^1(d) \xymatrix@C-10pt{\ar[r]^-{t}&}  Q^1(C)
\to H_0(C)\otimes_AC_1 \xymatrix@C-10pt{\ar[r]^-{F}&}
Q^0(d) \xymatrix@C-10pt{\ar[r]^-{t}&}Q^0(C) \to 0~.
\end{array}$$
In particular, if $A$ is an even commutative ring and
$$d~=~2~:~C_1~=~A^r \to C_0~=~A^r$$
then $d^{\%}=4$ and
$$\begin{array}{l}
Q^0(d)~=~\dfrac{\sym_r(A)}{4\sym_r(A)}~,~Q^1(d)~=~0~,\\[2ex]
Q^0(C)~=~{\rm coker}(2(1+T):M_r(A) \to\dfrac{\sym_r(A)}{4\sym_r(A)})~=~
\dfrac{\sym_r(A)}{2\qad_r(A)}~,\\[2ex]
Q^1(C)~=~{\rm ker}(2(1+T):\dfrac{M_r(A)}{2M_r(A)}
 \to\dfrac{\sym_r(A)}{4\sym_r(A)})\\[2ex]
 \hphantom{Q^1(C)~}=~\dfrac{\{(a_{ij}) \in M_r(A) \,\vert\, a_{ij}+a_{ji} \in 2A\}}
{2M_r(A)}~=~\dfrac{\sym_r(A)}{2\sym_r(A)}~.
\end{array}$$
\hfill$\qed$}
\end{Example}

We refer to Ranicki \cite{ranicki1} for the one-one correspondence
between highly-connected algebraic Poincar\'e complexes/pairs
and forms, lagrangians and formations.

\subsection{The Quadratic $Q$-Groups}

Given a f.g. projective $A$-module chain complex $C$ we set
$$ W_\% C~=~W \otimes_{\Z[\Z_2]} \Hom_A (C^{-\ast},C)~, $$
with $T \in \Z_2$ acting on $C\otimes_AC=\Hom_A (C^{-\ast},C)$ by
the $\epsilon$-duality involution $T_{\epsilon}$.
An {\it $n$-dimensional $\epsilon$-quadratic structure on $C$} is a cycle
$\psi \in (W_\% C)_n$, a collection
$\{\psi_s \in {\rm Hom}_A(C^r,C_{n-r-s})\,\vert\,r \in \Z,s\geqslant 0\}$
such that
$$d\psi_s+(-1)^r\psi_sd^*+(-1)^{n-s-1}(\psi_{s+1}+
(-1)^{s+1}T_{\epsilon}\psi_{s+1})~=~0~:~C^r \to C_{n-r-s-1}~.$$

\begin{Definition}
{\rm The {\it $n$-dimensional $\epsilon$-quadratic $Q$-group} $Q_n (C,\epsilon)$
is the abelian group of equivalence classes of
$n$-dimensional $\epsilon$-quadratic structures on $C$, that is,
$$Q_n (C,\epsilon)~=~H_n (W_\% C)~.$$
\hfill$\qed$}
\end{Definition}

Note that there are defined skew-suspension isomorphisms
$$\overline{S}~:~Q_n(C,\epsilon)
\xymatrix{\ar[r]^-{\displaystyle{\cong}}&}
Q_{n+2}(C_{*-1},-\epsilon)~;~\{\psi_s\} \mapsto \{\psi_s\}~.$$

\begin{Proposition} \label{quadQ}
The $\epsilon$-quadratic $Q$-groups of a f.g. projective $A$-module chain complex
concentrated in degree $k$
$$C~:~\dots \to 0 \to C_k \to 0 \to \dots$$
are given by
$$\begin{array}{ll}
Q_n(C,\epsilon)&=~H_{n-2k}(\Z_2;S(C^k),(-1)^kT_{\epsilon})\\[1ex]
&=~\begin{cases}
\widehat{H}^{n-2k+1}(\Z_2;S(C^k),(-1)^kT_{\epsilon})&\hbox{\rm if $n\geqslant 2k+1$}\\
H_0(\Z_2;S(C^k),(-1)^kT_{\epsilon})&\hbox{\rm if $n=2k$}\\
0&\hbox{\rm if $n\leqslant 2k-1$~.}
\end{cases}
\end{array}$$
\end{Proposition}
\begin{proof} The $\Z[\Z_2]$-module chain complex
$V=\Hom_A(C^{-\ast},C)$ is given by
$$V~:~\dots \to V_{2k+1}~=~0 \to
V_{2k}~=~\Hom_A (C^k, C_k) \to V_{2k-1}~=~0\to \dots$$
and
$$(W_\% C)_j~=~W_{j-2k}\otimes_{\Z[\Z_2]}V_{2k}~=~
\Hom_{\Z [\Z_2]} (W_{2k-j},S(C^k))$$
which vanishes for $j<2k$. Thus the chain complex $W_\% C$ is of the
form
$$\begin{array}{lcl}
(W_\% C)_{2k+2} &=&W_2 \otimes_{\Z [\Z_2]}V_{2k}~=~S(C^k) \\
\xymatrix{\ar[d]^-{d_{2k+2}=1+(-1)^kT_{\epsilon}}&\\&} & & \\
(W_\% C)_{2k+1} &=&W_1\otimes_{\Z [\Z_2]}V_{2k}~=~S(C^k) \\
\xymatrix{\ar[d]^-{d_{2k+1}=1+(-1)^{k+1}T_{\epsilon}}&\\&}  & & \\
(W_\% C)_{2k} &=&W_0\otimes_{\Z [\Z_2]}V_{2k}~=~S(C^k) \\
\xymatrix{\ar[d]&\\&} & & \\
(W_\% C)_{2k-1} &=&0 \\
\end{array}$$
and
$$Q_n(C,\epsilon)~=~H_n(W_{\%}C)~=~H_{n-2k}(\Z_2;S(C^k),(-1)^kT_{\epsilon})~.$$
\end{proof}

\begin{Example} \label{expl.absquadQ}
{\rm The $\epsilon$-quadratic $Q$-groups of the 0-dimensional f.g. free $A$-module
chain complex
$$C~:~\dots \to 0 \to C_0~=~A^r \to 0 \to \dots$$
are given by
$$Q_n(C)~=~\begin{cases}
\bigoplus\limits_r \widehat{H}^{n+1}(\Z_2;A)&\hbox{\rm if $n>0$}\\
\qad_r(A)&\hbox{\rm if $n=0$}\\
0&\hbox{\rm otherwise.}
\end{cases}$$
\hfill$\qed$}
\end{Example}

An $A$-module chain map $f:C\rightarrow D$ induces a chain map
$$ f_\%~=~1_W\otimes_{\Z[\Z_2]}\Hom_A (f^*, f)~:~
 W_\% C \longrightarrow W_\% D $$
which induces
$$ f_\% ~:~ Q_n (C,\epsilon) \longrightarrow Q_n (D,\epsilon) $$
on homology. The {\it relative $\epsilon$-quadratic $Q$-group} $Q_n(f,\epsilon)$ is
designed to fit into a long exact sequence
$$ \dots \longrightarrow Q_n (C,\epsilon) \stackrel{f_\%}{\longrightarrow}
   Q_n (D,\epsilon) \longrightarrow Q_n (f,\epsilon) \longrightarrow Q_{n-1} (C,\epsilon)
   \longrightarrow \dots~, $$
that is, $Q_n (f,\epsilon)$ is defined as the $n$-th homology group of the
mapping cone of $f_\%$,
$$ Q_n (f,\epsilon)~=~H_n (f_\%:W_\% C \longrightarrow W_\% D)~. $$

\begin{Proposition} \label{quadThom}
{\rm (i)} The relative $\epsilon$-quadratic $Q$-groups of $f:C \to D$ are related to the
absolute $\epsilon$-quadratic $Q$-groups of the algebraic mapping cone
$\Ca(f)$ by a long exact sequence
$$\dots \to H_n(\Ca(f)\otimes_AC) \xymatrix{\ar[r]^-{F}&}
Q_n(f,\epsilon) \xymatrix{\ar[r]^-{t}&} Q_n(\Ca(f),\epsilon) \to
H_{n-1}(\Ca(f)\otimes_AC) \to \dots~.$$
{\rm (ii)} If $f:C \to D$ is a chain equivalence the morphisms
$f_{\%}:Q_*(C) \to Q_*(D)$ are isomorphisms, and
$$Q_*(\Ca(f),\epsilon)~=~Q_*(f,\epsilon)~=~0~.$$
{\rm (iii)} For any collection $\{C(i)\,\vert\,i \in \Z\}$ of
f.g. projective $A$-module chain complexes $C(i)$
$$Q_n(\sum\limits_iC(i),\epsilon)~=~\bigoplus\limits_iQ_n(C(i),\epsilon) \oplus
  \bigoplus\limits_{i<j} H_n (C(i)\otimes_AC(j))~.$$
\hfill$\qed$
\end{Proposition}

\begin{Proposition} \label{quadcor}
Let $C$ be a f.g. projective $A$-module chain complex which is
concentrated in degrees $k,k+1$
$$C~:~\dots \to 0 \to C_{k+1} \xymatrix{\ar[r]^-{d}&} C_k \to 0 \to \dots~.$$
The absolute $\epsilon$-quadratic $Q$-groups $Q_*(C,\epsilon)$ and the
relative $\epsilon$-quadratic $Q$-groups $Q_*(d,\epsilon)$ of $d:C_{k+1} \to C_k$
regarded as a morphism of chain complexes concentrated in degree $k$
are given as follows.\\
{\rm (i)} For $n \neq 2k,2k+1,2k+2$
$$Q_n(C,\epsilon)~=~Q_n(d,\epsilon)~=~\begin{cases}\widehat{Q}^{n+1}(d,\epsilon)~=~\widehat{Q}^{n+1}(C,\epsilon)&
\hbox{if $n\geqslant 2k+3$}\\
0&\hbox{if $n\leqslant 2k-1$}\end{cases}$$
with $\widehat{Q}^n(C,\epsilon)$ as given by Proposition \ref{hypercor}.\\
{\rm (ii)} For $n=2k,2k+1,2k+2$ there are exact sequences
$$\begin{array}{l}
Q_{2k+2}(C_{k+1},\epsilon)=\widehat{H}^1(\Z_2;S(C^{k+1}),(-1)^kT_{\epsilon})
\xymatrix{\ar[r]^-{d_{\%}}&}
Q_{2k+2}(C_k,\epsilon)=\widehat{H}^1(\Z_2;S(C^k),(-1)^kT_{\epsilon})\\[2ex]
\hskip20pt \xymatrix{\ar[r]&} Q_{2k+2}(d,\epsilon)=\widehat{Q}^{2k+3}(C,\epsilon)
\xymatrix{\ar[r]^-{d_{\%}}&}
Q_{2k+1}(C_{k+1},\epsilon)=\widehat{H}^0(\Z_2;S(C^{k+1}),(-1)^kT_{\epsilon})\\[2ex]
\hskip20pt \xymatrix{\ar[r]^-{d_{\%}}&}
Q_{2k+1}(C_k,\epsilon)=\widehat{H}^0(\Z_2;S(C^k),(-1)^kT_{\epsilon})
\xymatrix{\ar[r]&} Q_{2k+1}(d,\epsilon)\\[2ex]
\hskip20pt \xymatrix{\ar[r]&} Q_{2k}(C_{k+1},\epsilon)=H_0(\Z_2;S(C^{k+1}),(-1)^kT_{\epsilon})\\[2ex]
\hskip20pt
\xymatrix{\ar[r]^-{d_{\%}}&} Q_{2k}(C_k,\epsilon)=H_0(\Z_2;S(C^k),(-1)^kT_{\epsilon})
\xymatrix{\ar[r]&}Q_{2k}(d,\epsilon) \xymatrix@C-10pt{\ar[r]&} Q_{2k-1}(C_{k+1})=0~,\\[2ex]
0 \to Q_{2k+2}(d,\epsilon)\xymatrix{\ar[r]^-{t}&}Q_{2k+2}(C,\epsilon)
\xymatrix{\ar[r]&}
H_{k+1}(C)\otimes_AC_{k+1} \xymatrix{\ar[r]^-{F}&} Q_{2k+1}(d,\epsilon)\\[2ex]
\hphantom{0 \to Q_{2k+2}(d,\epsilon)} \xymatrix{\ar[r]^-{t}&}Q_{2k+1}(C,\epsilon)
\xymatrix@C-10pt{\ar[r]&} C_{k+1}\otimes_A
H_k(C)\xymatrix@C-10pt{\ar[r]^-{F}&} Q_{2k}(d,\epsilon)
\xymatrix{\ar[r]^-{t}&}Q_{2k}(C,\epsilon) \to 0~.
\end{array}$$
\hfill$\qed$
\end{Proposition}

For $\epsilon=1$ we write
$$T_{\epsilon}~=~T~,~Q_n(C,\epsilon)~=~Q_n(C)~,~
\hbox{$\epsilon$-quadratic}~=~\hbox{quadratic}~.$$

\begin{Example} \label{expl.relquadQ}
{\rm Let $C$ be a 1-dimensional f.g. free $A$-module chain complex
$$C~:~\dots \to 0 \to C_1=A^q \xymatrix{\ar[r]^-{d}&} C_0=A^r \to 0 \to \dots~,$$
so that $C=\Ca(d)$ is the algebraic mapping cone of the
chain map $d:C_1 \rightarrow C_0$ of 0-dimensional complexes, with
$$d^{\%}~ :~ \Hom_A (C^1, C_1)~=~M_q(A) \to \Hom_A (C^0, C_0)~=~M_r(A)~;~
\phi \mapsto d\phi d^*~.$$
Example \ref{expl.absquadQ} and Proposition \ref{quadcor} give exact sequences
$$\begin{array}{l}
Q_1(C_1)=\bigoplus\limits_q \widehat{H}^0(\Z_2;A)
\xymatrix{\ar[r]^-{\widehat{d}^{\%}}&}
Q_1(C_0)=\bigoplus\limits_r \widehat{H}^0(\Z_2;A)
\to Q_1(d)\\[1ex]
\xymatrix@C-10pt{\ar[r]&} Q_0(C_1)=\qad_q(A)
\xymatrix{\ar[r]^-{d_{\%}}&} Q_0(C_0)=\qad_r(A)
\xymatrix@C-10pt{\ar[r]&} Q_0(d) \to Q_{-1}(C_1)=0~,\\[2ex]
H_1(C)\otimes_AC_1 \to Q_1(d) \to Q_1(C) \to H_0(C)\otimes_AC_1
\to Q_0(d) \to Q_0(C) \to 0~.
\end{array}$$
In particular, if $A$ is an even commutative ring  and
$$d~=~2~:~C_1~=~A^r \to C_0~=~A^r$$
then $d_{\%}=4$ and
$$\begin{array}{l}
Q_0(d)~=~\dfrac{\qad_r(A)}{4\qad_r(A)}~,\\[3ex]
Q_1(d)~=~\dfrac{\sym_r(A)}{\qad_r(A)+4\sym_r(A)}~,\\[3ex]
Q_0(C)~=~{\rm coker}(2(1+T):\dfrac{M_r(A)}{2M_r(A)} \to
\dfrac{\qad_r(A)}{4\qad_r(A)})~=~\dfrac{\qad_r(A)}{2\qad_r(A)}~,\\[3ex]
Q_1(C)~=~
\dfrac{
\{(\psi_0,\psi_1) \in M_r(A) \oplus M_r(A) \,\vert\, 2\psi_0=\psi_1-\psi_1^*\}}
{\{(2(\chi_0-\chi^*_0),4\chi_0+\chi_2+\chi^*_2)\,\vert\,
(\chi_0,\chi_2) \in M_r(A) \oplus M_r(A)\}}~=~
\bigoplus\limits_{\frac{r(r+1)}{2}} A_2~.
\end{array}$$
\hfill$\qed$}
\end{Example}

\subsection{$L$-groups}

An {\it $n$-dimensional} $\begin{cases} \hbox{\it $\epsilon$-symmetric}\cr
\hbox{\it $\epsilon$-quadratic}\end{cases}$
{\it Poincar\'e complex}
$\begin{cases}
(C,\phi)\\
(C,\psi)
\end{cases}$ over $A$
is an $n$-dimensional f.g. projective $A$-module chain complex
$$C~:~\dots \to 0 \to C_n \to C_{n-1} \to \dots \to C_1 \to C_0 \to 0 \to \dots$$
together with an element $\begin{cases}
\phi \in Q^n(C,\epsilon)\\
\psi \in Q_n(C,\epsilon)
\end{cases}$ such that the $A$-module chain map
$$\begin{cases}
\phi_0:C^{n-*} \to C\\
(1+T_{\epsilon})\psi_0:C^{n-*} \to C
\end{cases}$$
is a chain equivalence.
We refer to \cite{ranicki4} for the detailed definition of the
{\it $n$-dimensional} $\begin{cases} \hbox{\it $\epsilon$-symmetric}\cr
\hbox{\it $\epsilon$-quadratic}
\end{cases}$  {\it $L$-group}
$\begin{cases} L^n(A,\epsilon) \\
L_n(A,\epsilon) \end{cases}$
as the cobordism group of $n$-dimensional\break
$\begin{cases} \hbox{\rm $\epsilon$-symmetric}\\
\hbox{\rm $\epsilon$-quadratic}\\
\end{cases}$ Poincar\'e complexes over $A$.

\begin{Definition}
{\rm
(i) The {\it relative {\rm (}$\epsilon$-symmetric, $\epsilon$-quadratic{\rm )}
$Q$-group} $Q^n_n(f,\epsilon)$ of a chain map $f:C \to D$ of f.g. projective
$A$-module chain complexes is the relative group in the exact sequence
$$\dots \to Q_n(C,\epsilon) \xymatrix{\ar[r]^-{(1+T_{\epsilon})f_{\%}}&} Q^n(D,\epsilon) \to
Q^n_n(f,\epsilon) \to Q_{n-1}(C,\epsilon) \to \dots~.$$
An element $(\delta\phi,\psi) \in Q^n_n(f,\epsilon)$ is an equivalence class
of pairs
$$(\delta\phi,\psi) \in (W^{\%}D)_n \oplus (W_{\%}C)_{n-1}$$
such that
$$d(\psi)~=~0 \in (W_{\%}C)_{n-2}~,~
(1+T_{\epsilon})f_{\%}\psi~=~d(\delta\phi) \in (W^{\%}D)_{n-1}~.$$
(ii) An {\it $n$-dimensional {\rm (}$\epsilon$-symmetric,
$\epsilon$-quadratic{\rm )} pair over $A$}
$(f:C \to D,(\delta\phi,\psi))$ is a chain map $f$ together with
a class $(\delta\phi,\psi) \in Q^n_n(f,\epsilon)$ such that the chain map
$$(\delta\phi,(1+T_{\epsilon})\psi)_0~:~D^{n-*} \to \Ca(f)$$
defined by
$$(\delta\phi,(1+T_{\epsilon})\psi)_0~=~\begin{pmatrix}
\delta \phi_0 \\ (1+T_{\epsilon})\psi_0f^* \end{pmatrix}~:~
D^{n-r}\to \Ca(f)_r~=~D_r \oplus C_{r-1}$$
is a chain equivalence.
\hfill$\qed$}
\end{Definition}

\begin{Proposition}
The relative {\rm (}$\epsilon$-symmetric, $\epsilon$-quadratic{\rm )} $Q$-groups
$Q^n_n(f,\epsilon)$ of a chain map $f:C \to D$
fit into a commutative braid of exact sequences
\vskip2mm

$$\xymatrix@C-10pt{
Q_n(C,\epsilon)\ar[dr]^-{f_{\%}} \ar@/^2pc/[rr]^-{(1+T_{\epsilon})f_{\%}} &&
Q^n(D,\epsilon) \ar[dr] \ar@/^2pc/[rr]^-{J}  &&\widehat{Q}^n(D,\epsilon)  \\&
Q_n(D,\epsilon)\ar[ur]^-{1+T_{\epsilon}} \ar[dr] && Q^n_n(f,\epsilon) \ar[ur]^-{J_f} \ar[dr]&&\\
\widehat{Q}^{n+1}(D,\epsilon)  \ar[ur]^-{H} \ar@/_2pc/[rr]_-{}&&Q_n(f,\epsilon)
\ar[ur]\ar@/_2pc/[rr]_{}&&Q_{n-1}(C,\epsilon)}$$

\vskip2mm

\noindent with
$$\begin{array}{l}
J_f~:~Q^n_n(f,\epsilon) \to \widehat{Q}^n(D,\epsilon)~;~(\delta\phi,\psi) \mapsto \alpha~,\\[1ex]
\alpha_s~=~\begin{cases} \delta\phi_s&\text{if $s\geqslant 0$}\\
f\psi_{-s-1}f^*&\text{if $s \leqslant -1$}\end{cases}~:~D^r \to D_{n-r+s}~.
\end{array}$$
\hfill$\qed$
\end{Proposition}

The {\it $n$-dimensional $\epsilon$-hyperquadratic $L$-group}
$\widehat{L}^n(A,\epsilon)$ is the cobordism group of
$n$-dimensional ($\epsilon$-symmetric, $\epsilon$-quadratic)
Poincar\'e pairs $(f:C \to D,(\phi,\psi))$ over $A$. As in
\cite{ranicki1}, there is defined an exact sequence
$$\xymatrix@C-5pt{ \dots \ar[r] & L_n(A,\epsilon) \ar[r]^{1+T_{\epsilon}} & L^n(A,\epsilon)
    \ar[r] & \widehat{L}^n(A,\epsilon)\ar[r] &  L_{n-1}(A,\epsilon) \ar[r] & \dots~.}$$
The skew-suspension maps in the $\pm \epsilon$-quadratic $L$-groups
are isomorphisms
$$\overline{S}~:~L_n(A,\epsilon)
\xymatrix{\ar[r]^-{\displaystyle{\cong}}&}
L_{n+2}(A,-\epsilon)~;~(C,\{\psi_s\}) \mapsto (C_{*-1},\{\psi_s\})~,$$
so the $\epsilon$-quadratic $L$-groups are 4-periodic
$$L_n(A,\epsilon)~=~L_{n+2}(A,-\epsilon)~=~L_{n+4}(A,\epsilon)~.$$
\indent The skew-suspension maps in $\epsilon$-symmetric and
$\epsilon$-hyperquadratic $L$-groups
and $\pm\epsilon$-hyperquadratic $L$-groups
$$\begin{array}{l}
\overline{S}~:~L^n(A,\epsilon) \to L^{n+2}(A,-\epsilon)~;~(C,\{\phi_s\})
\mapsto (C_{*-1},\{\phi_s\})~,\\[1ex]
\overline{S}~:~\widehat{L}^n(A,\epsilon) \to \widehat{L}^{n+2}(A,-\epsilon)~;~
(f:C\to D,\{\psi_s,\phi_s\})
\mapsto (f:C_{*-1}\to D_{*-1},\{(\psi_s,\phi_s)\})
\end{array}$$
are not isomorphisms in general, so the $\epsilon$-symmetric and
$\epsilon$-hyperquadratic $L$-groups need not be 4-periodic.
We shall write the 4-periodic versions of the $\epsilon$-symmetric
and $\epsilon$-hyperquadratic $L$-groups of $A$ as
$$L^{n+4*}(A,\epsilon)~=~\lim\limits_{k\to\infty}L^{n+4k}(A,\epsilon)~,~
\widehat{L}^{n+4*}(A,\epsilon)~=~\lim\limits_{k\to\infty}\widehat{L}^{n+4k}
(A,\epsilon)~,$$
noting that there is defined an exact sequence
$$\dots \to L_n(A,\epsilon) \to L^{n+4*}(A,\epsilon) \to
\widehat{L}^{n+4*}(A,\epsilon) \to L_{n-1}(A,\epsilon) \to \dots~.$$

\begin{Definition} \label{Wu} {\rm
The {\it Wu classes} of an $n$-dimensional $\epsilon$-symmetric
complex $(C,\phi)$ over $A$ are the $A$-module morphisms
$$\widehat{v}_k(\phi)~:~H^{n-k}(C) \to \widehat{H}^k(\Z_2;A,\epsilon)~;~
x \mapsto \phi_{n-2k}(x)(x)~~(k \in \Z)~.$$
\hfill$\qed$}
\end{Definition}

For an $n$-dimensional $\epsilon$-symmetric Poincar\'e complex $(C,\phi)$
over $A$ the evaluation of the Wu class $\widehat{v}_k(\phi)(x) \in
\widehat{H}^k(\Z_2;A,\epsilon)$ is the obstruction to killing
$x \in H^{n-k}(C) \cong H_k(C)$ by algebraic surgery (\cite[\S4]{ranicki1}).

\begin{Proposition} \label{4period}
{\rm (i)} If $\widehat{H}^0(\Z_2;A,\epsilon)$
has a 1-dimensional f.g. projective $A$-module resolution
then the skew-suspension maps
$$\overline{S}~:~L^{n-2}(A,-\epsilon) \to L^n(A,\epsilon)~,~
\overline{S}~:~\widehat{L}^{n-2}(A,-\epsilon) \to \widehat{L}^n(A,\epsilon)
~~(n \geqslant 2)$$
are isomorphisms. Thus if $\widehat{H}^1(\Z_2;A,\epsilon)$ also
has a 1-dimensional f.g. projective $A$-module resolution
the $\epsilon$-symmetric and $\epsilon$-hyperquadratic $L$-groups of $A$ are
4-periodic
$$\begin{array}{l}
L^n(A,\epsilon)~=~L^{n+2}(A,-\epsilon)~=~L^{n+4}(A,\epsilon)~,\\[1ex]
\widehat{L}^n(A,\epsilon)~=~\widehat{L}^{n+2}(A,-\epsilon)~=~
\widehat{L}^{n+4}(A,\epsilon)~.
\end{array}$$
{\rm (ii)} If $A$ is a Dedekind ring then the $\epsilon$-symmetric 
$L$-groups are `homotopy invariant'
$$L^n(A[x],\epsilon)~=~L^n(A,\epsilon)$$
and the $\epsilon$-symmetric and $\epsilon$-hyperquadratic $L$-groups 
of $A$ and $A[x]$ are 4-periodic.
\end{Proposition}
\begin{proof} (i) Let $D$ be a 1-dimensional
f.g. projective $A$-module resolution of $\widehat{H}^0(\Z_2;A,\epsilon)$
$$0 \to D_1 \to D_0 \to \widehat{H}^0(\Z_2;A) \to 0~.$$
Given an $n$-dimensional $\epsilon$-symmetric Poincar\'e
complex $(C,\phi)$ over $A$ resolve the $A$-module morphism
$$\widehat{v}_n(\phi)(\phi_0)^{-1}~:~H_0(C) \cong H^n(C) \to
H_0(D)~=~\widehat{H}^0(\Z_2;A,\epsilon)~;~u \mapsto (\phi_0)^{-1}(u)(u)$$
by an $A$-module chain map $f:C \to D$,
defining an $(n+1)$-dimensional
$\epsilon$-symmetric pair $(f:C \to D,(\delta\phi,\phi))$.
The effect of algebraic surgery on $(C,\phi)$ using
$(f:C \to D,(\delta\phi,\phi))$ is a cobordant $n$-dimensional $\epsilon$-symmetric
Poincar\'e complex $(C',\phi')$ such that there are defined an exact sequence
$$0 \to H^n(C') \to H^n(C) \xymatrix{\ar[r]^-{\widehat{v}_n(\phi)}&}
\widehat{H}^0(\Z_2;A,\epsilon) \to H^{n+1}(C') \to 0$$
and an $(n+1)$-dimensional
$\epsilon$-symmetric pair $(f':C' \to D',(\delta\phi',\phi'))$ with
$f'$ the projection onto the quotient complex of $C'$ defined by
$$D'~:~\dots \to 0 \to D'_{n+1}~=~C'_{n+1} \to D'_n~=~C'_n \to 0 \to \dots ~.$$
The effect of algebraic surgery on $(C',\phi')$ using
$(f':C' \to D',(\delta\phi',\phi'))$ is a cobordant $n$-dimensional
$\epsilon$-symmetric Poincar\'e complex $(C'',\phi'')$ with $H_n(C'')=0$,
so that it is (homotopy equivalent to) the skew-suspension of an
$(n-2)$-dimensional $(-\epsilon)$-symmetric Poincar\'e complex.\\
(ii) The 4-periodicity $L^*(A,\epsilon)=L^{*+4}(A,\epsilon)$ 
was proved in \cite[\S7]{ranicki1}.
The `homotopy invariance' $L^*(A[x],\epsilon)=L^*(A,\epsilon)$ 
was proved in \cite[41.3]{ranicki3} and \cite[2.1]{connran}.
The 4-periodicity of the $\epsilon$-symmetric and $\epsilon$-hyperquadratic
$L$-groups for $A$ and $A[x]$ now follows from the 4-periodicity
of the $\epsilon$-quadratic $L$-groups $L_*(A,\epsilon)=L_{*+4}(A,\epsilon)$.
\end{proof}

\section{Chain Bundle Theory}

\subsection{Chain Bundles}

\begin{Definition}
{\rm (i) An {\it $\epsilon$-bundle} over an $A$-module chain complex $C$
is a $0$-dimensional $\epsilon$-hyperquadratic structure $\gamma$ on
$C^{0-\ast}$, that is, a cycle
$$ \gamma \in (\widehat{W}^\% C^{0-\ast})_0$$
as given by a collection of $A$-module morphisms
$$\{\gamma_s\in \Hom_A(C_{r-s},C^{-r})\,\vert\, r,s\in \Z\}$$
such that
$$(-1)^{r+1}d^*\gamma_s +(-1)^s\gamma_sd +
(-1)^{s-1}(\gamma_{s-1}+(-1)^sT_{\epsilon}\gamma_{s-1})~=~0~:~C_{r-s+1}\to C^{-r}~.$$
(ii) An {\it equivalence of $\epsilon$-bundles} over $C$,
$$ \chi~:~ \gamma \longrightarrow \gamma' $$
is an equivalence of $\epsilon$-hyperquadratic structures. \\
(iii) A {\it chain $\epsilon$-bundle} $(C,\gamma)$ over
$A$ is an $A$-module chain complex $C$
together with an $\epsilon$-bundle $\gamma \in (\widehat{W}^\% C^{0-\ast})_0$.
$$\eqno{\qed}$$}
\end{Definition}

Let $(D,\delta)$ be a chain $\epsilon$-bundle and $f:C\to D$ a chain map.
The dual of $f$
$$ f^*~:~ D^{0-\ast} \longrightarrow C^{0-\ast}$$
induces a map
$$ (\widehat{f^*})^\%_0~:~ (\widehat{W}^\% D^{0-\ast})_0 \longrightarrow
    (\widehat{W}^\% C^{0-\ast})_0~.$$

\begin{Definition}
{\rm (i) The {\it pullback chain $\epsilon$-bundle} $(C, f^* \delta)$ is defined to be
$$  f^* \delta~=~(\widehat{f^*})^\%_0 (\delta) \in
    (\widehat{W}^\% C^{0-\ast})_0. $$
(ii) A {\it map of chain $\epsilon$-bundles}
$$ (f,\chi)~:~ (C,\gamma) \longrightarrow (D,\delta) $$
is a chain map $f:C\to D$ together with an equivalence of $\epsilon$-bundles
over $C$
$$ \chi~:~ \gamma \longrightarrow f^* \delta~. \eqno{\qed}$$}
\end{Definition}

The $\epsilon$-hyperquadratic $Q$-group $\widehat{Q}^0(C^{0-*},\epsilon)$ is thus the
group of equivalence classes of chain $\epsilon$-bundles on the chain complex $C$,
the algebraic analogue of the topological $K$-group of a space.
The Tate $\Z_2$-cohomology groups
$$\widehat{H}^n(\Z_2;A,\epsilon)~=~
\frac{\{a \in A \,\vert\, \overline{a}=(-1)^n\epsilon a\}}
{\{b+(-1)^n\epsilon\overline{b}\,\vert\, b \in A\}}$$
are $A$-modules via
$$ A\times \widehat{H}^n (\Z_2; A,\epsilon) \rightarrow
\widehat{H}^n (\Z_2; A,\epsilon)~;~   (a,x)\mapsto ax \overline{a}. $$

\begin{Definition}{\rm
The {\it Wu classes} of a chain $\epsilon$-bundle $(C,\gamma)$ are the
$A$-module morphisms
$$ \widehat{v}_k (\gamma)~:~ H_k(C) \to \widehat{H}^k (\Z_2;A,\epsilon)~;~
   x\mapsto \gamma_{-2k}(x)(x)~~(k \in \Z)~. $$
\hfill$\qed$}
\end{Definition}

An $n$-dimensional $\epsilon$-symmetric Poincar\'e complex $(C,\phi)$
with Wu classes (\ref{Wu})
$$\widehat{v}_k(\phi)~:~H^{n-k}(C) \to \widehat{H}^k (\Z_2;A,\epsilon)~;~
y \mapsto  \phi_{n-2k}(y)(y)~~(k \in \Z)$$
has a Spivak normal $\epsilon$-bundle (\cite{ranicki1})
$$\gamma~=~S^{-n}(\phi_0^{\%})^{-1}(J(\phi)) \in \widehat{Q}^0(C^{0-*},\epsilon)$$
such that
$$\widehat{v}_k(\phi)~=~\widehat{v}_k (\gamma)\phi_0~:~
H^{n-k}(C)~\cong~H_k(C)  \to \widehat{H}^k (\Z_2;A,\epsilon)~~(k \in \Z)~,$$
the abstract analogue of the formulae of Wu and Thom.

For any $A$-module chain map $f:C \to D$ Proposition \ref{hyperThom} (i)
gives an exact sequence
$$\dots \to \widehat{Q}^1(C^{0-*},\epsilon) \to
\widehat{Q}^0(\Ca(f)^{0-*},\epsilon) \to \widehat{Q}^0(D^{0-*},\epsilon)
\xymatrix{\ar[r]^-{(\widehat{f}^*)^{\%}}&} \widehat{Q}^0(C^{0-*},\epsilon)
\to \dots~,$$
motivating the following construction of chain $\epsilon$-bundles~:

\begin{Definition} \label{cone}
{\rm The {\it cone} of a chain $\epsilon$-bundle map $(f,\chi):(C,0) \to (D,\delta)$
is the chain $\epsilon$-bundle
$$(B,\beta)~=~\Ca(f,\chi)$$
with $B=\Ca(f)$ the algebraic mapping cone of $f:C\to D$ and
$$\beta_s~=~\begin{pmatrix} \delta_s & 0 \\
f^*\delta_{s+1} & \chi_{s+1} \end{pmatrix}~:~
B_{r-s}~=~D_{r-s} \oplus C_{r-s-1} \to B^{-r}~=~D^{-r} \oplus C^{-r-1}~.$$
Note that $(D,\delta)=g^*(B,\beta)$
is the pullback of $(B,\beta)$ along the inclusion $g:D \to B$.}
$$\eqno{\qed}$$
\end{Definition}

\begin{Proposition} \label{001Wu}
 For a f.g. projective $A$-module chain complex concentrated in degree $k$
$$C~:~\dots \to 0 \to C_k \to 0 \to \dots$$
the $k$th Wu class defines an isomorphism
$$\widehat{v}_k~:~\widehat{Q}^0(C^{0-*},\epsilon) \xymatrix{\ar[r]^-{\cong}&}
{\rm Hom}_A(C_k,\widehat{H}^k(\Z_2;A,\epsilon))~;~
\gamma \mapsto \widehat{v}_k(\gamma)~.$$
\end{Proposition}
\begin{proof} By construction.
\end{proof}

\begin{Proposition} \label{01Wu}
For a f.g. projective $A$-module chain complex concentrated in degrees $k,k+1$
$$C~:~\dots \to 0 \to C_{k+1} \xymatrix@C-5pt{\ar[r]^-{d}&} C_k \to 0 \to \dots$$
there is defined an exact sequence
$$\begin{array}{l}
{\rm Hom}_A(C_k,\widehat{H}^{k+1}(\Z_2;A,\epsilon))
\xymatrix@C-5pt{\ar[r]^-{d^*}&}
{\rm Hom}_A(C_{k+1},\widehat{H}^{k+1}(\Z_2;A,\epsilon)) \\[1ex]
\xymatrix@C-5pt{\ar[r]&}
\widehat{Q}^0(C^{0-*},\epsilon) \xymatrix@C-5pt{\ar[r]^-{p^*\widehat{v}_k}&}
{\rm Hom}_A(C_k,\widehat{H}^k(\Z_2;A,\epsilon)) \xymatrix@C-5pt{\ar[r]^-{d^*}&}
{\rm Hom}_A(C_{k+1},\widehat{H}^k(\Z_2;A,\epsilon))
\end{array}$$
with $p:C_k \to H_k(C)$ the projection. Thus every chain $\epsilon$-bundle
$(C,\gamma)$ is equivalent to the cone $\Ca(d,\chi)$ $($\ref{cone}$)$
of a chain $\epsilon$-bundle map
$(d,\chi):(C_{k+1},0) \to (C_k,\delta)$, regarding $d:C_{k+1} \to C_k$
as a map of chain complexes concentrated in degree $k$,
with
$$\begin{array}{l}
\delta^*~=~(-1)^k\delta~:~C_k \to C^k~,~d^*\delta d~=~\chi+(-1)^k\chi^*~:~
C_{k+1} \to C^{k+1}~,\\[1ex]
\gamma_{-2k}~=~\delta~:~C_k \to C^k~,~
\gamma_{-2k-1}~=~\begin{cases}
d^*\delta~:~C_k \to C^{k+1}& \\
0~:~C_{k+1} \to C^k&
\end{cases}~,\\[1ex]
\gamma_{-2k-2}~=~\chi~:~C_{k+1} \to C^{k+1}~.
\end{array}$$
\end{Proposition}
\begin{proof} This follows from \ref{01Wu} and the algebraic Thom isomorphisms
$$\widehat{t}~:~\widehat{Q}^*(d,\epsilon)~ \cong~\widehat{Q}^*(C,\epsilon)$$
of Proposition \ref{hyperThom}.
\end{proof}

\subsection{The Twisted Quadratic $Q$-Groups}

For any f.g. projective $A$-module chain complex $C$ there is defined
a short exact sequence of abelian group chain complexes
$$0 \to W_{\%}C \xymatrix{\ar[r]^-{1+T_{\epsilon}}&} W^{\%}C
\xymatrix{\ar[r]^-{J}&} \widehat{W}^{\%}C \to 0$$
with $1+T_{\epsilon}$, $J$ the chain maps
$$\begin{array}{l}
1+T_{\epsilon}~:~W_{\%}C \to W^{\%}C~;~\psi \mapsto (1+T_{\epsilon})\psi~,~
((1+T_{\epsilon})\psi)_s~=~\begin{cases}(1+T_{\epsilon})(\psi_0)&\hbox{if $s=0$}\\
0&\hbox{if $s \geqslant 1$}~,\end{cases}\\
J~:~W^{\%}C \to \widehat{W}^{\%}C~;~\phi \mapsto J\phi~,~
(J\phi)_s~=~\begin{cases}\phi_s&\hbox{if $s\geqslant 0$}\\
0&\hbox{if $s \geqslant -1$}~.\end{cases}
\end{array}$$
The $\epsilon$-symmetric, $\epsilon$-quadratic and $\epsilon$-hyperquadratic $Q$-groups of $C$
are thus related by the exact sequence of Ranicki \cite{ranicki1}
$$\dots \to \widehat{Q}^{n+1} (C,\epsilon) \xymatrix{\ar[r]^-{H}&}
Q_n(C,\epsilon) \xymatrix{\ar[r]^-{1+T_{\epsilon}}&} Q^n (C,\epsilon) \xymatrix{\ar[r]^-{J}&}
\widehat{Q}^n (C,\epsilon) \to \dots$$
with
$$H~:~\widehat{W}^{\%}C \to (W_{\%}C)_{*-1}~;~\theta \mapsto H\theta~,~
(H\theta)_s~=~\theta_{-s-1}~~(s\geqslant 0)~.$$
Weiss \cite{weiss1} used simplicial abelian groups to defined
the twisted quadratic $Q$-groups $Q_*(C,\gamma,\epsilon)$ of
a chain $\epsilon$-bundle $(C,\gamma)$, to fit into the exact sequence
$$\dots \to \widehat{Q}^{n+1} (C,\epsilon) \xymatrix{\ar[r]^-{H_{\gamma}}&}
Q_n(C,\gamma,\epsilon) \xymatrix{\ar[r]^-{N_{\gamma}}&}
Q^n (C,\epsilon) \xymatrix{\ar[r]^-{J_{\gamma}}&} \widehat{Q}^n (C,\epsilon) \to \dots~.$$
The morphisms
$$J_\gamma~:~ Q^n (C,\epsilon) \longrightarrow \widehat{Q}^n (C,\epsilon)~;
\phi \mapsto J_{\gamma}\phi~,~
(J_{\gamma}\phi)_s~=~J(\phi) - (\phi_0)^\% (S^n \gamma)$$
are induced by a morphism of simplicial abelian groups, where
$$ S^n~:~ \widehat{Q}^0 (C^{0-\ast},\epsilon) \xymatrix{\ar[r]^-{\cong}&}
   \widehat{Q}^n (C^{n-\ast},\epsilon)~;~\{\theta_s\} \mapsto \{(S^n\theta)_s=\theta_{s-n}\}$$
are the $n$-fold suspension isomorphisms.

The Kan-Dold theory associates to a chain complex $C$ a simplicial
abelian group $K(C)$ such that
$$\pi_*(K(C))~=~H_*(C)~.$$
For any chain complexes $C,D$ a simplicial map $f:K(C) \to K(D)$
has a mapping fibre $K(f)$. The relative homology groups of $f$ are defined by
$$H_*(f)~=~\pi_{*-1}(K(f))$$
and the fibration sequence of simplicial abelian groups
$$\xymatrix{K(f) \ar[r] & K(C) \ar[r]^-{f} & K(D)}$$
induces a long exact sequence in homology
$$\dots \to H_n(C) \to H_n(D) \to H_n(f) \to H_{n-1}(C) \to \dots~.$$
For a chain map $f:C \to D$
$$K(f)~=~K(\Ca(f))~.$$
The applications involve simplicial maps which are not chain maps,
and the {\it triad homology groups}: given a homotopy-commutative
square of simplicial abelian groups
$$\Phi~:~\raise20pt\hbox{
$\xymatrix{K(C) \ar[r] \ar@{~>}[dr] \ar[d] & K(D) \ar[d] \\
K(E) \ar[r] & K(F)}$}$$
(with $\xymatrix{\ar@{~>}[r]&}$ denoting an explicit homotopy) the
triad homology groups of $\Phi$ are the homotopy groups of the
mapping fibre of the map of mapping fibres
$$H_*(\Phi)~=~\pi_{*-1}(K(C \to D) \to K(E \to F))$$
which fit into a commutative diagram of exact sequences
$$\xymatrix@C-10pt@R-10pt{
& \vdots \ar[d]&\vdots \ar[d] & \vdots \ar[d]&\vdots \ar[d]& \\
\dots \ar[r] & H_{n+1}(D) \ar[r]\ar[d]
& H_{n+1}(C \to D)  \ar[r] \ar[d]& H_n(C)\ar[r] \ar[d]
& H_n(D) \ar[r] \ar[d] & \dots \\
\dots \ar[r] & H_{n+1}(F) \ar[r]\ar[d]
& H_{n+1}(E \to F)  \ar[r] \ar[d]& H_n(E) \ar[r]  \ar[d]
& H_n(F) \ar[r] \ar[d] & \dots \\
\dots \ar[r] & H_{n+1}(D \to F) \ar[r] \ar[d]
& H_{n+1}(\Phi) \ar[r] \ar[d]& H_n(C \to E)\ar[r] \ar[d]
& H_n(D \to F) \ar[r] \ar[d] & \dots \\
\dots \ar[r] & H_n(D) \ar[r] \ar[d]
& H_n(C \to D)  \ar[r] \ar[d]& H_{n-1}(C)  \ar[r] \ar[d]
& H_{n-1}(D) \ar[r] \ar[d] & \dots \\
& \vdots& \vdots & \vdots &\vdots & }$$
If $H_*(\Phi)=0$ there is a commutative braid of exact sequences
\vskip5mm

$$\xymatrix@C-20pt{
H_{n+1}(C \to D)\ar[dr]\ar@/^2pc/[rr]&&
H_n(E)  \ar[dr]\ar@/^2pc/[rr]  &&H_n(C \to E)  \\&
H_n(C)\ar[ur]\ar[dr] &&
\hbox{$H_n(F)$} \ar[ur] \ar[dr]&&\\
H_{n+1}(C \to E) \ar[ur]\ar@/_2pc/[rr]_-{}&&H_n(D)
\ar[ur]\ar@/_2pc/[rr]_{}&&H_n(C \to D)}$$

\indent
The twisted $\epsilon$-quadratic $Q$-groups were defined in \cite{weiss1} to be
the relative homology groups of a simplicial map
$$J_{\gamma}~:~K(W^{\%}C) \to K(\widehat{W}^{\%}C)~,$$
with
$$Q_n(C,\gamma,\epsilon)~=~\pi_{n+1}(J_{\gamma})~.$$
A more explicit description of the twisted quadratic $Q$-groups
was then obtained in Ranicki \cite{ranicki4}, as equivalence
classes of $\epsilon$-symmetric structures on the chain $\epsilon$-bundle.

\begin{Definition}
{\rm
(i) An {\it $\epsilon$-symmetric structure on a chain $\epsilon$-bundle} $(C,\gamma)$ is a
pair $(\phi, \theta)$ with $\phi \in (W^\% C)_n$ a cycle and
$\theta \in (\widehat{W}^\% C)_{n+1}$ such that
$$ d\theta~=~J_\gamma (\phi)~, $$
or equivalently
$$\begin{array}{l}
d\phi_s+(-1)^r\phi_sd^*+
(-1)^{n+s-1}(\phi_{s-1}+(-1)^sT_{\epsilon}\phi_{s-1})~=~0~:~C^r \to C_{n-r+s-1}~,\\[1ex]
\phi_s - \phi_0^* \gamma_{s-n} \phi_0~=~
d\theta_s+(-1)^r\theta_sd^*+(-1)^{n+s}(\theta_{s-1}+(-1)^sT_{\epsilon}
\theta_{s-1})~:~C^r \to C_{n-r+s}\\[1ex]
\hskip120pt (r,s \in \Z,\phi_s=0~{\rm for}~s \leqslant -1)~.
\end{array}$$
(ii) Two structures $(\phi,\theta)$ and $(\phi', \theta')$ are {\it equivalent}
if there exist $\xi \in (W^\% C)_{n+1},$ $\eta \in
(\widehat{W}^\% C)_{n+2}$ such that
$$ d\xi~=~\phi' - \phi,~ d\eta~=~\theta' - \theta + J(\xi) +
   (\xi_0, \phi_0, \phi'_0)^\% (S^n \gamma)~, $$
where $(\xi_0, \phi_0, \phi'_0)^\%:(\widehat{W}^\% C^{-\ast})_n
\rightarrow (\widehat{W}^\% C)_{n+1}$ is the chain homotopy from
$(\phi_0)^\%$ to $(\phi'_0)^\%$ induced by $\xi_0.$ (See
\cite[1.1]{ranicki1} for the precise formula).\\
(iii) The
{\it $n$-dimensional twisted $\epsilon$-quadratic $Q$-group} $Q_n
(C,\gamma,\epsilon)$ is the abelian group of equivalence classes of
$n$-dimensional $\epsilon$-symmetric structures on $(C,\gamma)$ with addition by
$$(\phi,\theta) + (\phi', \theta')~=~(\phi + \phi', \theta + \theta'
+ \zeta)~,~ \text{where } \zeta_s~=~\phi_0 \gamma_{s-n+1} \phi'_0~.$$
\hfill$\qed$}
\end{Definition}

As for the $\pm \epsilon$-symmetric and $\pm\epsilon$-quadratic
$Q$-groups, there are defined skew-suspension isomorphisms of
twisted $\pm \epsilon$-quadratic $Q$-groups
$$\overline{S}~:~
Q_n(C,\gamma,\epsilon) \xymatrix{\ar[r]^-{\displaystyle{\cong}}&}
Q_{n+2}(C_{*-1},\gamma,-\epsilon)~;~
(\{\phi_s\},\{\theta_s\}) \mapsto (\{\phi_s\},\{\theta_s\})~.$$

\begin{Proposition} \label{seq}
{\rm (i)} The twisted $\epsilon$-quadratic
$Q$-groups $Q_*(C,\gamma,\epsilon)$ are related to the $\epsilon$-symmetric
$Q$-groups $Q^*(C,\epsilon)$ and the $\epsilon$-hyperquadratic $Q$-groups
$\widehat{Q}^*(C,\epsilon)$ by the exact sequence
$$ \dots \rightarrow
   \widehat{Q}^{n+1} (C,\epsilon) \stackrel{H_\gamma}{\longrightarrow} Q_n (C,\gamma,\epsilon)
   \stackrel{N_\gamma}{\longrightarrow} Q^n (C,\epsilon)
   \stackrel{J_\gamma}{\longrightarrow} \widehat{Q}^n (C,\epsilon)
   \rightarrow \dots $$
with
$$\begin{array}{l}
H_{\gamma}~:~\widehat{Q}^{n+1}(C,\epsilon) \to Q_n(C,\gamma,\epsilon)~;~\theta
\mapsto (0,\theta)~,\\[1ex]
N_{\gamma}~:~Q_n(C,\gamma,\epsilon) \to Q^n(C,\epsilon)~;~(\phi,\theta) \mapsto
\phi~.
\end{array}$$
{\rm (ii)} For a chain $\epsilon$-bundle $(C, \gamma)$ such that $C$ splits as
$$C~=~\sum_{i=-\infty}^\infty C(i)~,$$
the $\epsilon$-hyperquadratic $Q$-groups split as
$$\widehat{Q}^n(C,\epsilon)~=~\sum_{i=-\infty}^\infty \widehat{Q}^n(C(i),\epsilon)$$
and
$$\gamma ~=~\sum_{i=-\infty}^\infty \gamma (i) \in
\widehat{Q}^0(C^{-*},\epsilon)~=~\sum_{i=-\infty}^\infty \widehat{Q}^0(C(i)^{-*},\epsilon)~.$$
The twisted $\epsilon$-quadratic $Q$-groups of $(C,\gamma)$ fit into the exact
sequence
$$ \begin{array}{l}
\dots \to \sum\limits_iQ_n (C(i),\gamma (i),\epsilon) \xymatrix{\ar[r]^-{q}&}
  Q_n (C,\gamma,\epsilon) \xymatrix{\ar[r]^-{p}&}
  \sum\limits_{i<j} H_n (C(i)\otimes_AC(j))\\[2ex]
\hphantom{\dots \to \sum\limits_i Q_n (C(i),\gamma (i),\epsilon)
  \to Q_n (C,\gamma,\epsilon)}
\xymatrix{\ar[r]^-{\partial}&} \sum\limits_i Q_{n-1} (C(i),\gamma (i),\epsilon)
   \to \dots
\end{array}$$
with
$$\begin{array}{l}
p~:~Q_n (C,\gamma,\epsilon) \to \sum\limits_{i<j} H_n (C(i)\otimes_AC(j))~;~
(\phi,\theta) \mapsto \sum\limits_{i<j} (p(i)\otimes p(j))(\phi_0)\\[1ex]
\hskip50pt (p(i)={\rm projection}:C \to C(i))~, \\[1ex]
q~=~\sum\limits_i q(i)_{\%}~:~
\sum\limits_iQ_n (C(i),\gamma (i),\epsilon) \to Q_n (C,\gamma,\epsilon)~,\\[1ex]
\hskip50pt (q(i)={\rm inclusion}:C(i) \to C))~,\\[1ex]
\partial~:~\sum\limits_{i<j} H_n (C(i)\otimes_AC(j))
\to \sum\limits_i Q_{n-1}(C(i),\gamma(i),\epsilon)~;\\
\hskip50pt \sum\limits_{i<j}h(i,j) \mapsto \sum\limits_i
(0,\sum\limits_{i<j}\widehat{h(i,j)}^{\%}(S^n\gamma(j)))~~
(h(i,j):C(j)^{n-*} \to C(i))~.
\end{array}$$
\end{Proposition}
\begin{proof} (i) See Weiss \cite{weiss1}.\\
(ii) See Ranicki \cite[p.26]{ranicki4}.
\end{proof}

\begin{Example} {\rm The twisted $\epsilon$-quadratic $Q$-groups of the
zero chain $\epsilon$-bundle $(C,0)$ are just the $\epsilon$-quadratic
$Q$-groups of $C$, with isomorphisms
$$Q_n(C,\epsilon) \to Q_n(C,0,\epsilon)~;~ \psi \mapsto ((1+T)\psi,\theta)$$
defined by
$$\theta_s~=~\begin{cases}
\psi_{-s-1}:C^{n-r+s+1} \to C_r&\hbox{\rm if}~ s \leqslant -1 \\
0&\hbox{\rm if}~ s \geqslant 0~,
\end{cases}$$
and with an exact sequence
$$ \dots \rightarrow
   \widehat{Q}^{n+1} (C,\epsilon) \stackrel{H}{\longrightarrow} Q_n (C,\epsilon)
   \stackrel{N}{\longrightarrow} Q^n (C,\epsilon)
   \stackrel{J}{\longrightarrow} \widehat{Q}^n (C,\epsilon)
   \rightarrow \dots~. $$
\hfill$\qed$}
\end{Example}

For $\epsilon=1$ we write
$$\hbox{\rm chain 1-bundle}~=~\hbox{\rm chain bundle}~,~
Q_n(C,\gamma,1)~=~Q_n(C,\gamma)~.$$

\subsection{The Algebraic Normal Invariant}

Fix a chain $\epsilon$-bundle $(B,\beta)$  over $A$.

\begin{Definition} {\rm
(i) A {\it $(B,\beta)$-structure} $(\gamma,\phi,\theta,g,\chi)$ on
an $n$-dimensional ($\epsilon$-symmetric, $\epsilon$-quadratic) Poincar\'e pair
$(f:C \to D, (\delta\phi,\psi))$ is a Spivak normal structure
$(\gamma,\phi,\theta)$ together with a chain $\epsilon$-bundle map
$$(g,\chi)~:~(\Ca(f),\gamma) \to (B,\beta)~.$$
(ii) The {\it $n$-dimensional $(B,\beta)$-structure $\epsilon$-symmetric $L$-group}
$L\langle B,\beta \rangle^n(A,\epsilon)$ is the cobordism group of $n$-dimensional
$\epsilon$-symmetric Poincar\'e complexes $(D,\delta\phi)$ over $A$ together with a
$(B,\beta)$-structure $(\gamma,\phi,\theta,g,\chi)$ (so $(C,\psi)=(0,0)$).\\
(iii) The {\it $n$-dimensional $(B,\beta)$-structure $\epsilon$-hyperquadratic $L$-group}
$\widehat{L}\langle B,\beta \rangle^n(A,\epsilon)$ is the cobordism group of $n$-dimensional
($\epsilon$-symmetric, $\epsilon$-quadratic) Poincar\'e pairs $(f:C \to D, (\delta\phi,\psi))$
over $A$ together with a
$(B,\beta)$-structure $(\gamma,\phi,\theta,g,\chi)$.}
\hfill$\qed$
\end{Definition}

There are defined skew-suspension maps in the $(B,\beta)$-structure $\epsilon$-symmetric
and $\epsilon$-hyperquadratic $L$-groups
$$\begin{array}{l}
\overline{S}~:~L\langle B,\beta \rangle ^n(A,\epsilon) \to
L\langle B_{*-1},\beta_{*-1} \rangle ^{n+2}(A,-\epsilon)~,\\[1ex]
\overline{S}~:~\widehat{L}\langle B,\beta \rangle ^n(A,\epsilon) \to
\widehat{L}\langle B_{*-1},\beta_{*-1} \rangle ^{n+2}(A,-\epsilon)
\end{array}$$
given by $C \mapsto C_{*-1}$ on the chain complexes, with
$(B_{*-1},\beta_{*-1})$ a chain $(-\epsilon)$-bundle. We shall write the
$4$-periodic versions of the $(B,\beta)$-structure $L$-groups as
$$\begin{array}{l}
L\langle B,\beta \rangle ^{n+4*}(A,\epsilon)~=~
\lim\limits_{k \to \infty}L\langle B,\beta \rangle ^{n+4k}(A,\epsilon)~,\\[1ex]
\widehat{L}\langle B,\beta \rangle ^{n+4*}(A,\epsilon)~=~
\lim\limits_{k \to \infty}\widehat{L}\langle B,\beta \rangle ^{n+4k}
(A,\epsilon)~.
\end{array}$$

\begin{Example} {\rm An ($\epsilon$-symmetric, $\epsilon$-quadratic)
Poincar\'e pair with a $(0,0)$-structure is essentially the same as an
$\epsilon$-quadratic Poincar\'e pair.  In particular, an
$\epsilon$-symmetric Poincar\'e complex with a $(0,0)$-structure is
essentially the same as an $\epsilon$-quadratic Poincar\'e complex.
The (0,0)-structure $L$-groups are given by
$$L\langle 0,0 \rangle ^n(A,\epsilon)~=~L_n(A,\epsilon)~,~
\widehat{L}\langle 0,0 \rangle ^n(A,\epsilon)~=~0~.$$
\hfill$\qed$}
\end{Example}

\begin{Proposition}\label{norm}
{\rm (Ranicki \cite[\S7]{ranicki4})}\\
{\rm (i)} An $n$-dimensional $\epsilon$-symmetric structure
$(\phi,\theta) \in Q_n(B,\beta,\epsilon)$ on a chain $\epsilon$-bundle
$(B,\beta)$ determines an $n$-dimensional {\rm (}$\epsilon$-symmetric,
$\epsilon$-quadratic{\rm )} Poincar\'e pair $(f:C \to D,(\delta
\phi,\psi))$ with
$$\begin{array}{l}
f~=~{\rm proj.}~:~C~=~\Ca(\phi_0:B^{n-*} \to B)_{*+1} \to D~=~B^{n-*}~,\\[1ex]
\psi_0~=~\begin{pmatrix} \theta_0 & 0 \\
1+ \beta_{-n}\phi^*_0 & \beta^*_{-n-1}\end{pmatrix}~:\\[3ex]
\hskip75pt C^r~=~B^{r+1} \oplus B_{n-r} \to
C_{n-r-1}~=~B_{n-r} \oplus B^{r+1}~,\\[1ex]
\psi_s~=~\begin{pmatrix} \theta_{-s} & 0 \\
\beta_{-n-s}\phi^*_0 & \beta^*_{-n-s-1}\end{pmatrix}~:\\[3ex]
\hskip75pt C^r~=~B^{r+1} \oplus B_{n-r}
\to C_{n-r-s-1}~=~B_{n-r-s} \oplus B^{r+s+1}~~(s \geqslant 1)~,\\[2ex]
\delta \phi_s~=~\beta_{s-n}~:~D^r~=~B_{n-r} \to D_{n-r+s}~=~B^{r-s}
~~(s \geqslant 0)
\end{array}$$
(up to signs) such that $(\Ca(f),\gamma) \simeq (B,\beta)$.\\
{\rm (ii)} An $n$-dimensional {\rm (}$\epsilon$-symmetric, $\epsilon$-quadratic{\rm )} Poincar\'e pair
$(f:C \to D,(\delta\phi,\psi) \in Q_n^n(f,\epsilon))$
has a canonical equivalence class of `algebraic Spivak normal structures'
$(\gamma,\phi,\theta)$ with $\gamma$ a chain $\epsilon$-bundle over $\Ca(f)$
and $(\phi,\theta)$ an $n$-dimensional $\epsilon$-symmetric structure on $\gamma$
representing an element
$$(\phi,\theta) \in Q_n(\Ca(f),\gamma,\epsilon)~.$$
The construction of {\rm (i)} applied to $(\phi,\theta)$ gives an
$n$-dimensional {\rm (}$\epsilon$-symmetric, $\epsilon$-quadratic{\rm )} Poincar\'e pair homotopy
equivalent to $(f:C \to D,(\delta\phi,\psi) \in Q_n^n(f,\epsilon))$.
\end{Proposition}
\begin{proof} (i) By construction.\\
(ii) The equivalence class $\phi=\delta\phi/(1+T_{\epsilon})\psi \in Q^n(\Ca(f))$ is
given by the algebraic Thom construction
$$\begin{array}{l}
\phi_s~=~\begin{cases}
\begin{pmatrix} \delta \phi_0 & 0 \\
(1+T_{\epsilon})\psi_0f^* & 0\end{pmatrix} &\text{if $s=0$} \\[3ex]
\begin{pmatrix} \delta \phi_1 & 0 \\
0 & (1+T_{\epsilon}) \psi_0\end{pmatrix} &\text{if $s=1$} \\[3ex]
\begin{pmatrix} \delta \phi_s & 0 \\ 0 & 0\end{pmatrix}
&\text{if $s \geqslant 2$}\end{cases}\\[12ex]
\hphantom{\phi_s~=~}
:~\Ca(f)^r~=~D^r \oplus C^{r-1} \to \Ca(f)_{n-r+s}~=~D_{n-r+s} \oplus
C_{n-r+s-1}~,
\end{array}$$
such that
$$\phi_0~:~\Ca(f)^{n-*} \to D^{n-*}
\xymatrix@C+20pt{\ar[r]^-{(\delta\phi,(1+T_{\epsilon})\psi)_0}_-{\simeq}&} \Ca(f)~.$$
The equivalence class $\gamma \in \widehat{Q}^0(\Ca(f)^{0-*},\epsilon)$
of the Spivak normal chain bundle is
the image of $(\delta\phi,\psi) \in Q^n_n(f,\epsilon)$ under the composite
$$Q^n_n(f,\epsilon) \xymatrix{\ar[r]^-{J_f}&} \widehat{Q}^n(D,\epsilon)
\xymatrix@C+30pt{\ar[r]^-{((\delta\phi,(1+T_{\epsilon})\psi)^{\%}_0)^{-1}}_-
{\displaystyle{\cong}}&} \widehat{Q}^n(\Ca(f)^{n-*},\epsilon)
\xymatrix{\ar[r]^-{S^{-n}}_-{\displaystyle{\cong}}&}
\widehat{Q}^0(\Ca(f)^{0-*},\epsilon)~.$$
\end{proof}

\begin{Definition} \label{norminv} {\rm
(i) The {\it boundary} of an $n$-dimensional $\epsilon$-symmetric
structure $(\phi,\theta) \in Q_n(B,\beta,\epsilon)$ on a chain $\epsilon$-bundle
$(B,\beta)$ over $A$ is the $\epsilon$-symmetric null-cobordant
$(n-1)$-dimensional $\epsilon$-quadratic Poincar\'e complex over $A$
$$\partial(\phi,\theta)~=~(C,\psi)$$
defined in Proposition \ref{norm} (i) above, with $C=\Ca(\phi_0:B^{n-*} \to B)_{*+1}$.\\
(ii) The {\it algebraic normal invariant} of an $n$-dimensional ($\epsilon$-symmetric,
$\epsilon$-quadratic) Poincar\'e pair over $A$
$(f:C \to D,(\delta\phi,\psi) \in Q_n^n(f,\epsilon))$ is the class
$$(\phi,\theta) \in Q_n(\Ca(f),\gamma,\epsilon)$$
defined in Proposition \ref{norm} (ii) above.
\hfill$\qed$}
\end{Definition}

\begin{Proposition} \label{bformation}
Let $(B,\beta)$ be a chain $\epsilon$-bundle over $A$
such that $B$ is concentrated in degree $k$
$$B~:~\dots \to 0 \to B_k \to 0 \to \dots~.$$
The boundary map
$\partial:Q_{2k}(B,\beta,\epsilon) \to L_{2k-1}(A,\epsilon)$
sends an $\epsilon$-symmetric structure $(\phi,\theta) \in Q_{2k}(B,\beta,\epsilon)$
to the Witt class of the $(-1)^{k-1}\epsilon$-quadratic formation
$$\partial(\phi,\theta)~=~
(H_{(-1)^{k-1}\epsilon}(B^k);B^k,{\rm im}(\begin{pmatrix} 1-\beta\phi \\
\phi \end{pmatrix}:B^k \to B^k \oplus B_k))$$
with
$$H_{(-1)^{k-1}\epsilon}(B^k)~=~(B^k \oplus B_k,\begin{pmatrix} 0 & 1 \\ 0 & 0
\end{pmatrix})$$
the hyperbolic $(-1)^{k-1}\epsilon$-quadratic form.
\end{Proposition}
\begin{proof} The chain $\epsilon$-bundle (equivalence class)
$$\beta \in \widehat{Q}^0(B^{0-*},\epsilon)~=~
\widehat{H}^0(\Z_2;S(B_k),\epsilon)$$
is represented by an $\epsilon$-symmetric
form $(B_k,\beta)$. An $\epsilon$-symmetric structure
$(\phi,\theta) \in Q_{2k}(B,\beta,\epsilon)$
is represented by an $(-1)^k\epsilon$-symmetric form $(B^k,\phi)$
together with $\theta \in S(B_k)$ such that
$$\phi - \phi \beta \phi~=~\theta+(-1)^k\epsilon \theta^* \in
H^0(\Z_2;S(B^k),(-1)^k\epsilon)~.$$
The boundary of $(\phi,\theta)$ is the $\epsilon$-symmetric null-cobordant
$(2k-1)$-dimensional $\epsilon$-quadratic Poincar\'e complex
$\partial(\phi,\theta)=(C,\psi)$ concentrated in degrees $k-1,k$
corresponding to the formation in the statement.
\end{proof}

\begin{Proposition} \label{bform}
Let $(B,\beta)$ be a chain $\epsilon$-bundle over $A$
such that $B$ is concentrated in degrees $k,k+1$
$$B~:~\dots \to 0 \to B_{k+1} \xymatrix{\ar[r]^d&} B_k \to 0 \to \dots~.$$
The boundary map
$\partial:Q_{2k+1}(B,\beta,\epsilon) \to L_{2k}(A,\epsilon)$
sends an $\epsilon$-symmetric structure $(\phi,\theta) \in
Q_{2k+1}(B,\beta,\epsilon)$ to the Witt class of the
nonsingular $(-1)^k\epsilon$-quadratic form over $A$
$$({\rm coker}(\begin{pmatrix} -d^* \\ \phi^*_0 \\ 1-\beta_{-2k}d\phi^*_0
\end{pmatrix}:B^k \to B^{k+1} \oplus B_{k+1} \oplus B^k),
\begin{pmatrix} \theta_0 & 0 & \phi_0 \\ 1 & \beta^*_{-2k-2} &
d^* \\ 0 & 0 & 0 \end{pmatrix})~.$$
\end{Proposition}
\begin{proof} This is an application of the instant surgery obstruction
of \cite[4.3]{ranicki1}, which identifies the cobordism class
$(C,\psi) \in L_{2k}(A,\epsilon)$ of a $2k$-dimensional $\epsilon$-quadratic
Poincar\'e complex $(C,\psi)$ with the Witt class of the nonsingular
$\epsilon$-quadratic form
$$I(C,\psi)~=~
({\rm coker}(\begin{pmatrix} d^* \\ (-1)^{k+1}(1+T_{\epsilon})\psi_0 \end{pmatrix}:
C^{k-1} \to C^k \oplus C_{k+1}),\begin{pmatrix} \psi_0 & d \\ 0 & 0 \end{pmatrix})~.$$
\indent By Proposition \ref{01Wu} the chain $\epsilon$-bundle
$\beta$ can be taken to be the cone of a chain $\epsilon$-bundle map
$$(d,\beta_{-2k-2})~:~(B_{k+1},0) \to (B_k,\beta_{-2k})$$
with
$$\begin{array}{l}
\beta^*_{-2k}~=~(-1)^k\epsilon \beta_{-2k}~:~B_k \to B^k~,\\[1ex]
d^* \beta_{-2k} d~=~\beta_{-2k-2}+ (-1)^k\epsilon\beta^*_{-2k-2}~:~B_{k+1} \to
B^{k+1}~,\\[1ex]
\beta_{-2k-1}~=~\begin{cases}
\beta_{-2k}d~:~B_{k+1} \to B^k&\\
0~:~B_k \to B^{k+1}~.&
\end{cases}$$
\end{array}$$
An  $\epsilon$-symmetric structure
$(\phi,\theta) \in Q_{2k+1}(B,\beta,\epsilon)$ is represented
by $A$-module morphisms
$$\begin{array}{l}
\phi_0~:~B^k \to B_{k+1}~,~\widetilde{\phi}_0~:~B^{k+1} \to B_k~,~
\phi_1~:~B^{k+1} \to B_{k+1}~,\\[1ex]
\theta_0~:~B^{k+1} \to B_{k+1}~,~\theta_{-1}~:~B^k \to B_{k+1}~,~
\widetilde{\theta}_{-1}~:~B^{k+1} \to B_k~,~\theta_{-2}~:~B^k \to B_k
\end{array}$$
such that
$$\begin{array}{l}
d\phi_0+(-1)^k\widetilde{\phi}_0d^*~=~0~:~B^k \to B_k~,\\[1ex]
\phi_0 - \epsilon\widetilde{\phi}_0^* + (-1)^{k+1}\phi_1d^*~=~0~:~B^k \to B_{k+1}~,\\[1ex]
\phi_1+(-1)^{k+1}\epsilon\phi_1^*~=~0~:~B^{k+1} \to B_{k+1}~,\\[1ex]
\phi_0 - \phi_0 \beta_{-2k}d \widetilde{\phi}_0^*~=~
(-1)^k\theta_0d^*-\theta_{-1}-\epsilon\widetilde{\theta}_{-1}^*~:~
B^k \to B_{k+1}~,\\[1ex]
\widetilde{\phi}_0~=~d\theta_0-\widetilde{\theta}_{-1}-
\epsilon\theta^*_{-1}~:~B^{k+1} \to B_k~,\\[1ex]
-\widetilde{\phi}_0 \beta_{-2k-2} \widetilde{\phi}_0^*~=~
\theta_{-2}+(-1)^{k+1}\epsilon\theta_{-2}^*~:~B^k \to B_k~,\\[1ex]
\phi_1- \phi_0 \beta_{-2k} \phi_0^*~=~\theta_0+(-1)^k\epsilon\theta_0^*~:~
B^{k+1} \to B_{k+1}~.
\end{array}$$
The boundary of $(\phi,\theta)$ given by \ref{norminv} (i)
is an $\epsilon$-symmetric null-cobordant
$2k$-dimensional $\epsilon$-quadratic Poincar\'e complex
$\partial(\phi,\theta)=(C,\psi)$ concentrated in degrees $k-1,k,k+1$,
with $I(C,\psi)$ the instant surgery obstruction form (\ref{bform})
in the statement.
\end{proof}

The $\epsilon$-quadratic $L$-groups and the
$(B,\beta)$-structure $L$-groups fit into an evident exact sequence
$$\dots \to L_n(A,\epsilon) \to L\langle B,\beta \rangle^n(A,\epsilon)
\to \widehat{L}\langle B,\beta \rangle^n(A,\epsilon)
\xymatrix{\ar[r]^-{\partial}&} L_{n-1}(A,\epsilon) \to \dots~,$$
and similarly for the 4-periodic versions
$$\dots \to L_n(A,\epsilon) \to L\langle B,\beta \rangle^{n+4*}(A,\epsilon)
\to \widehat{L}\langle B,\beta \rangle^{n+4*}(A,\epsilon)
\xymatrix{\ar[r]^-{\partial}&}  L_{n-1}(A,\epsilon) \to \dots~.$$

\begin{Proposition} \label{cob} {\rm (Weiss \cite{weiss1})}\\
{\rm (i)} The function
$$Q_n(B,\beta,\epsilon) \to \widehat{L}\langle B,\beta \rangle^{n+4*}(A,\epsilon)~;~
(\phi,\theta) \mapsto (f:C \to D,(\delta\phi,\psi))~(\ref{norm}~{\rm (ii)})$$
is an isomorphism, with inverse given by the algebraic normal invariant.
The $\epsilon$-quadratic $L$-groups of $A$, the 4-periodic
$(B,\beta)$-structure $\epsilon$-symmetric $L$-groups of $A$ and the twisted
$\epsilon$-quadratic $Q$-groups of $(B,\beta)$ are thus related by an exact
sequence
$$\dots \rightarrow L_n (A,\epsilon) \xymatrix{\ar[r]^-{1+T}&}
   L\langle B,\beta \rangle^{n+4*} (A,\epsilon) \rightarrow
   Q_n (B,\beta,\epsilon) \stackrel{\partial}{\rightarrow}
   L_{n-1} (A,\epsilon) \rightarrow \dots~.$$
{\rm (ii)} The cobordism class of an $n$-dimensional {\rm (}$\epsilon$-symmetric, $\epsilon$-quadratic{\rm )}
Poincar\'e pair $(f:C \to D,(\delta\phi,\psi))$ over $A$ with a
$(B,\beta)$-structure $(\gamma,\phi,\theta,g,\chi)$ is the
image of the algebraic normal invariant $(\phi,\theta) \in Q_n(\Ca(f),\gamma,\epsilon)$
$$(f:C \to D,(\delta\phi,\psi))~=~(g,\chi)_{\%}(\phi,\theta)
\in Q_n (B,\beta,\epsilon)~.$$
\end{Proposition}
\begin{proof} The $\epsilon$-symmetrization of an $n$-dimensional $\epsilon$-quadratic Poincar\'e complex $(C,\psi)$
is an $n$-dimensional $\epsilon$-symmetric Poincar\'e complex $(C,(1+T_{\epsilon})\psi)$ with
$(B,\beta)$-structure $(0,(1+T)\psi,\theta,0,0)$ given by
$$\theta_s~=~\begin{cases}
\psi_{-s-1} \in {\rm Hom}_A(C^{-*},C)_{n+s+1} &\hbox{\rm if $s \leqslant -1$}\\
0&\hbox{\rm if $s \geqslant 0$}~.
\end{cases}$$
The relative groups of the symmetrization map
$$1+T_{\epsilon}~:~L_n(A,\epsilon) \to L\langle B,\beta \rangle^n(A,\epsilon)~;~(C,\psi) \mapsto
(C,(1+T_{\epsilon})\psi)$$
are the cobordism groups of $n$-dimensional ($\epsilon$-symmetric, $\epsilon$-quadratic) Poincar\'e
pairs $(f:C \to D,(\delta\phi,\psi))$ together with
a $(B,\beta)$-structure $(\gamma,\phi,\theta,g,\chi)$.
\end{proof}

\begin{Proposition} \label{0dim}
Let $(B,\beta)$ be a chain $\epsilon$-bundle over $A$ with $B$ concentrated in degree $k$
$$B~:~\dots \to 0 \to B_k \to 0 \to \dots$$
so that
$\beta \in \widehat{Q}^0(B^{0-*},\epsilon)=\widehat{H}^0(\Z_2;S(B^k),(-1)^kT_{\epsilon})$
is represented by an element
$$\beta_{-2k}~=~(-1)^k\epsilon\beta^*_{-2k} \in S(B^k)~.$$
The twisted $\epsilon$-quadratic $Q$-groups $Q_n(B,\beta,\epsilon)$ are
given as follows.\\
{\rm (i)} For $n\neq 2k-1,2k$
$$\begin{array}{ll}
Q_n(B,\beta,\epsilon)&=~Q_n(B,\epsilon)\\[1ex]
&=~\begin{cases}
\widehat{Q}^{n+1}(B,\epsilon)=\widehat{H}^{n-2k+1}(\Z_2;S(B^k),(-1)^kT_{\epsilon})
&\hbox{\it if $n \geqslant 2k+1$}
\\[2ex]
0&\hbox{\it if $n \leqslant 2k-2$}~.
\end{cases}
\end{array}$$
{\rm (ii)} For $n=2k$
$$Q_{2k}(B,\beta,\epsilon)~=~\dfrac{
\{(\phi,\theta) \in S(B^k) \oplus S(B^k)\,\vert\, \phi=(-1)^k\epsilon\phi^*,
\phi-\phi\beta_{-2k}\phi^*=\theta+(-1)^k\epsilon\theta^*\}}
{\{(0,\eta+(-1)^{k+1}\epsilon\eta^*)\,\vert\, \eta \in S(B^k)\}}$$
with addition by
$$(\phi,\theta)+(\phi',\theta')~=~(\phi+\phi',\theta+\theta'+
\phi' \beta_{-2k} \phi^*)~.$$
The boundary of $(\phi,\theta) \in Q_{2k}(B,\beta,\epsilon)$ is the
$(2k-1)$-dimensional $\epsilon$-quadratic Poincar\'e complex over $A$
concentrated in degrees $k-1,k$
corresponding to the $(-1)^{k+1}\epsilon$-quadratic formation over $A$
$$\partial(\phi,\theta)~=~
(H_{(-1)^{k+1}\epsilon}(B^k);B^k,{\rm im}(\begin{pmatrix} 1-\beta_{-2k}\phi \\
\phi \end{pmatrix}:B^k \to B^k \oplus B_k))~.$$
{\rm (iii)} For $n=2k-1$
$$\begin{array}{ll}
Q_{2k-1}(B,\beta,\epsilon)~&=~{\rm coker}(J_{\beta}:Q^{2k}(B,\epsilon)
\to \widehat{Q}^{2k}(B,\epsilon))\\[1ex]
&=~\dfrac{
\{\sigma  \in S(B^k)\,\vert\, \sigma=(-1)^k\epsilon\sigma^*\}}
{\{\phi-\phi\beta_{-2k}\phi^*-(\theta+(-1)^k\epsilon\theta^*)
\,\vert\,\phi=(-1)^k\epsilon\phi^*,\theta \in S(B^k)\}}~.
\end{array}$$
The boundary of $\sigma \in Q_{2k-1}(B,\beta,\epsilon)$ is the
$(2k-2)$-dimensional $\epsilon$-quadratic Poincar\'e complex over $A$
concentrated in degree $k-1$
corresponding to the $(-1)^{k+1}\epsilon$-quadratic form over $A$
$$\partial(\sigma)~=~(B^k \oplus B_k,\begin{pmatrix} \sigma & 1 \\
0 & \beta_{-2k}\end{pmatrix})$$
with
$$(1+T_{(-1)^{k+1}\epsilon})\partial(\sigma)~=~
(B^k \oplus B_k,\begin{pmatrix} 0 & 1 \\
(-1)^{k+1}\epsilon & 0 \end{pmatrix})~.$$
{\rm (iv)} The maps in the exact sequence
$$\begin{array}{l}
0 \to \widehat{Q}^{2k+1}(B,\epsilon)
\xymatrix@C-5pt{\ar[r]^-{H_{\beta}}&} Q_{2k}(B,\beta,\epsilon)
\xymatrix@C-5pt{\ar[r]^-{N_{\beta}}&} Q^{2k}(B,\epsilon)\\[1ex]
\hphantom{0 \to \widehat{Q}^{2k+1}(B,\epsilon)
\xymatrix@C-5pt{\ar[r]^-{H_{\beta}}&} Q_{2k}(B,\beta,\epsilon)}
\xymatrix@C-5pt{\ar[r]^-{J_{\beta}}&} \widehat{Q}^{2k}(B,\epsilon)
\xymatrix@C-5pt{\ar[r]^-{H_{\beta}}&} Q_{2k-1}(B,\beta,\epsilon)
\to 0
\end{array}$$
are given by
$$\begin{array}{l}
H_{\beta}~:~\widehat{Q}^{2k+1}(B,\epsilon)=\widehat{H}^1(\Z_2;S(B^k),(-1)^kT_{\epsilon})
\to Q_{2k}(B,\beta,\epsilon)~;~ \theta \mapsto (0,\theta)~,\\[1ex]
N_{\beta}~:~Q_{2k}(B,\beta,\epsilon) \to Q^{2k}(B,\epsilon)=H^0(\Z_2;S(B^k),(-1)^kT_{\epsilon})~;~
(\phi,\theta) \mapsto \phi~,\\[1ex]
J_{\beta}~:~Q^{2k}(B,\epsilon)=H^0(\Z_2;S(B^k),(-1)^kT_{\epsilon}) \to
\widehat{Q}^{2k}(B,\epsilon)=\widehat{H}^0(\Z_2;S(B^k),(-1)^kT_{\epsilon})~;\\[1ex]
\hphantom{N_{\beta}~:~Q_{2k}(B,\beta,\epsilon) \to Q^{2k}(B,\epsilon)=H^0(\Z_2;S(B^k),(-1)^kT_{\epsilon})~;~}
\phi \mapsto \phi - \phi \beta_{-2k}\phi^*~,\\[1ex]
H_{\beta}~:~\widehat{Q}^{2k}(B,\epsilon)=\widehat{H}^0(\Z_2;S(B^k),(-1)^kT_{\epsilon}) \to
Q_{2k-1}(B,\beta,\epsilon)~;~ \sigma \mapsto \sigma~.
\end{array}$$
\hfill$\qed$
\end{Proposition}

\begin{Example} \label{Spiv}
{\rm Let $(K,\lambda)$ be a nonsingular $\epsilon$-symmetric
form over $A$, which may be regarded as a 0-dimensional $\epsilon$-symmetric
Poincar\'e complex $(D,\phi)$ over $A$ with
$$\phi_0~=~\lambda~:~D^0~=~K \to D_0~=~K^*~.$$
The composite
$$Q^0(D,\epsilon)~=~H^0(\Z_2;S(K),\epsilon)
\xymatrix{\ar[r]^J&} \widehat{Q}^0(D,\epsilon)
\xymatrix{\ar[r]^{(\phi_0)^{-1}}&} \widehat{Q}^0(D^{0-*},\epsilon)$$
sends $\phi \in Q^0(D,\epsilon)$ to the algebraic Spivak normal chain bundle
$$\gamma \in \widehat{Q}^0(D^{0-*},\epsilon)~=~\widehat{H}^0(\Z_2;S(K^*),\epsilon) $$
with
$$\gamma_0~=~\epsilon\lambda^{-1} ~:~D_0~=~K^* \to D^0~=~K~.$$
By Proposition \ref{0dim}
$$Q_0(D,\gamma,\epsilon)~=~\dfrac{
\{(\kappa,\theta) \in S(K) \oplus S(K)\,\vert\, \kappa=\epsilon\kappa^*,
\kappa-\kappa\gamma_0\kappa^*=\theta+\epsilon\theta^*\}}
{\{(0,\eta-\epsilon\eta^*)\,\vert\, \eta \in S(K)\}}$$
with addition by
$$(\kappa,\theta)+(\kappa',\theta')~=~(\kappa+\kappa',\theta+\theta'+
\kappa' \gamma_0 \kappa^*)~.$$
The algebraic normal invariant of $(D,\phi)$ is given by
$$(\phi,0) \in Q_0(D,\gamma,\epsilon)~.$$
\hfill$\qed$}
\end{Example}

\begin{Example} \label{0dim.twq}
{\rm Let $A$ be a ring with even involution (\ref{even}),
and let $C$ be concentrated in degree $k$ with $C_k=A^r$.
For odd $k=2j+1$
$$\widehat{Q}^0(C^{0-*})~=~0$$
and there is only one chain $\epsilon$-bundle $\gamma=0$ over $C$, with
$$Q_n(C,\gamma)~=~Q_n(C)~=~
\begin{cases}
\bigoplus\limits_r
\widehat{H}^0(\Z_2;A)&\hbox{if $n\geqslant 4j+2$, $n \equiv 0(\bmod\, 2)$}\\
0&\hbox{otherwise~.}
\end{cases}$$
For even $k=2j$
$$\widehat{Q}^0(C^{0-*})~=~\bigoplus\limits_r\widehat{H}^0(\Z_2;A)$$
a chain $\epsilon$-bundle
$\gamma \in \widehat{Q}^0(C^{0-*})$ is represented by a diagonal matrix
$$\gamma~=~X~=~\begin{pmatrix} x_1 & 0 & 0 & \dots & 0 \\
0 & x_2 & 0 & \dots & 0 \\
0 & 0 & x_3 & \dots & 0 \\
\vdots & \vdots & \vdots & \ddots & \vdots \\
0 & 0 & 0 & \dots & x_r \end{pmatrix} \in {\rm Sym}_r(A)$$
with $\overline{x}_i=x_i \in A$, and there is defined an exact sequence
$$\widehat{Q}^{4j+1}(C)=0 \to Q_{4j}(C,\gamma) \to Q^{4j}(C)
\xymatrix{\ar[r]^-{J_{\gamma}}&} \widehat{Q}^{4j}(C) \to
Q_{4j-1}(C,\gamma) \to Q^{4j-1}(C)=0$$
with
$$J_{\gamma}~:~Q^{4j}(C)~=~{\rm Sym}_r(A) \to
\widehat{Q}^{4j}(C)~=~\dfrac{{\rm Sym}_r(A)}{{\rm Quad}_r(A)}~;~
M \mapsto M - M X M~,$$
so that
$$\begin{array}{l}
Q_n(C,\gamma)~=\\[1ex]
\begin{cases}
\bigoplus\limits_r
\widehat{H}^0(\Z_2;A) &
\hbox{\rm if $n \geqslant 4j+1$}\\[-1ex]
&\hskip10pt\hbox{\rm and $n \equiv 1(\bmod\, 2)$}\\
\{M \in {\rm Sym}_r(A)\,\vert\, M-MXM \in {\rm Quad}_r(A)\}
&\hbox{\rm if $n=4j$}\\
M_r(A)/\{M-MXM-(N+N^t)
\,\vert\,M\in {\rm Sym}_r(A),N \in M_r(A)\}
&\hbox{\rm if $n=4j-1$}\\
0&\hbox{\rm otherwise}~.
\end{cases}
\end{array}
$$
Moreover, Proposition \ref{seq}  (ii) gives an exact sequence
$$0 \to \bigoplus\limits_{i=1}^rQ_{4j}(B,x_i) \to
Q_{4j}(C,\gamma) \to \bigoplus\limits_{r(r-1)/2}A
\to \bigoplus\limits_{i=1}^rQ_{4j-1}(B,x_i) \to Q_{4j-1}(C,\gamma) \to 0$$
with $B$ concentrated in degree $2j$ with $B_{2j}=A$.\hfill$\qed$}
\end{Example}

\subsection{The Relative Twisted Quadratic $Q$-groups}

Let $(f,\chi):(C,\gamma)\to (D,\delta)$ be a map of chain
$\epsilon$-bundles, and let $(\phi, \theta)$ be an $n$-dimensional
$\epsilon$-symmetric structure on $(C,\gamma)$, so that $\chi \in
(\widehat{W}^{\%}C)_1$, $\phi \in (W^\% C)_n$ and $\theta \in
(\widehat{W}^\% C)_{n+1}$.  Composing the chain map $\phi_0:C^{n-\ast}
\to C$ with $f$, we get an induced map
$$(\widehat{f \phi_0})^\%~:~ \widehat{W}^\% C^{n-\ast} \to \widehat{W}^\% D~.$$
The morphisms of twisted quadratic $Q$-groups
$$ (f,\chi)_\%~:~ Q_n (C,\gamma,\epsilon) \to Q_n (D, \delta,\epsilon)~;~
(\phi,\theta) \mapsto (f^\% (\phi), \widehat{f}^\% (\theta)
   + (\widehat{f \phi_0})^\% (S^n \chi))$$
are induced by a simplicial map of simplicial abelian groups.
The relative homotopy groups are the {\it relative twisted
$\epsilon$-quadratic $Q$-groups}
$Q_n (f,\chi,\epsilon)$, designed to fit into a long exact sequence
$$ \dots \longrightarrow Q_n (C,\gamma,\epsilon)
   \stackrel{(f,\chi)_\%}{\longrightarrow}
   Q_n (D, \delta,\epsilon) \longrightarrow Q_n (f,\chi,\epsilon)
   \longrightarrow Q_{n-1} (C,\gamma,\epsilon)
   \longrightarrow \dots. $$

\begin{Proposition}
For any chain $\epsilon$-bundle map $(f,\chi):(C,\gamma) \to (D,\delta)$ the
various $Q$-groups fit into a commutative diagram with exact rows and
columns
$$\xymatrix{
& \vdots \ar[d]&\vdots \ar[d] & \vdots \ar[d]&\vdots \ar[d]& \\
\dots \ar[r] & \widehat{Q}^{n+1}(C,\epsilon)
\ar[r]^-{\displaystyle{H_\gamma}}  \ar[d]^-{\displaystyle{\widehat{f}^{\%}}}
& Q_n(C,\gamma,\epsilon)  \ar[r]^-{\displaystyle{N_\gamma}} \ar[d]^-{\displaystyle{(f,\chi)}_{\%}}&
Q^n(C,\epsilon)  \ar[r]^-{\displaystyle{J_\gamma}} \ar[d]^-{\displaystyle{f^{\%}}}
& \widehat{Q}^n(C,\epsilon) \ar[r]
\ar[d]^-{\displaystyle{\widehat{f}^{\%}}} & \dots \\
\dots \ar[r] & \widehat{Q}^{n+1}(D,\epsilon) \ar[r]^-{\displaystyle{H_\delta}} \ar[d]
& Q_n(D,\delta,\epsilon)  \ar[r]^-{\displaystyle{N_\gamma}} \ar[d]&
Q^n(D,\epsilon) \ar[r]^-{\displaystyle{J_\delta}}  \ar[d]
& \widehat{Q}^n(D,\epsilon) \ar[r] \ar[d] & \dots \\
\dots \ar[r] & \widehat{Q}^{n+1}(f) \ar[r]^-{\displaystyle{H_\chi}} \ar[d]
& Q_n(f,\chi,\epsilon)  \ar[r]^-{\displaystyle{N_\chi}} \ar[d]&
Q^n(f,\epsilon)  \ar[r]^-{\displaystyle{J_\chi}} \ar[d]
& \widehat{Q}^n(f,\epsilon) \ar[r] \ar[d] & \dots \\
\dots \ar[r] & \widehat{Q}^n(C,\epsilon) \ar[r]^-{\displaystyle{H_\gamma}} \ar[d]
& Q_{n-1}(C,\gamma,\epsilon)  \ar[r]^-{\displaystyle{N_\delta}} \ar[d]&
Q^{n-1}(C,\epsilon)  \ar[r]^-{\displaystyle{J_\gamma}} \ar[d]
& \widehat{Q}^{n-1}(C,\epsilon) \ar[r] \ar[d] & \dots \\
& \vdots& \vdots & \vdots &\vdots & }$$
\end{Proposition}
\begin{proof} These are the exact sequences of the homotopy groups
of the simplicial abelian groups in the commutative diagram of fibration
sequences
$$\xymatrix{K(J_{\gamma}) \ar[r] \ar[d]
& K(W^{\%}C) \ar[r]^-{J_{\gamma}} \ar[d]^-{f^{\%}}
& K(\widehat{W}^{\%}C) \ar[d]^-{\widehat{f}^{\%}} \\
K(J_{\delta}) \ar[r] \ar[d]
& K(W^{\%}D) \ar[r]^-{J_{\delta}} \ar[d]
& K(\widehat{W}^{\%}D) \ar[d] \\
K(J_{\chi}) \ar[r]
& K(W^{\%}\Ca(f)) \ar[r]^-{J_{\chi}}
& K(\widehat{W}^{\%}\Ca(f)) }$$
with
$$\pi_n(K(J_{\chi}))~=~Q_n(f,\chi,\epsilon)~.$$
\end{proof}

There is also a twisted $\epsilon$-quadratic $Q$-group version of the algebraic
Thom constructions (\ref{hyperThom}, \ref{symmThom}, \ref{quadThom})~:

\begin{Proposition} \label{twquThom}
Let $(f,\chi):(C,0) \to (D,\delta)$ be a chain $\epsilon$-bundle map,
and let $(B,\beta)=\Ca(f,\chi)$ be the cone chain $\epsilon$-bundle (\ref{cone}).
The relative twisted $\epsilon$-quadratic $Q$-groups
$Q_*(f,\chi,\epsilon)$ are related to the {\rm (}absolute{\rm )}
twisted $\epsilon$-quadratic $Q$-groups $Q_*(B,\beta,\epsilon)$ by a
commutative braid of exact sequences \vskip1mm

$$\xymatrix@!C@C-40pt@R-10pt{
\widehat{Q}^{n+1}(B,\epsilon) \ar[dr]\ar@/^2pc/[rr]&&
Q_n(B,\beta,\epsilon) \ar[dr]\ar@/^2pc/[rr]&& H_{n-1}(B\otimes_AC)\\
& Q_n(f,\chi,\epsilon) \ar[ur]^-{t}\ar[dr]&& Q^n(B,\epsilon)\ar[ur]\ar[dr]&\\
H_n(B\otimes_AC)\ar[ur]^-{F}\ar@/_2pc/[rr]^-{F} &&Q^n(f,\epsilon)\ar[ur]^-{t}
\ar@/_2pc/[rr] && \widehat{Q}^n(B,\epsilon)}$$

\vskip5mm

\noindent involving the exact sequence of \ref{symmThom}
$$\dots \to H_n(B\otimes_AC) \xymatrix{\ar[r]^-{F}&}
Q^n(f,\epsilon)
\xymatrix{\ar[r]^-{t}&} Q^n(B,\epsilon) \to H_{n-1}(B\otimes_AC) \to \dots~.$$
\end{Proposition}
\begin{proof} The identity
$$\widehat{f}^{*\%}(\delta)~=~d \chi \in (\widehat{W}C^{0-*})_0$$
determines a homotopy $\xymatrix{\ar@{~>}[r]&}$ in the square
$$\raise20pt\hbox{$\xymatrix{K(W^{\%}C) \ar[r]^-{\displaystyle{J}} \ar[d]^-{\displaystyle{f^{\%}}}
\ar@{~>}[dr]&
K(\widehat{W}^{\%}C) \ar[d]^-{\displaystyle{\widehat{f}^{\%}}}\\
K(W^{\%}D) \ar[r]^-{\displaystyle{J_{\delta}}} &
K(\widehat{W}^{\%}D)}$}$$
(with $J=J_0$) and hence maps of the mapping fibres
$$J_{\chi}~:~K(\Ca(f^{\%})) \to K(\Ca(\widehat{f}^{\%}))~,~
(f,\chi)_{\%}~:~K(J) \to K(J_{\delta})~.$$
The map $J_{\chi}$ is related to $J_\beta:K(W^{\%}B) \to K(\widehat{W}^{\%}B)$
by a homotopy commutative diagram
$$\raise20pt\hbox{$\xymatrix{K(\Ca(f^{\%}))
\ar[r]^-{\displaystyle{J_{\chi}}} \ar[d]^-{\displaystyle{t}} \ar@{~>}[dr]&
K(\Ca(\widehat{f}^{\%})) \ar[d]_{\displaystyle{\simeq}}^-{\displaystyle{\widehat{t}}}\\
K(W^{\%}B) \ar[r]^{\displaystyle{J_\beta}} & K(\widehat{W}^{\%}B)}$}$$
with
$\widehat{t}:K(\Ca(\widehat{f}^{\%})) \simeq K(\widehat{W}^{\%}B)$
a simplicial homotopy equivalence inducing the algebraic Thom
isomorphisms $\widehat{t}:\widehat{Q}^*(f,\epsilon) \cong \widehat{Q}^*(B,\epsilon)$ of
Proposition \ref{hyperThom}, and $t:K(\Ca(f^{\%})) \to K(W^{\%}B)$ a
simplicial map inducing the algebraic Thom maps $t:Q^*(f,\epsilon) \to Q^*(B,\epsilon)$
of Proposition \ref{symmThom}, with mapping fibre $K(t)\simeq K(B \otimes_AC)$.
The braid in the statement is the commutative braid of homotopy groups
induced by the homotopy commutative braid of fibrations

$$\xymatrix@!C@C-40pt@R-10pt{
&& K(J_\beta) \ar[dr]&& \\
& K(J_\chi) \ar[ur]\ar[dr]\ar@{~>}[rr]&& K(W^{\%}B)\ar[dr]^-{J_\beta}&\\
K(B\otimes_AC)\ar[ur]^-{F}\ar@/_2pc/[rr]^-{F} &&K(\Ca(f^{\%}))\ar[ur]^-{t}
\ar@/_2pc/[rr]^-{\widehat{t}J_{\chi}} && K(\widehat{W}^{\%}B)}$$

\end{proof}

\begin{Proposition} \label{twqcor}
Let $(C,\gamma)$ be a chain $\epsilon$-bundle over a
f.g. projective $A$-module chain complex which is concentrated in degrees $k,k+1$
$$C~:~\dots \to 0 \to C_{k+1} \xymatrix{\ar[r]^d&} C_k \to 0 \to \dots~,$$
so that $(C,\gamma)$ can be taken $($up to equivalence$)$ to be the cone
$\Ca(d,\chi)$ of a chain $\epsilon$-bundle map $(d,\chi):(C_{k+1},0) \to (C_k,\delta)$
$($\ref{01Wu}$)$,
regarding $C_k$, $C_{k+1}$ as chain complexes concentrated in degree $k$.
The relative twisted $\epsilon$-quadratic $Q$-groups $Q_*(d,\chi,\epsilon)$ and the
absolute twisted $\epsilon$-quadratic $Q$-groups $Q_*(C,\gamma,\epsilon)$ are given as follows.\\
{\rm (i)} For $n \neq 2k-1,2k,2k+1,2k+2$
$$Q_n(C,\gamma,\epsilon)~=~Q_n(d,\chi,\epsilon)~=~Q_n(C,\epsilon)~=~
\begin{cases}
\widehat{Q}^{n+1}(C,\epsilon)&\text{\it if $n \geqslant 2k+3$}\\
0&\text{\it if $n\leqslant 2k-2$}
\end{cases}$$
with
$$\widehat{Q}^{n+1}(C,\epsilon)~=~
\dfrac{\{(\phi,\theta) \in S(C^{k+1})\oplus S(C^k)\,\vert\,
\phi=(-1)^{n+k}\epsilon\phi^*,d\phi d^*=\theta+(-1)^{n+k}\epsilon\theta^*\}}
{\{(\sigma+(-1)^{n+k}\epsilon\sigma^*,d\sigma d^*+\tau+(-1)^{n+k+1}\epsilon\tau^*\,\vert\,
(\sigma,\tau) \in S(C^{k+1})\oplus S(C^k)\}}$$
as given by Proposition \ref{hypercor}.\\
{\rm (ii)} For $n=2k-1,2k,2k+1,2k+2$
the relative twisted $\epsilon$-quadratic $Q$-groups are given by
$$\begin{array}{l}
Q_n(d,\chi,\epsilon)~=\\[2ex]
\begin{cases}
\dfrac{\{(\phi,\theta) \in S(C^{k+1})\oplus S(C^k)\,\vert\,
\phi=(-1)^k\epsilon\phi^*,d\phi d^*=\theta+(-1)^k\epsilon\theta^*\}}
{\{(\sigma+(-1)^k\epsilon\sigma^*,d\sigma d^*+\tau+(-1)^{k+1}\epsilon\tau^*\,\vert\,
(\sigma,\tau) \in S(C^{k+1})\oplus S(C^k)\}}&\text{\it if $n=2k+2$}\\[3ex]
\dfrac{\{(\psi,\eta) \in S(C^{k+1})\oplus S(C^k)\,\vert\,
(d,\chi)_{\%}(\psi)=(0,\eta+(-1)^{k+1}\epsilon\eta^*)\}}
{\{(\sigma+(-1)^{k+1}\epsilon\sigma^*,d\sigma d^* + \tau+(-1)^k\epsilon\tau^*)\,\vert\,
(\sigma,\tau) \in S(C^{k+1})\oplus S(C^k)\}}&\text{\it if $n=2k+1$}\\[4ex]
{\rm coker}((d,\chi)_{\%}:Q_{2k}(C_{k+1}) \to Q_{2k}(C_k,\delta))&
\text{\it if $n=2k$}\\[2ex]
Q_{2k-1}(C_k,\delta,\epsilon)&\text{\it if $n=2k-1$}
\end{cases}
\end{array}$$
with
$$\begin{array}{l}
Q_{2k}(C_k,\delta,\epsilon)~=~\dfrac{
\{(\phi,\theta) \in S(C^k) \oplus S(C^k)\,\vert\, \phi=(-1)^k\epsilon\phi^*,
\phi-\phi\delta\phi^*=\theta+(-1)^k\epsilon\theta^*\}}
{\{(0,\eta+(-1)^{k+1}\epsilon\eta^*)\,\vert\, \eta \in S(C^k)\}}~,\\[3ex]
Q_{2k-1}(C_k,\delta,\epsilon)~=~\dfrac{
\{\sigma  \in S(C^k)\,\vert\, \sigma=(-1)^k\epsilon\sigma^*\}}
{\{\phi-\phi\delta\phi^*-(\theta+(-1)^k\epsilon\theta^*)
\,\vert\,\phi=(-1)^k\epsilon\phi^*,\theta \in S(C^k)\}}~,\\[3ex]
(d,\chi)_{\%}~:~Q_{2k}(C_{k+1},\epsilon)~=~H_0(\Z_2;S(C^{k+1}),
(-1)^kT_{\epsilon}) \to
Q_{2k}(C_k,\delta,\epsilon)~;\\[1ex]
\hphantom{(d,\chi)_{\%}~:~} \psi \mapsto (d(\psi+(-1)^k\epsilon\psi^*)d^*,
d\psi d^* -
d(\psi+(-1)^k\epsilon\psi^*) \chi (\psi^*+(-1)^k\epsilon\psi)d^*)~.
\end{array}$$
The absolute twisted quadratic $Q$-groups are such that
$$Q_{2k-1}(C,\gamma,\epsilon)~=~Q_{2k-1}(d,\chi,\epsilon)~=~Q_{2k-1}(C_k,\delta,\epsilon)$$
and there is defined an exact sequence
$$\begin{array}{l}
0 \to Q_{2k+2}(d,\chi,\epsilon)  \xymatrix{\ar[r]^{t}&} Q_{2k+2}(C,\gamma,\epsilon) \\[1ex]
\hskip20pt \to H_{k+1}(C) \otimes_AC_{k+1} \xymatrix{\ar[r]^{F}&}
Q_{2k+1}(d,\chi,\epsilon)  \xymatrix{\ar[r]^{t}&}
Q_{2k+1}(C,\gamma,\epsilon)\\[1ex]
\hskip20pt \to H_k(C) \otimes_AC_{k+1} \xymatrix{\ar[r]^{F}&}
Q_{2k}(d,\chi,\epsilon)  \xymatrix{\ar[r]^{t}&}
Q_{2k}(C,\gamma,\epsilon) \to 0
\end{array}$$
with
$$\begin{array}{l}
F~:~H_k(C)\otimes_A C_{k+1}~=~{\rm coker}(d^*:{\rm Hom}_A(C^{k+1},C_{k+1})
\to {\rm Hom}_A(C^k,C_{k+1})) \\[1ex]
\hskip5pt \to Q_{2k}(d,\chi)~;~
\lambda \mapsto (\lambda d^*+(-1)^k\epsilon
d \lambda^*-d\lambda^*\delta \lambda d^*,\\[1ex]
\hskip5pt \lambda d^*-\lambda \chi \lambda^* - d \lambda^* \delta \lambda \chi
\lambda^* \delta \lambda d^* - d\lambda^*\delta(\lambda d^* +(-1)^k\epsilon
d \lambda^*) - (\lambda d^* +(-1)^k\epsilon d \lambda^*)\delta
d \lambda^* \delta\lambda d^*)~.
\end{array}$$
\end{Proposition}
\begin{proof}  The absolute and relative twisted $\epsilon$-quadratic
$Q$-groups are related by the exact sequence of \ref{twquThom}
$$\dots \to Q_n(d,\chi,\epsilon) \xymatrix{\ar[r]^{t}&} Q_n(C,\gamma,\epsilon) \to
H_{n-k-1}(C) \otimes_AC_{k+1} \xymatrix{\ar[r]^{F}&}
Q_{n-1}(d,\chi,\epsilon) \to \dots~.$$
The twisted $\epsilon$-quadratic $Q$-groups of
$(C_{k+1},0)$ are given by Proposition \ref{quadQ}
$$\begin{array}{ll}
Q_n(C_{k+1},0,\epsilon)~=~Q_n(C_{k+1},\epsilon)
&=~H_{n-2k}(\Z_2;S(C^{k+1}),(-1)^kT_{\epsilon})\\[1ex]
&=~\begin{cases}
\widehat{H}^{n-2k+1}(\Z_2;S(C^{k+1}),(-1)^kT_{\epsilon})&\hbox{\rm if $n\geqslant 2k+1$}\\
H_0(\Z_2;S(C^{k+1}),(-1)^kT_{\epsilon})&\hbox{\rm if $n=2k$}\\
0&\hbox{\rm if $n\leqslant 2k-1$}~.
\end{cases}
\end{array}$$
The twisted $\epsilon$-quadratic $Q$-groups of $(C_k,\delta)$ are given by Proposition \ref{0dim}
$$\begin{array}{l}
Q_n(C_k,\delta,\epsilon)\\[2ex]
=~\begin{cases}
\widehat{H}^{n-2k+1}(\Z_2;S(C^k),(-1)^kT_{\epsilon})
&\hbox{\rm if $n \geqslant 2k+1$}\\[2ex]
\dfrac{
\{(\phi,\theta) \in S(C^k) \oplus S(C^k)\,\vert\, \phi=(-1)^k\epsilon\phi^*,
\phi-\phi\delta\phi^*=\theta+(-1)^k\epsilon\theta^*\}}
{\{(0,\eta+(-1)^{k+1}\epsilon\eta^*)\,\vert\, \eta \in S(C^k)\}}&
\hbox{\rm if $n=2k$}\\[2ex]
\dfrac{
\{\sigma  \in S(C^k)\,\vert\, \sigma=(-1)^k\epsilon\sigma^*\}}
{\{\phi-\phi\delta\phi^*-(\theta+(-1)^k\epsilon\theta^*)
\,\vert\,\phi=(-1)^k\epsilon\phi^*,\theta \in S(C^k)\}}&\hbox{\rm if $n=2k-1$}\\[2ex]
0&\hbox{\rm if $n \leqslant 2k-2$}~.
\end{cases}
\end{array}$$
The twisted $\epsilon$-quadratic $Q$-groups of $(d,\chi)$ fit into the exact sequence
$$ \dots \longrightarrow Q_n (C_{k+1},\epsilon)
   \stackrel{(d,\chi)_\%}{\longrightarrow}
   Q_n (C_k, \delta,\epsilon) \longrightarrow Q_n (d,\chi,\epsilon)
   \longrightarrow Q_{n-1} (C_{k+1},\epsilon)
   \longrightarrow \dots$$
giving the expressions in the statements of (i) and (ii).
\end{proof}

\subsection{The Computation of $Q_*(C(X),\gamma(X))$}

In this we compute the twisted quadratic $Q$-groups
$Q_*(C(X),\gamma(X))$ of the following chain bundles
over an even commutative ring $A$.

\begin{Definition} \label{cx} {\rm
For $X \in {\rm Sym}_r(A)$ let
$$(C(X),\gamma(X))~=~\Ca(d,\chi)$$
be the cone of the chain bundle map over $A$
$$(d,\chi)~:~(C(X)_1,0) \to (C(X)_0,\delta)$$
defined by
$$\begin{array}{l}
d~=~2~:~C(X)_1~=~A^r \to C(X)_0~=~A^r~,\\[1ex]
\delta~=~X~:~C(X)_0~=~A^r \to C(X)^0~=~A^r~,\\[1ex]
\chi~=~2X~:~C(X)_1~=~A^r \to C(X)^1~=~A^r~.
\end{array}$$
\hfill$\qed$}
\end{Definition}

By Proposition \ref{01Wu} every chain bundle $(C,\gamma)$ with
$C_1=A^r \xymatrix{\ar[r]^2&} C_0=A^r$
is of the form $(C(X),\gamma(X))$ for some $X=(x_{ij}) \in {\rm Sym}_r(A)$, with
the equivalence class given by
$$\begin{array}{l}
\gamma~=~\gamma(X)~=~X~=~(x_{11},x_{22},\dots,x_{rr})\\[1ex]
\hskip50pt \in \widehat{Q}^0(C(X)^{-*})~=~
\dfrac{{\rm Sym}_r(A)}{{\rm Quad}_r(A)}~=~
\bigoplus\limits_r \widehat{H}^0(\Z_2;A)~(\ref{expl.relhyperQ})~.
\end{array}
$$
The 0th Wu class of $(C(X),\gamma(X))$ is the $A$-module morphism
$$\begin{array}{ll}
\widehat{v}_0(\gamma(X))~:&H_0(C(X))~=~(A_2)^r \to \widehat{H}^0(\Z_2;A)~;\\[1ex]
&a~=~(a_1,a_2,\dots,a_r) \mapsto aXa^t~=~\sum\limits^r_{i=1}a_ix_{ij}a_j~
=~\sum\limits^r_{i=1}(a_i)^2x_{ii}~.
\end{array}$$
In Theorem \ref{equ.ba-ba0} below the universal chain bundle
$(B^A,\beta^A)$ of a commutative even ring $A$ with $\widehat{H}^0(\Z_2;A)$
a f.g.  free $A_2$-module will be constructed from $(C(X),\gamma(X))$
for a diagonal $X \in {\rm Sym}_r(A)$ with $\widehat{v}_0(\gamma(X))$ an
isomorphism, and the twisted quadratic $Q$-groups $Q_*(B^A,\beta^A)$
will be computed using the following computation of
$Q_*(C(X),\gamma(X))$ (which holds for arbitrary $X$).

\begin{Theorem} \label{expl.Q0fchi}
Let $A$ be an even commutative ring, and let $X \in {\rm Sym}_r(A)$.\\
{\rm (i)} The twisted quadratic $Q$-groups of $(C(X),\gamma(X))$ are given by
$$\begin{array}{l}
Q_n(C(X),\gamma(X))~=\\[2ex]
\begin{cases}
0&\text{if $n \leqslant -2$}\\[2ex]
\dfrac{{\rm Sym}_r(A)}
{{\rm Quad}_r(A) + \{M-MXM\,\vert\,M \in {\rm Sym}_r(A)\}}
&\text{if $n=-1$}\\[2ex]
\dfrac{\{M \in {\rm Sym}_r(A)\,\vert\, M-MXM \in {\rm Quad}_r(A)\}}
{4{\rm Quad}_r(A) + \{2(N+N^t)-4N^tXN\,\vert\,N \in M_r(A)\}}
&\text{if $n=0$}\\[2ex]
\dfrac{\{N \in M_r(A)\,\vert\,N+N^t \in 2{\rm Sym}_r(A),
\dfrac{1}{2}(N+N^t)-N^tXN \in {\rm Quad}_r(A)\}}{2M_r(A)}&\text{if $n=1$}\\[2ex]
\dfrac{{\rm Sym}_r(A)}{{\rm Quad}_r(A)}&\text{if $n\geqslant 2$}~.
\end{cases}
\end{array}$$
{\rm (ii)} The boundary maps $\partial:Q_n(C(X),\gamma(X)) \to L_{n-1}(A)$
are given by
$$\begin{array}{l}
\partial~:~Q_{-1}(C(X),\gamma(X)) \to L_{-2}(A)~;~M \mapsto
(A^r \oplus (A^r)^*,\begin{pmatrix} M & 1 \\ 0 & X \end{pmatrix})~,\\[2ex]
\partial~:~Q_0(C(X),\gamma(X)) \to L_{-1}(A)~;~
M \mapsto (H_{-}(A^r);A^r,{\rm im}(\begin{pmatrix}
1-XM \\ M \end{pmatrix}:A^r \to A^r\oplus (A^r)^*))~,\\[2ex]
\partial~:~Q_1(C(X),\gamma(X)) \to L_0(A)~;~N \mapsto
(A^r \oplus A^r, \begin{pmatrix}
\dfrac{1}{4}(N+N^t-2N^tXN) & 1-2NX \\
0 & -2X \end{pmatrix})~.
\end{array}$$
{\rm (iii)} The twisted quadratic $Q$-groups of the chain bundles
$$(B(i),\beta(i))~=~(C(X),\gamma(X))_{*+2i}~~(i \in \Z)$$
are just the twisted quadratic $Q$-groups of $(C(X),\gamma(X))$ with a
dimension shift
$$Q_n(B(i),\beta(i))~=~Q_{n-4i}(C(X),\gamma(X))~.$$
\end{Theorem}
\begin{proof} (i)
Proposition \ref{twqcor} (i) and Example \ref{expl.relhyperQ} (ii) give
$$Q_n(C(X),\gamma(X))~=~\begin{cases}
0&\text{\rm if $n \leqslant -2$}\\
\widehat{Q}^{n+1}(C(X))=\dfrac{{\rm Sym}_r(A)}{{\rm Quad}_r(A)}
&\text{\rm if $n \geqslant 3$}~.
\end{cases}$$
For $-1 \leqslant n \leqslant 2$
Examples \ref{expl.relhyperQ}, \ref{expl.relsymmQ}, \ref{0dim.twq} and
Proposition \ref{twqcor} (ii) show that the commutative diagram with exact
rows and columns

$$\xymatrix{
Q^1(C(X)_1) \ar[d]^-{\displaystyle{d^{\%}}} \ar[r]^-{\displaystyle{J}}
& \widehat{Q}^1(C(X)_1) \ar[r]^-{\displaystyle{H}}
\ar[d]^-{\displaystyle{\widehat{d}^{\%}}}
& Q_0(C(X)_1)  \ar[r]^-{\displaystyle{1+T}} \ar[d]^-{\displaystyle{(d,\chi)}_{\%}}&
Q^0(C(X)_1)  \ar[r]^-{\displaystyle{J}} \ar[d]^-{\displaystyle{d^{\%}}}
& \widehat{Q}^0(C(X)_1) \ar[d]^-{\displaystyle{ \widehat{d}^{\%} }} \\
Q^1(C(X)_0) \ar[d] \ar[r]^-{\displaystyle{ J_{\delta} }}
& \widehat{Q}^1(C(X)_0) \ar[r]^-{\displaystyle{H_{\delta}}}  \ar[d] &
Q_0(C(X)_0,\delta)  \ar[r]^-{\displaystyle{N_{\delta}}}   \ar[d]&
Q^0(C(X)_0) \ar[r]^-{\displaystyle{J_\delta}} \ar[d]
& \widehat{Q}^0(C(X)_0) \ar[d]  \\
Q^1(d) \ar[r]^-{\displaystyle{J_\chi}} \ar[d] & \widehat{Q}^1(d)
\ar[r]^-{\displaystyle{H_{\chi}}} \ar[d]
& Q_0(d,\chi)  \ar[r]^-{\displaystyle{N_\chi}}\ar[d]& Q^0(d)  \ar[r]^-{\displaystyle{J_{\chi}}} \ar[d]
& \widehat{Q}^0(d) \ar[d] \\
Q^0(C(X)_1) \ar[r]^-{\displaystyle{J}} & \widehat{Q}^0(C(X)_1)
\ar[r]^-{\displaystyle{H}} &
Q_{-1}(C(X)_1)  \ar[r]^-{\displaystyle{1+T}} &
Q^{-1}(C(X)_1)  \ar[r]^-{\displaystyle{J}} & \widehat{Q}^{-1}(C(X)_1) }$$

\noindent is given by

$$\xymatrix@C-10pt{
0 \ar[d] \ar[r] & 0 \ar[r] \ar[d]
& {\rm Quad}_r(A)  \ar[r] \ar[d]^-{\displaystyle{4}}&
{\rm Sym}_r(A) \ar[r] \ar[d]^-{\displaystyle{4}}
&  \dfrac{{\rm Sym}_r(A)}{{\rm Quad}_r(A)} \ar[d]^-{\displaystyle{0}} \\
0 \ar[d] \ar[r] & 0  \ar[r] \ar[d] & Q_0(C(X)_0,\delta)  \ar[r]
\ar[d]& {\rm Sym}_r(A) \ar[r]^-{\displaystyle{J_\delta}} \ar[d]
& \dfrac{{\rm Sym}_r(A)}{{\rm Quad}_r(A)} \ar[d]^-{\displaystyle{1}}  \\
0 \ar[r] \ar[d] & \dfrac{{\rm Sym}_r(A)}{{\rm Quad}_r(A)}
\ar[r]^-{\displaystyle{4}} \ar[d]^-{\displaystyle{1}}
& Q_0(d,\chi)  \ar[r]\ar[d]& \dfrac{{\rm Sym}_r(A)}{4{\rm Sym}_r(A)}
\ar[r] \ar[d]
& \dfrac{{\rm Sym}_r(A)}{{\rm Quad}_r(A)} \ar[d] \\
{\rm Sym}_r(A) \ar[r] & \dfrac{{\rm Sym}_r(A)}{{\rm Quad}_r(A)} \ar[r] & 0  \ar[r] &
0  \ar[r] & 0 }$$

\noindent with
$$\begin{array}{l}
J_{\delta}~:~{\rm Sym}_r(A) \to \dfrac{{\rm Sym}_r(A)}{{\rm Quad}_r(A)}~;~
M \mapsto M-MXM~,\\[2ex]
Q_0(C(X)_0,\delta)~=~{\rm ker}(J_{\delta})~=~
\{M \in {\rm Sym}_r(A)\,\vert\,M-MXM \in {\rm Quad}_r(A)\}~,\\[2ex]
Q_0(d,\chi)~=~{\rm coker}((d,\chi)_{\%}:Q_0(C(X)_1) \to Q_0(C(X)_0,\delta))\\[2ex]
\hphantom{Q_0(d,\chi)~}=~\dfrac{\{M \in {\rm Sym}_r(A)\,\vert\,M-MXM \in {\rm Quad}_r(A)\}}
{4{\rm Quad}_r(A)}~,\\[2ex]
N_{\chi}~:~Q_0(d,\chi) \to Q^0(d)~=~\dfrac{{\rm Sym}_r(A)}{4{\rm Sym}_r(A)}~;~
M \mapsto M~.
\end{array}$$
\noindent Furthermore, the commutative braid of exact sequences
\vskip4mm

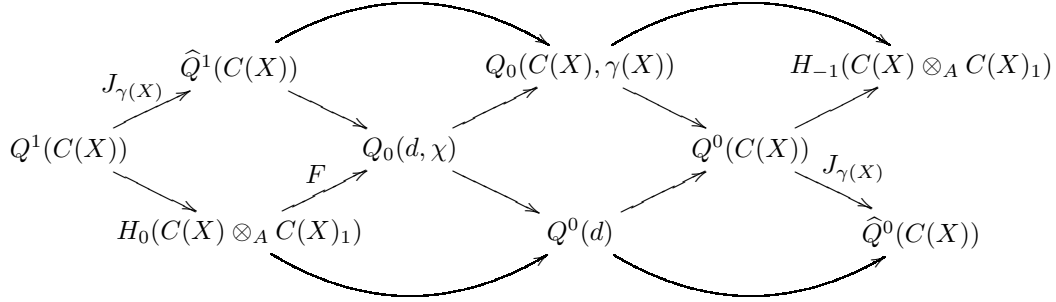
$$\xymatrix@!C@C-65pt@R-10pt{
&\widehat{Q}^1(C(X)) \ar[dr]\ar@/^2pc/[rr]&& Q_0(C(X),\gamma(X)) \ar[dr]\ar@/^2pc/[rr]&&
H_{-1}(C(X)\otimes_AC(X)_1)\\
Q^1(C(X)) \ar[ur]^-{\displaystyle{J_{\gamma(X)}}} \ar[dr] && Q_0(d,\chi)
\ar[ur]\ar[dr]&& Q^0(C(X))\ar[ur]\ar[dr]^-{\displaystyle{J_{\gamma(X)}}}&\\
&H_0(C(X)\otimes_AC(X)_1)\ar[ur]^{\displaystyle{F}}\ar@/_2pc/[rr] &&Q^0(d)\ar[ur]\ar@/_2pc/[rr] &&
\widehat{Q}^0(C(X))}$$

\vskip4mm

\noindent is given by

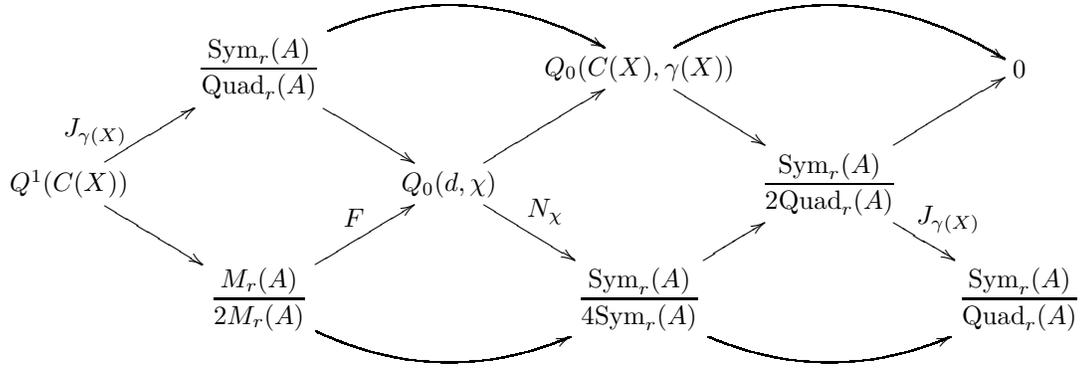
$$\xymatrix@!C@C-30pt@R-10pt{
&\dfrac{{\rm Sym}_r(A)}{{\rm Quad}_r(A)}
\ar[dr]\ar@/^2pc/[rr]&& Q_0(C(X),\gamma(X)) \ar[dr]\ar@/^2pc/[rr]&&0\\
Q^1(C(X)) \ar[ur]^-<<{\displaystyle{J_{\gamma(X)}}} \ar[dr] && Q_0(d,\chi)
\ar[ur]\ar[dr]^-{\displaystyle{N_{\chi}}}&&
\dfrac{{\rm Sym}_r(A)}{2{\rm Quad}_r(A)}
\ar[ur]\ar[dr]^-{\displaystyle{J_{\gamma(X)}}}&\\
&\dfrac{M_r(A)}{2M_r(A)} \ar[ur]^-{\displaystyle{F}}\ar@/_2pc/[rr] &&
\dfrac{{\rm Sym}_r(A)}{4{\rm Sym}_r(A)} \ar[ur]\ar@/_2pc/[rr] &&
\dfrac{{\rm Sym}_r(A)}{{\rm Quad}_r(A)}}$$

\noindent with
$$\begin{array}{l}
\dfrac{{\rm Sym}_r(A)}{2{\rm Quad}_r(A)}~\cong~Q^0(C(X))~;~M \mapsto \phi~
({\rm where}~\phi_0=M:C^0\to C(X)_0)~,\\[2ex]
J_{\gamma(X)}~:~\dfrac{{\rm Sym}_r(A)}{2{\rm Quad}_r(A)} \to
\dfrac{{\rm Sym}_r(A)}{{\rm Quad}_r(A)}~;~M \mapsto M-MXM~,\\[2ex]
F~:~H_0(C(X)\otimes_AC(X)_1)~=~\dfrac{M_r(A)}{2M_r(A)} \to Q_0(d,\chi)~;~N \mapsto
2(N+N^t)-4N^tXN~,\\[2ex]
Q^1(C(X))~=~{\rm ker}(N_{\chi}F:H_0(C(X)\otimes_AC(X)_1) \to Q^0(d))\\[2ex]
\hphantom{Q^1(C(X))~}=~\dfrac{\{N \in M_r(A)\,\vert\,N+N^t \in 2{\rm Sym}_r(A)\}}
{2M_r(A)}\\[2ex]
\hphantom{Q^1(C(X))~=~}({\rm where}~\phi \in Q^1(C(X))~{\rm corresponds ~to}~N=\phi_0 \in M_r(A))~,
\\[2ex]
J_{\gamma(X)}~:~Q^1(C(X)) \to \widehat{Q}^1(C(X))~=~
\dfrac{{\rm Sym}_r(A)}{{\rm Quad}_r(A)}~;~N \mapsto
\displaystyle{\frac{1}{2}}(N+N^t) - N^tXN~.
\end{array}$$
It follows that
$$\begin{array}{l}
Q_0(C(X),\gamma(X))~=~{\rm coker}(
F:\dfrac{M_r(A)}{2M_r(A)} \to Q_0(d,\chi))\\[2ex]
\hphantom{Q_0(C(X),\gamma(X))~}
=~\dfrac{\{M \in {\rm Sym}_r(A)\,\vert\, M-MXM \in {\rm Quad}_r(A)\}}
{4{\rm Quad}_r(A) + \{2(N+N^t)-4N^tXN\,\vert\,N \in M_r(A)\}}~,\\[3ex]
Q_{-1}(C(X),\gamma(X))~=~Q_{-1}(d,\chi)\\[2ex]
\hphantom{Q_{-1}(C(X),\gamma(X))~}
=~{\rm coker}(J_{\gamma(X)}:Q^0(C(X)) \to \widehat{Q}^0(C(X)))\\[2ex]
\hphantom{Q_{-1}(C(X),\gamma(X))~}
=~\displaystyle{
\frac{{\rm Sym}_r(A)}{{\rm Quad}_r(A)+ \{ M-MXM\,\vert\,M \in {\rm Sym}_r(A)\}}}
\end{array}$$
with
$$\begin{array}{l}
\widehat{Q}^1(C(X))~=~\dfrac{{\rm Sym}_r(A)}{{\rm Quad}_r(A)}
\to Q_0(C(X),\gamma(X))~;~M \mapsto 4M~,\\[2ex]
\widehat{Q}^0(C(X))~=~\dfrac{{\rm Sym}_r(A)}{{\rm Quad}_r(A)}
\to Q_{-1}(C(X),\gamma(X))~;~M \mapsto M~.
\end{array}$$
\noindent Also
$$\begin{array}{l}
(d,\chi)_{\%}~=~0~:~Q_2(d,\chi)~=~
Q_1(C(X)_1)~=~\dfrac{{\rm Sym}_r(A)}{{\rm Quad}_r(A)}\\[2ex]
\hskip150pt
\to Q_1(C(X)_0,\delta)~=~\widehat{Q}^2(C(X)_0)~=~
\dfrac{{\rm Sym}_r(A)}{{\rm Quad}_r(A)}~,\\[2ex]
Q^1(C(X))~=~{\rm ker}(N_{\chi}F:H_0(C(X)\otimes_AC(X)_1) \to Q^0(d))\\[2ex]
\hphantom{Q^1(C(X))~}=~\dfrac{\{N \in M_r(A)\,\vert\,N+N^t \in 2{\rm Sym}_r(A)\}}
{2M_r(A)}~,\\[2ex]
J_{\gamma(X)}~:~Q^1(C(X))
\to \widehat{Q}^1(C(X))~=~\dfrac{{\rm Sym}_r(A)}{{\rm Quad}_r(A)}~;~
N \mapsto \dfrac{1}{2}(N+N^t)-N^tXN~,\\[2ex]
Q_1(C(X),\gamma(X))~=~{\rm ker}(J_{\gamma(X)}:Q^1(C(X)) \to \widehat{Q}^1(C(X)))\\[2ex]
\hskip20pt
=~\dfrac{\{N \in M_r(A)\,\vert\,N+N^t \in 2{\rm Sym}_r(A),
\dfrac{1}{2}(N+N^t)-N^tXN \in {\rm Quad}_r(A)\}}{2M_r(A)}~,\\[2ex]
Q_2(C(X),\gamma(X))~=~Q_2(d,\chi)~=~\dfrac{{\rm Sym}_r(A)}{{\rm Quad}_r(A)}~.
\end{array}$$
(ii) The expressions for $\partial:Q_n(C(X),\gamma(X)) \to L_{n-1}(A)$
are given by the boundary construction of Proposition
\ref{norminv} and its expression in terms of forms and
formations (\ref{bformation}, \ref{bform}).
The form in the case $n=-1$ (resp. the formation in the case $n=0$)
is given by \ref{bform} (resp. \ref{bformation}) applied to the
$n$-dimensional symmetric structure $(\phi,\theta) \in Q_n(C(X),\gamma(X))$
corresponding to $M \in {\rm Sym}_r(A)$. For $n=1$
the boundary of the $1$-dimensional symmetric structure
$(\phi,\theta) \in Q_1(C(X),\gamma(X))$ corresponding to $N \in
M_r(A)$ with
$$N+N^t \in 2{\rm Sym}_r(A)~,~\dfrac{1}{2}(N+N^t) - N^tXN \in
{\rm Quad}_r(A)$$
is a 0-dimensional quadratic Poincar\'e complex $(C,\psi)$ with
$$C~=~\Ca(N:C(X)^{1-*} \to C(X))_{*+1}~.$$
The instant surgery obstruction (\ref{bform}) is the nonsingular
quadratic form
$$\begin{array}{l}
I(C,\psi)~=\\[1ex]
({\rm coker}( \begin{pmatrix} -2 \\ N^t \\ 1+2XN^t \end{pmatrix}
:A^r \to A^r \oplus A^r \oplus A^r),
\begin{pmatrix} \dfrac{1}{4}(N+N^t-2NXN^t) & 1 & N \\
0 & -2X & 2 \\ 0 & 0 & 0
\end{pmatrix})
\end{array}
$$
such that there is defined an isomorphism
$$\begin{pmatrix} 1 & -4X & 2 \\ N^t & 1-2N^t X & N^t \end{pmatrix}~:~
I(C,\psi) \to (A^r \oplus A^r, \begin{pmatrix}
\dfrac{1}{4}(N+N^t-2N^tXN) & 1-2NX \\
0 & -2X \end{pmatrix})~.$$
(iii) The even multiple skew-suspension isomorphisms of
the symmetric $Q$-groups
$$\overline{S}^{2i}~:~Q^{n-4i}(C(X)_{*+2i}) \xymatrix{\ar[r]^-{\cong}&}
Q^n(C(X))~;~
\{\phi_s\,\vert\,s \geqslant 0\} \mapsto \{\phi_s\,\vert\,s \geqslant 0\}~~(i \in \Z)$$
are defined also for
the hyperquadratic, quadratic and twisted quadratic $Q$-groups.
\end{proof}

\subsection{The Universal Chain Bundle}

For any $A$-module chain complexes $B,C$ the additive group
$H_0({\rm Hom}_A(C,B))$ consists of the
chain homotopy classes of $A$-module chain maps
$f:C \to B$. For a chain $\epsilon$-bundle $(B,\beta)$ there is thus defined a
morphism
$$H_0({\rm Hom}_A(C,B)) \to \widehat{Q}^0(C^{0-*},\epsilon)~;~(f:C \to B)
\mapsto \widehat{f}^*(\beta)~.$$

\begin{Proposition} {\rm (Weiss \cite{weiss1})}
{\rm (i)} For every ring
with involution $A$ and $\epsilon=\pm 1$ there exists a universal chain $\epsilon$-bundle
$(B^{A,\epsilon},\beta^{A,\epsilon})$ over $A$ such that for any finite f.g. projective
$A$-module chain complex $C$ the morphism
$$H_0({\rm Hom}_A(C,B^{A,\epsilon})) \to \widehat{Q}^0(C^{0-*},\epsilon)~;~(f:C \to B^{A,\epsilon})
\mapsto \widehat{f}^*(\beta^{A,\epsilon})$$
is an isomorphism. Thus every chain $\epsilon$-bundle $(C,\gamma)$ is classified
by a chain $\epsilon$-bundle map
$$(f,\chi)~:~(C,\gamma) \to (B^{A,\epsilon},\beta^{A,\epsilon})~.$$
{\rm (ii)} The universal chain $\epsilon$-bundle $(B^{A,\epsilon}, \beta^{A,\epsilon})$ is
characterized {\rm (}uniquely up to equivalence{\rm )}
by the property that its Wu classes are $A$-module isomorphisms
$$ \widehat{v}_k (\beta^{A,\epsilon})~:~ H_k (B^{A,\epsilon}) \stackrel{\cong}{\longrightarrow}
  \widehat{H}^k (\Z_2;A,\epsilon)~~(k \in \Z)~. $$
{\rm (iii)} An $n$-dimensional
{\rm (}$\epsilon$-symmetric, $\epsilon$-quadratic{\rm )} Poincar\'e pair
over $A$ has a canonical universal $\epsilon$-bundle
$(B^{A,\epsilon},\beta^{A,\epsilon})$-structure.\\
{\rm (iv)}  The 4-periodic $(B^{A,\epsilon},\beta^{A,\epsilon})$-structure
$L$-groups are the 4-periodic versions of the
$\epsilon$-symmetric and $\epsilon$-hyperquadratic $L$-groups of
$A$
$$\begin{array}{l}
L\langle B^{A,\epsilon}, \beta^{A,\epsilon} \rangle^{n+4*}
(A,\epsilon)~=~L^{n+4*} (A,\epsilon)~,\\[1ex]
\widehat{L}\langle B^{A,\epsilon},
\beta^{A,\epsilon} \rangle^{n+4*} (A,\epsilon)~=~
\widehat{L}^{n+4*}(A,\epsilon)~.
\end{array}$$
{\rm (v)} The twisted $\epsilon$-quadratic $Q$-groups of
$(B^{A,\epsilon},\beta^{A,\epsilon})$ fit into an exact sequence
$$\dots \rightarrow L_n (A,\epsilon) \xymatrix{\ar[r]^-{1+T_{\epsilon}}&}
L^{n+4*}(A,\epsilon)
\rightarrow Q_n(B^{A,\epsilon},\beta^{A,\epsilon},\epsilon)
\stackrel{\partial}{\rightarrow} L_{n-1} (A,\epsilon) \rightarrow
\dots$$
with
$$\partial~:~Q_n (B^{A,\epsilon},\beta^{A,\epsilon},\epsilon) \to
L_{n-1} (A,\epsilon)~;~(\phi,\theta) \mapsto (C,\psi)$$
given by the construction of Proposition \ref{norm} (ii), with
$$C~=~\Ca(\phi_0:(B^{A,\epsilon})^{n-*} \to B^{A,\epsilon})_{*+1}~
{\it  etc.}$$
\hfill$\qed$
\end{Proposition}

For $\epsilon=1$ write
$$(B^{A,1},\beta^{A,1})~=~(B^A,\beta^A)$$
and note that
$$(B^{A,-1},\beta^{A,-1})~=~(B^A,\beta^A)_{*-1}~.$$

In general, the chain $A$-modules $B^{A,\epsilon}$ are not finitely
generated, although $B^{A,\epsilon}$ is a direct limit of f.g.  free
$A$-module chain complexes.  In our applications the involution on $A$
will satisfy the following conditions~:

\begin{Proposition} \label{cr}
{\rm (Connolly and Ranicki \cite[Section 3.6]{connran})}\\
Let $A$ be a ring with an even involution such that $\widehat{H}^0(\Z_2;A)$
has a 1-dimensional f.g. projective $A$-module resolution
$$0 \to C_1 \xymatrix{\ar[r]^d&} C_0  \xymatrix{\ar[r]^x&}
\widehat{H}^0(\Z_2;A) \to 0~.$$
Let $(C,\gamma)=\Ca(d,\chi)$ be the cone of a chain bundle map
$(d,\chi):(C_1,0) \to (C_0,\delta)$ with
$$\widehat{v}_0(\delta)~=~x~:~C_0 \to \widehat{H}^0(\Z_2;A)$$
and set
$$(B^A(i),\beta^A(i))~=~(C,\gamma)_{*+2i}~~(i \in \Z)~.$$
{\rm (i)} The chain bundle over $A$
$$(B^A,\beta^A)~=~\bigoplus\limits_i(B^A(i),\beta^A(i))$$
is universal.\\
{\rm (ii)} The twisted quadratic $Q$-groups of $(B^A,\beta^A)$ are given by
$$Q_n(B^A,\beta^A)~=~
\begin{cases}
Q_0(C,\gamma)&\hbox{\it if $n\equiv 0(\bmod\, 4)$}\\
{\rm ker}(J_{\gamma}:Q^1(C) \to \widehat{Q}^1(C))&\text{if $n\equiv 1(\bmod\, 4)$}\\
0&\hbox{\it if $n\equiv 2(\bmod\, 4)$}\\
Q_{-1}(C,\gamma)&\hbox{\it if $n\equiv 3(\bmod\, 4)$}~.
\end{cases}$$
The projection $(B^A,\beta^A) \to (B^A(2j),\beta^A(2j))$ induces isomorphisms
$$Q_n(B^A,\beta^A)~\cong~
\begin{cases} Q_n(B^A(2j),\beta^A(2j))&\text{\it if $n=4j,4j-1$}\\
{\rm ker}(J_{\beta^A(2j)}:Q^n(B^A(2j)) \to \widehat{Q}^n(B^A(2j)))
&\text{\it if $n=4j+1$~.}
\end{cases} $$
\end{Proposition}
\begin{proof} (i)
The Wu classes of the chain bundle $(C,\gamma)_{*+2i}$ are isomorphisms
$$\widehat{v}_k(\gamma)~:~H_k(C_{*+2i})
\xymatrix{\ar[r]^-{\displaystyle{\cong}}&} \widehat{H}^k(\Z_2;A)$$
for $k=2i,2i+1$.\\
(ii) See \cite{connran} for the
detailed analysis of the exact sequence of \ref{seq} (ii)
$$ \begin{array}{l}
\dots \rightarrow \sum\limits_{i=-\infty}^\infty Q_n (B^A(i),\beta^A (i))
  \rightarrow Q_n (B^A,\beta^A) \rightarrow
  \sum\limits_{i<j} H_n (B^A(i)\otimes_AB^A(j))\\[2ex]
\hphantom{\dots \rightarrow \sum\limits_{i=-\infty}^\infty Q_n (B^A(i),\beta^A (i))
  \rightarrow Q_n (B^A,\beta^A)}
\rightarrow \sum\limits_{i=-\infty}^\infty Q_{n-1} (B^A(i),\beta^A (i))
   \rightarrow \dots~.
\end{array}$$
\end{proof}

As in the Introduction~:

\begin{Definition} {\rm A ring with involution $A$ is {\it $r$-even} for some
$r \geqslant 1$ if
\begin{itemize}
\item[(i)] $A$ is commutative, with the identity involution,
\item[(ii)] $2 \in A$ is a non-zero divisor,
\item[(iii)] $\widehat{H}^0(\Z_2;A)$ is a f.g. free $A_2$-module of rank $r$
with a basis $\{x_1=1,x_2,\dots,x_r\}$.
\end{itemize}}
\hfill$\qed$
\end{Definition}

\begin{Example} {\rm $\Z$ is 1-even.}
\hfill$\qed$
\end{Example}

\begin{Proposition} \label{2even}
If $A$ is 1-even the polynomial extension $A[x]$
is 2-even, with $A[x]_2=A_2[x]$ and $\{1,x\}$ an $A_2[x]$-module
basis of $\widehat{H}^0(\Z_2;A[x])$.
\end{Proposition}
\begin{proof} For any $a=\sum\limits^{\infty}_{i=0}a_ix^i \in A[x]$
$$\begin{array}{ll}
a^2&=~\sum\limits^{\infty}_{i=0}(a_i)^2x^{2i}+
2\sum\limits_{0 \leqslant i < j <\infty}a_ia_jx^{i+j}\\[1ex]
&=~\sum\limits^{\infty}_{i=0}a_ix^{2i} \in A_2[x]~.
\end{array}$$
The $A_2[x]$-module morphism
$$A_2[x] \oplus A_2[x] \to
\widehat{H}^0(\Z_2;A[x])~;~(p,q) \mapsto p^2+q^2x$$
is thus an isomorphism, with inverse
$$\widehat{H}^0(\Z_2;A[x]) \xymatrix{\ar[r]^-{\displaystyle{\cong}}&}
A_2[x] \oplus A_2[x]~;~a~=~\sum\limits^{\infty}_{i=0}a_ix^i
\mapsto (\sum\limits^{\infty}_{j=0}a_{2j}x^j,
\sum\limits^{\infty}_{j=0}a_{2j+1}x^j)~.$$
\end{proof}

Proposition \ref{2even} is the special case $k=1$ of a general result:
if $A$ is 1-even and $t_1,t_2,\dots,t_k$ are commuting
indeterminates over $A$ then the polynomial ring $A[t_1,t_2,\dots,t_k]$
is $2^k$-even with
$$\{x_1=1,x_2,x_3,\dots,x_{2^k}\}~=~\{(t_1)^{i_1}(t_2)^{i_2}\dots
(t_k)^{i_k}\,\vert\,i_j=0~{\rm or}~1,1 \leqslant j \leqslant k\}$$
an $A_2[t_1,t_2,\dots,t_k]$-module basis of
$\widehat{H}^0(\Z_2;A[t_1,t_2,\dots,t_k])$.

We can now prove Theorem \ref{main3}~:

\begin{Theorem} \label{equ.ba-ba0}
Let $A$ be an $r$-even ring with involution.\\
{\rm (i)} The $A$-module morphism
$$x~:~A^r \to \widehat{H}^0(\Z_2;A)~;~(a_1,a_2,\dots,a_r)\mapsto
\sum\limits^r_{i=1}(a_i)^2x_i$$
fits into a 1-dimensional f.g.  free $A$-module resolution of
$\widehat{H}^0(\Z_2;A)$
$$0 \to C_1~=~A^r \xymatrix{\ar[r]^2&} C_0~=~A^r  \xymatrix{\ar[r]^x&}
\widehat{H}^0(\Z_2;A) \to 0~.$$
The symmetric and hyperquadratic $L$-groups of $A$ are 4-periodic
$$L^n(A)~=~L^{n+4}(A)~,~\widehat{L}^n(A)~=~\widehat{L}^{n+4}(A)~.$$
{\rm (ii)} Let $(C(X),\gamma(X))$ be the chain bundle over $A$
given by the construction of {\rm (}\ref{cx}{\rm )} for
$$X~=~\begin{pmatrix} x_1 & 0 & 0 & \dots & 0 \\
0 & x_2 & 0 & \dots & 0 \\
0 & 0 & x_3 & \dots & 0 \\
\vdots & \vdots & \vdots & \ddots & \vdots \\
0 & 0 & 0 & \dots & x_r \end{pmatrix} \in {\rm Sym}_r(A)~,$$
with $C(X)=\Ca(2:A^r \to A^r)$.
The chain bundle over $A$ defined by
$$(B^A,\beta^A)~=~\bigoplus\limits_i(C(X),\gamma(X))_{*+2i}~=~
\bigoplus\limits_i(B^A(i),\beta^A(i))$$
is universal. The hyperquadratic $L$-groups of $A$ are given by
$$\begin{array}{l}
\widehat{L}^n(A)~=~Q_n(B^A,\beta^A)~=\\[2ex]
\begin{cases}
Q_0(C(X),\gamma(X))~=~\dfrac{\{M \in {\rm Sym}_r(A)\,\vert\, M-MXM \in {\rm Quad}_r(A)\}}
{4{\rm Quad}_r(A) + \{2(N+N^t)-N^tXN\,\vert\,N \in M_r(A)\}}&
\hbox{\it if $n=0$}\\[3ex]
{\rm im}(N_{\gamma(X)}:Q_1(C(X),\gamma(X)) \to Q^1(C(X)))~=~
{\rm ker}(J_{\gamma(X)}:Q^1(C(X)) \to \widehat{Q}^1(C(X)))&\\[2ex]
=~\dfrac{\{N \in M_r(A)\,\vert\,N+N^t \in 2{\rm Sym}_r(A),
\dfrac{1}{2}(N+N^t)-N^tXN \in {\rm Quad}_r(A)\}}{2M_r(A)}&\text{if $n=1$}\\[2ex]
0&\hbox{\it if $n=2$}\\[2ex]
Q_{-1}(C(X),\gamma(X))~=~\dfrac{{\rm Sym}_r(A)}
{{\rm Quad}_r(A) + \{L-LXL\,\vert\,L \in {\rm Sym}_r(A)\}}&
\hbox{\it if $n=3$}
\end{cases}
\end{array}$$
with
$$\begin{array}{l}
\partial~:~\widehat{L}^0(A) \to L_{-1}(A)~;~
M \mapsto (H_{-}(A^r);A^r,{\rm im}(\begin{pmatrix}
1-XM \\ M \end{pmatrix}:A^r \to A^r\oplus (A^r)^*))~,\\[2ex]
\partial~:~\widehat{L}^1(A) \to L_0(A)~;~N \mapsto
(A^r \oplus A^r, \begin{pmatrix}
\dfrac{1}{4}(N+N^t-2N^tXN) & 1-2NX \\
0 & -2X \end{pmatrix})~,\\[1ex]
\partial~:~\widehat{L}^3(A) \to L_2(A)~;~M \mapsto
(A^r \oplus (A^r)^*,\begin{pmatrix} M & 1 \\ 0 & X \end{pmatrix})~.
\end{array}$$
\end{Theorem}
\begin{proof} Combine Proposition \ref{4period},
Theorem \ref{expl.Q0fchi} and Proposition \ref{cr}.
\end{proof}

We can now prove Theorem \ref{main1}~:

\begin{Corollary} \label{cor1}
Let $A$ be a 1-even ring.\\
{\rm (i)} The universal chain bundle $(B^A,\beta^A)$ over $A$ is given by
$$\begin{array}{l}
\xymatrix@C-10pt{B^A~:~\dots \ar[r]&
B^A_{2k+2}=A\ar[r]^-{0}&B^A_{2k+1}=A\ar[r]^-{2}&
B^A_{2k}=A\ar[r]^-{0}&B^A_{2k-1}=A\ar[r]&\dots}~,\\[1ex]
(\beta^A)_{-4k}~=~1~:~B^A_{2k}~=~A \to (B^A)^{2k}~=~A~~(k \in \Z)~.
\end{array}$$
{\rm (ii)} The hyperquadratic $L$-groups of $A$ are given by
$$\widehat{L}^n(A)~=~Q_n(B^A,\beta^A)~=~
\begin{cases}
A_8&\hbox{\it if $n\equiv 0(\bmod\, 4)$}\\
A_2&\text{if $n\equiv 1,3(\bmod\, 4)$}\\
0&\hbox{\it if $n\equiv 2(\bmod\, 4)$}
\end{cases}$$
with
$$\begin{array}{l}
\partial~:~\widehat{L}^0(A)=A_8 \to L_{-1}(A)~;~
a \mapsto (H_-(A);A,{\rm im}(\begin{pmatrix} 1-a \\ a \end{pmatrix}:A
\to A \oplus A))~,\\[2ex]
\partial~:~\widehat{L}^1(A)=A_2 \to L_0(A)~;~
a \mapsto (A\oplus A,\begin{pmatrix} a(1-a)/2 & 1-2a  \\ 0 & -2 \end{pmatrix})~;\\[2ex]
\partial~:~\widehat{L}^3(A)=A_2 \to L_2(A)~;~
a \mapsto (A\oplus A,\begin{pmatrix} a & 1 \\ 0 & 1 \end{pmatrix})~.
\end{array}$$
{\rm (iii)} The map $L^0(A) \to \widehat{L}^0(A)$ sends the
Witt class $(K,\lambda) \in L^0(A)$ of a nonsingular symmetric form
$(K,\lambda)$ over $A$ to
$$[K,\lambda]~=~\lambda(v,v) \in \widehat{L}^0(A)~=~A_8$$
for any $v \in K$ such that
$$\lambda(x,x)~=~\lambda(x,v) \in A_2~~(x \in K)~.$$
\end{Corollary}
\begin{proof} (i)+(ii) The $A$-module morphism
$$\widehat{v}_0(\beta^A)~:~H_0(B^A)~=~A_2 \to \widehat{H}^0(\Z_2;A)~;~a \mapsto a^2$$
is an isomorphism. Apply Theorem \ref{equ.ba-ba0} with $r=1$, $x_1=1$.  \\
(ii) The computation of $\widehat{L}^*(A)=Q_*(B^A,\beta^A)$
is given by Theorem \ref{equ.ba-ba0}, using the fact that $a-a^2 \in 2A$
$(a \in A)$ for a 1-even $A$. The explicit descriptions of $\partial$
are special cases of the formulae in Theorem \ref{expl.Q0fchi} (ii).\\
(iii) As in Example \ref{Spiv} regard $(K,\lambda)$ as a 0-dimensional
symmetric Poincar\'e complex $(D,\phi)$ with
$$\phi_0~=~\epsilon \lambda^{-1}~:~D^0~=~K \to D^0~=~K^*~.$$
The Spivak normal chain bundle $\gamma=\lambda^{-1} \in \widehat{Q}^0(D^{0-*})$
is classified by the chain bundle map $(v,0):(D,\gamma) \to (B^A,\beta^A)$
with
$$g~:~D_0~=~K^* \to \widehat{H}^0(\Z_2;A)~;~x \mapsto \lambda^{-1}(x,x)~=~x(v)~.$$
The algebraic normal invariant $(\phi,0) \in Q_0(D,\gamma)$ has image
$$g_{\%}(\phi,0)~=~\lambda(v,v) \in Q_0(B^A,\beta^A)~=~A_8~.$$
\end{proof}

\begin{Example} \label{expl.QBZ}
{\rm For $R=\Z$
$$\widehat{L}^n(\Z)~=~Q_n(B^\Z,\beta^\Z)~=~\begin{cases}
\Z_8&\text{if $n\equiv 0(\bmod\, 4)$}\\
\Z_2&\text{if $n\equiv 1,3(\bmod\, 4)$}\\
0&\text{if $n\equiv 2(\bmod\, 4)$}
\end{cases}$$
as recalled (from \cite{ranicki1}) in the Introduction.\hfill$\qed$}
\end{Example}

\section{The Generalized Arf Invariant for Forms}

A nonsingular $\epsilon$-quadratic form $(K,\psi)$ over $A$ corresponds
to a 0-dimensional $\epsilon$-quadratic Poincar\'e complex over $A$.
The 0-dimensional $\epsilon$-quadratic $L$-group $L_0(A,\epsilon)$ is
the Witt group of nonsingular $\epsilon$-quadratic forms, and similarly
for $L^0(A,\epsilon)$ and $\epsilon$-symmetric forms.  In this section
we define the `generalized Arf invariant'
$$(K,\psi;L) \in Q_1(B^{A,\epsilon},\beta^{A,\epsilon})~=~
\widehat{L}^{4*+1}(A,\epsilon)$$
for a nonsingular $\epsilon$-quadratic form $(K,\psi)$ over $A$ with a
lagrangian $L$ for the $\epsilon$-symmetric form $(K,\psi+\epsilon\psi^*)$,
so that
$$\begin{array}{l}
(K,\psi)~=~\partial(K,\psi;L) \in
{\rm ker}(1+T:L_0(A,\epsilon) \to L^{4*}(A,\epsilon))\\[1ex]
\hphantom{(K,\psi)~=~\partial(K,\psi;L) \in}~=~
{\rm im}(\partial:Q_1(B^{A,\epsilon},\beta^{A,\epsilon},\epsilon) \to
L_0(A,\epsilon))~.
\end{array}$$

\subsection{Forms and Formations}

Given a f.g. projective $A$-module $K$ and the inclusion
$j:L \to K$ of a direct summand, let $f:C \to D$ be the chain map defined by
$$\begin{array}{l}
C~:~\dots \to 0 \to C_k~=~K^* \to 0 \to \dots~,\\[1ex]
D~:~\dots \to 0 \to D_k~=~L^* \to 0 \to \dots~,\\[1ex]
f~=~j^*~:~C_k~=~K^* \to D_k~=~L^*~.
\end{array}$$
The symmetric $Q$-group
$$Q^{2k}(C)~=~H^0(\Z_2;S(K),(-1)^kT)~=~\{\phi \in S(K)\,\vert\, \phi^*=(-1)^k\phi\}$$
is the additive group of $(-1)^k$-symmetric pairings on $K$, and
$$f^{\%}~=~S(j)~:~Q^{2k}(C) \to Q^{2k}(D)~;~\phi \mapsto f \phi f^*~=~j^*\phi j~=~
\phi\vert_L$$
sends such a pairing to its restriction to $L$.
A $2k$-dimensional symmetric (Poincar\'e) complex $(C,\phi\in Q^{2k}(C))$
is the same as a (nonsingular) $(-1)^k$-symmetric form $(K,\phi)$.
The relative symmetric $Q$-group of $f$
$$\begin{array}{ll}
Q^{2k+1}(f)&=~{\rm ker}(f^{\%}:Q^{2k}(C) \to Q^{2k}(D))\\[1ex]
&=~\{\phi \in S(K)\,\vert\, \phi^*=(-1)^k\phi \in S(K),
\phi\vert_L =0 \in S(L)\}~,
\end{array}$$
consists of the $(-1)^k$-symmetric pairings on $K$
which restrict to 0 on $L$. The submodule $L \subset K$ is a
lagrangian for $(K,\phi)$ if and only if $\phi$ restricts to 0 on $L$ and
$$L^{\perp}~=~\{x \in K \,\vert\, \phi(x)(L)=\{0\}\subset A\}~=~L~,$$
if and only if $(f:C \to D,(0,\phi)\in Q^{2k+1}(f))$
defines a $(2k+1)$-dimensional symmetric Poincar\'e pair,
with an exact sequence
$$\xymatrix{0 \ar[r] & D^k=L \ar[r]^{f^*=j} & C^k=K \ar[r]^{f\phi=j^*\phi} &
D_k=L^* \ar[r] &0~.}$$
Similarly for the quadratic case, with
$$\begin{array}{l}
Q_{2k}(C)~=~H_0(\Z_2;S(K),(-1)^kT)~,\\[1ex]
Q_{2k+1}(f)~=~\dfrac{\{(\psi,\chi) \in S(K) \oplus S(L)\,\vert\,
f^*\psi f=\chi+(-1)^{k+1}\chi^* \in S(L)\}}
{\{(\theta+(-1)^{k+1}\theta^*,f\theta f^*+\nu+(-1)^k\nu^*)\,\vert\,
\theta \in S(K),\nu \in S(L)\}}
\end{array}$$
A quadratic structure $\psi \in Q_{2k}(C)$ determines and is determined
by the pair $(\lambda,\mu)$ with $\lambda=\psi+(-1)^k\psi^* \in Q^{2k}(C)$
and
$$\mu~:~K \to H_0(\Z_2;A,(-1)^k)~;~x \mapsto \psi(x)(x)~.$$
A $(2k+1)$-dimensional (symmetric, quadratic) Poincar\'e pair
$(f:C \to D,(\delta\phi,\psi))$ is a nonsingular $(-1)^k$-quadratic
form $(K,\psi)$ together with a lagrangian $L\subset K$ for the nonsingular
$(-1)^k$-symmetric form $(K,\psi+(-1)^k\psi^*)$.

\begin{Lemma}
Let $(K,\psi)$ be a nonsingular $(-1)^k$-quadratic
form over $A$, and let $L \subset K$ be a lagrangian for
$(K,\psi+(-1)^k\psi^*)$.  There exists a direct complement
for $L \subset K$ which is also a lagrangian for $(K,\psi+(-1)^k\psi^*)$.
\end{Lemma}
\begin{proof} Choosing a direct complement $L' \subset K$ to $L \subset K$
write
$$\psi~=~ \begin{pmatrix} \mu & \lambda \\ 0 & \nu' \end{pmatrix}
~:~K~=~L \oplus L' \to K^*~=~L^* \oplus (L')^*$$
with $\lambda:L' \to L^*$ an isomorphism and
$$\mu+(-1)^k\mu^*~=~0~:~L \to L^*~.$$
In general $\nu'+(-1)^k(\nu')^* \neq 0:L^* \to L$, but if the
direct complement $L'$ is replaced by
$$L''~=~\{(-(\lambda^{-1})^*(\nu')^*(x),x) \in L \oplus L'
\,\vert\, x \in L'\} \subset K$$
and the isomorphism
$$\lambda''~:~L'' \to L^*~;~(-(\lambda^{-1})^*(\nu')^*(x),x) \mapsto
\lambda(x)$$
is used as an identification then
$$\psi~=~ \begin{pmatrix} \mu & 1 \\ 0 & \nu \end{pmatrix}
~:~K~=~L \oplus L^* \to K^*~=~L^* \oplus L$$
with $\nu=(\nu')^*\mu\nu':L^* \to L$ such that
$$\nu+(-1)^k\nu^*~=~0~:~L^* \to L~.$$
Thus $L''=L^*\subset K$ is a direct complement for $L$ which is
a lagrangian for $(K,\psi+(-1)^k\psi^*)$, with
$$\psi+(-1)^k\psi^*~=~\begin{pmatrix} 0 & 1 \\ (-1)^k & 0 \end{pmatrix}~:~
K~=~L \oplus L^* \to K^*~=~L^* \oplus L~.$$
\end{proof}

A lagrangian $L$ for the $(-1)^k$-symmetrization $(K,\psi+(-1)^k\psi^*)$
is a lagrangian for the $(-1)^k$-quadratic form
$(K,\psi)$ if and only if $\psi\vert_L=\mu$
is a $(-1)^{k+1}$-symmetrization, i.e.
$$\mu~=~\theta+(-)^{k+1}\theta^*~:~L \to L^*$$
for some $\theta \in S(L)$, in which case the inclusion $j:(L,0) \to (K,\psi)$
extends to an isomorphism of $(-1)^k$-quadratic forms
$$\begin{pmatrix} 1 & -\nu^* \\ 0 & 1 \end{pmatrix}~:~
H_{(-1)^k}(L)~=~(L\oplus L^*,\begin{pmatrix} 0 & 1 \\ 0 & 0 \end{pmatrix})
~\xymatrix{\ar[r]^-{\cong}&} ~(K,\psi)$$
with $\nu=\psi\vert_{L^*}$.
The $2k$-dimensional quadratic $L$-group $L_{2k}(A)$ is the
Witt group of stable isomorphism classes of nonsingular
$(-1)^k$-quadratic forms over $A$, such that
$$\begin{array}{ll}
(K,\psi)~=~(K',\psi') \in L_{2k}(A)&\hbox{\rm if and only
if there exists an isomorphism}\\[1ex]
&(K,\psi)\oplus H_{(-1)^k}(L)~\cong~(K',\psi')\oplus H_{(-1)^k}(L')~.
\end{array}$$

\begin{Proposition}\label{norminvform}
Given a $(-1)^k$-quadratic form $(L,\mu)$ over $A$ such that
$$\mu+(-1)^k\mu^*~=~0~:~L \to L^*$$
let $(B,\beta)$ be the chain bundle over $A$ given by
$$\begin{array}{l}
B~:~\dots \to 0 \to B_{k+1}~=~L \to 0 \to \dots~,\\[1ex]
\beta~=~\mu \in \widehat{Q}^0(B^{0-*})~=~{\rm Hom}_A(L,
\widehat{H}^{k+1}(\Z_2;A))~=~\widehat{H}^0(\Z_2;S(L),(-1)^{k+1}T)~.
\end{array}$$
{\rm (i)} The $(2k+1)$-dimensional twisted quadratic $Q$-group of $(B,\beta)$
$$\begin{array}{ll}
Q_{2k+1}(B,\beta)&=~\dfrac{
\{\nu \in S(L^*)\,\vert\, \nu+(-1)^k\nu^*=0\}}
{\{\phi - \phi \mu \phi^* -(\theta+(-1)^{k+1}\theta^*)\,\vert\,
\phi^*=(-1)^{k+1}\phi,\theta \in S(L^*)\}}\\[3ex]
&=~{\rm coker}(J_{\mu}:H^0(\Z_2;S(L^*),(-1)^{k+1}T) \to
\widehat{H}^0(\Z_2;S(L^*),(-1)^{k+1}T))
\end{array}$$
classifies nonsingular $(-1)^k$-quadratic forms $(K,\psi)$
over $A$ for which there exists a lagrangian $L$ for $(K,\psi+(-1)^k\psi^*)$ such that
$$\begin{array}{ll}
\psi\vert_L~=~\mu \in &{\rm im}(\widehat{H}^1(\Z_2;S(L),(-1)^kT)
\to H_0(\Z_2;S(L),(-1)^kT))\\[1ex]
&=~{\rm ker}(1+(-1)^kT:H_0(\Z_2;S(L),(-1)^kT) \to H^0(\Z_2;S(L),(-1)^kT))~.
\end{array}$$
Specifically, for any $(-1)^k$-quadratic form $(L^*,\nu)$ such that
$$\nu+(-1)^k\nu^*~=~0~:~L^* \to L$$
the nonsingular $(-1)^k$-quadratic form $(K,\psi)$ defined by
$$\psi~=~\begin{pmatrix} \mu & 1 \\ 0 & \nu \end{pmatrix}
~:~K~=~L \oplus L^* \to K^*~=~L^* \oplus L$$
is such that $L$ is a lagrangian of $(K,\psi+(-1)^k\psi^*)$, and
$$\partial~:~Q_{2k+1}(B,\beta) \to L_{2k}(A)~;~\nu \mapsto (K,\psi)~.$$
{\rm (ii)} The algebraic normal invariant of a $(2k+1)$-dimensional
{\rm (}symmetric, quadratic{\rm )}
Poincar\'e pair $(f:C \to D,(\delta\phi,\psi))$ concentrated in degree $k$ with
$$\begin{array}{l}
C_k~=~K^*~,~D_k~=~L^*~,\\[1ex]
f\psi_0f^* =~\mu \in {\rm ker}(1+(-1)^kT:H_0(\Z_2;S(L),(-1)^kT) \to
H^0(\Z_2;S(L),(-1)^kT))
\end{array}$$
is given by
$$(\phi,\theta)~=~\nu \in Q_{2k+1}(\Ca(f),\gamma)~=~Q_{2k+1}(B,\beta)$$
with
$$\widehat{v}_{k+1}(\gamma)~=~\widehat{v}_{k+1}(\beta)~:~
L~=~H_{k+1}(f)~=~H_{k+1}(B) \to \widehat{H}^{k+1}(\Z_2;A)~;~
x \mapsto \mu(x)(x)$$
and $\nu=\psi\vert_{L^*}$ the restriction of $\psi$ to any
lagrangian $L^* \subset K$ of $(K,\psi+(-1)^k\psi^*)$ complementary
to $L$.
\end{Proposition}
\begin{proof} (i) Given $(-1)^{k+1}$-symmetric forms $(L^*,\nu)$,
$(L^*,\phi)$ and $\theta \in S(L^*)$ replacing $\nu$ by
$$\nu'~=~\nu+\phi - \phi \mu \phi^* -(\theta+(-1)^{k+1}\theta^*)~:~
L^* \to L$$
results in a $(-1)^k$-quadratic form $(K,\psi')$ such that there is defined
an isomorphism
$$\begin{pmatrix} 1 & \phi^* \\ 0 & 1 \end{pmatrix}~:~(K,\psi') \to (K,\psi)$$
which is the identity on $L$. \\
(ii) This is the translation of Proposition \ref{norm} (iii)
into the language of forms and lagrangians.
\end{proof}

More generally~:

\begin{Proposition}
Given  $(-1)^k$-quadratic forms $(L,\mu)$, $(L^*,\nu)$ over $A$ such that
$$\mu+(-1)^k\mu^*~=~0~:~L \to L^*~,~\nu+(-1)^k\nu^*~=~0~:~L^* \to L$$
define a nonsingular $(-1)^k$-quadratic form
$$(K,\psi)~=~(L \oplus L^*,\begin{pmatrix} \mu & 1 \\ 0 & \nu \end{pmatrix})$$
such that $L$ and $L^*$ are complementary lagrangians of the
nonsingular $(-1)^k$-symmetric form
$$(K,\psi+(-1)^k\psi^*)~=~
(L \oplus L^*,\begin{pmatrix} 0 & 1 \\ (-1)^k & 0 \end{pmatrix})~,$$
and let $(f:C\to D,(\delta\phi,\psi))$ be the $(2k+1)$-dimensional
{\rm (}symmetric, quadratic{\rm )} Poincar\'e pair concentrated in
degree $k$ defined by
$$f~=~\begin{pmatrix} 1 & 0 \end{pmatrix}~:~C_k~=~K^*~=~L^*\oplus L
\to D_k~=~L^*~,~\delta \phi~=~0~,$$
with $\Ca(f) \simeq L_{*-k-1}$.\\
{\rm (i)} The Spivak normal bundle of $(f:C\to D,(\delta\phi,\psi))$
is given by
$$\gamma~=~\mu  \in \widehat{Q}^0(\Ca(f)^{0-*})~=~
\widehat{H}^0(\Z_2;S(L),(-1)^{k+1}T)~,$$
and
$$\begin{array}{ll}
Q_{2k+1}(\Ca(f),\gamma)&=~\dfrac{
\{\lambda \in S(L^*)\,\vert\, \lambda+(-1)^k\lambda^*=0\}}
{\{\phi - \phi \mu \phi^* -(\theta+(-1)^{k+1}\theta^*)\,\vert\,
\phi^*=(-1)^{k+1}\phi,\theta \in S(L^*)\}}\\[3ex]
&=~{\rm coker}(J_{\mu}:H^0(\Z_2;S(L^*),(-1)^{k+1}T) \to
\widehat{H}^0(\Z_2;S(L^*),(-1)^{k+1}T))~.
\end{array}$$
The algebraic normal invariant of $(f:C \to D,(\delta\phi,\psi))$ is
$$(\phi,\theta)~=~\nu \in Q_{2k+1}(\Ca(f),\gamma)~.$$
{\rm (ii)} Let $(B,\beta)$ be a chain bundle concentrated in degree $k+1$
$$\begin{array}{l}
B~:~\dots \to 0 \to B_{k+1} \to 0 \to \dots~,\\[1ex]
\beta \in \widehat{Q}^0(B^{0-*})~=~{\rm Hom}_A(B_{k+1},
\widehat{H}^{k+1}(\Z_2;A))~=~\widehat{H}^0(\Z_2;S(B_{k+1}),(-1)^{k+1}T)~,
\end{array}$$
so that
$$\begin{array}{ll}
Q_{2k+1}(B,\beta)&=~\dfrac{
\{\lambda \in S(B^{k+1})\,\vert\, \lambda+(-1)^k\lambda^*=0\}}
{\{\phi - \phi \beta \phi^* -(\theta+(-1)^{k+1}\theta^*)\,\vert\,
\phi^*=(-1)^{k+1}\phi,\theta \in S(B^{k+1})\}}\\[3ex]
&=~{\rm coker}(J_{\beta}:H^0(\Z_2;S(B^{k+1}),(-1)^{k+1}T) \to
\widehat{H}^0(\Z_2;S(B^{k+1}),(-1)^{k+1}T))~.
\end{array}$$
A $(B,\beta)$-structure on $(f:C\to D,(\delta\phi,\psi))$ is
given by a chain bundle map $(g,\chi):(\Ca(f),\gamma) \to (B,\beta)$,
corresponding to an $A$-module morphism $g:L \to B_{k+1}$ such that
$$g^*\beta g~=~\mu \in \widehat{H}^0(\Z_2;S(L),(-1)^{k+1}T)~,$$
with
$$(g,\chi)_{\%}~:~Q_{2k+1}(\Ca(f),\gamma) \to Q_{2k+1}(B,\beta)~;~
\lambda \mapsto g \lambda g^*~.$$
The 4-periodic $(B,\beta)$-structure cobordism class is thus given by
$$\begin{array}{ll}
(K,\psi;L)~=~(f:C\to D,(\delta\phi,\psi))&=~(g,\chi)_{\%}(\phi,\theta)~=~
g\nu g^*\\[1ex]
& \in \widehat{L}\langle B,\beta \rangle^{4*+2k+1}(A)~=~Q_{2k+1}(B,\beta)~,
\end{array}
$$
with
$$\begin{array}{ll}
(K,\psi)&=~(B_{k+1} \oplus B^{k+1}, \begin{pmatrix} \beta & 1 \\
0 & g\nu g^*\end{pmatrix})\\[2ex]
&\in {\rm im}(\partial:Q_{2k+1}(B,\beta) \to L_{2k}(A))~=~
{\rm ker}(L_{2k}(A) \to L\langle B,\beta \rangle^{4*+2k}(A))~.
\end{array} $$
\hfill$\qed$
\end{Proposition}

\subsection{The Generalized Arf Invariant}

\begin{Definition} \label{genArf}
{\rm The {\it generalized Arf invariant}
of a nonsingular $(-1)^k$-quadratic form $(K,\psi)$ over $A$
together with a lagrangian $L\subset K$ for the $(-1)^k$-symmetric form
$(K,\psi+(-1)^k\psi^*)$ is the image
$$(K,\psi;L)~=~(g,\chi)_{\%}(\phi,\theta) \in
\widehat{L}^{4*+2k+1}(A)~=~Q_{2k+1}(B^A,\beta^A)$$
of the algebraic normal invariant $(\phi,\theta) \in Q_{2k+1}(\Ca(f),\gamma)$
$($\ref{norminv}$)$ of the corresponding
$(2k+1)$-dimensional (symmetric, quadratic) Poincar\'e pair
$(f:C \to D,(\delta\phi,\psi)\in Q^{2k+1}_{2k+1}(f))$
$$\begin{array}{l}
(\phi,\theta)~=~\nu \in Q_{2k+1}(\Ca(f),\gamma)\\[1ex]
\hphantom{(\phi,\theta)~=~\nu \in}=~
{\rm coker}(J_{\mu}:H^0(\Z_2;S(L^*),(-1)^{k+1}T) \to
\widehat{H}^0(\Z_2;S(L^*),(-1)^{k+1}T))
\end{array}$$
under the morphism $(g,\chi)_{\%}$ induced by the classifying
chain bundle map $(g,\chi):(\Ca(f),\gamma) \to (B^A,\beta^A)$.
As in \ref{norminvform} $\nu=\psi\vert_{L^*}$ is the restriction of
$\psi$ to a lagrangian $L^* \subset K$ of $(K,\psi+(-1)^k\psi^*)$
complementary to $L$. \hfill$\qed$}
\end{Definition}

A nonsingular $(-1)^k$-symmetric formation $(K,\phi;L,L')$ is a
nonsingular $(-1)^k$-symmetric form $(K,\phi)$ together with two
lagrangians $L$, $L'$.  This type of formation is essentially the same
as a $(2k+1)$-dimensional symmetric Poincar\'e complex concentrated in
degrees $k,k+1$, and represents an element of $L^{4*+2k+1}(A)$.

\begin{Proposition}
{\rm (i)} The generalized Arf invariant is such that
$$(K,\psi;L)~=~0 \in Q_{2k+1}(B^A,\beta^A)~=~\widehat{L}^{4*+2k+1}(A)$$
if and only if there exists an isomorphism of $(-1)^k$-quadratic forms
$$(K,\psi) \oplus H_{(-1)^k}(L')~\cong~H_{(-1)^k}(L'')$$
such that
$$((K,\psi+(-1)^k\psi^*)\oplus (1+T)H_{(-1)^k}(L');L\oplus L',L'')~=~0
\in L^{4*+2k+1}(A)~.$$
{\rm (ii)} If $(K,\psi)$ is a nonsingular $(-1)^k$-quadratic form over $A$
and $L,L'\subset K$ are lagrangians for $(K,\psi+(-1)^k\psi^*)$ then
$$\begin{array}{l}
(K,\psi;L)-(K,\psi;L')~=~(K,\psi+(-1)^k\psi^*;L,L') \\[1ex]
\hskip25pt \in {\rm im}(L^{4*+2k+1}(A) \to \widehat{L}^{4*+2k+1}(A))
~=~{\rm ker}(\widehat{L}^{4*+2k+1}(A) \to L_{2k}(A))~.
\end{array}$$
\end{Proposition}
\begin{proof} This is the translation of the isomorphism
$Q_{2k+1}(B^A,\beta^A) \cong \widehat{L}^{4*+2k+1}(A)$ given
by \ref{cob} into the language of forms and formations.
\end{proof}

\begin{Example} {\rm Let $A$ be a field, so that
each $\widehat{H}^n(\Z_2;A)$ is a free $A$-module, and the universal
chain bundle over $A$ can be taken to be
$$\xymatrix@C-7pt{B^A~=~\widehat{H}^*(\Z_2;A)~:~ \dots \ar[r] &
B^A_n~=~\widehat{H}^n(\Z_2;A) \ar[r]^-{0} &
B^A_{n-1}~=~\widehat{H}^{n-1}(\Z_2;A) \ar[r]^-{0} & \dots~.}$$
If $A$ is a perfect field of characteristic 2 with the identity
involution squaring defines an $A$-module isomorphism
$$A \xymatrix{\ar[r]^-{\displaystyle{\cong}}&} \widehat{H}^n(\Z_2;A)~;~
a \mapsto a^2~.$$
Every nonsingular $(-1)^k$-quadratic
form over $A$ is isomorphic to one of the type
$$(K,\psi)~=~(L \oplus L^*,\begin{pmatrix} \mu & 1 \\ 0 & \nu \end{pmatrix})$$
with $L=A^{\ell}$ f.g. free and
$$\mu~=~(-1)^{k+1}\mu^*~:~L \to L^*~,~\nu~=~(-1)^{k+1}\nu^*~:~L^* \to L~.$$
For $j=1,2,\dots,\ell$ let
$$\begin{array}{l}
e_j~=~(0,\dots,0,1,0,\dots,0) \in L~,~g_j~=~\mu(e_j)(e_j) \in A~,\\[1ex]
e^*_j~=~(0,\dots,0,1,0,\dots,0) \in L^*~,~h_j~=~\nu(e^*_j)(e^*_j) \in A~.
\end{array}$$
The generalized Arf invariant in this case was identified in
Ranicki \cite[\S11]{ranicki4} with
the original invariant of Arf \cite{arf}
$$(K,\psi;L)~=~\sum\limits^{\ell}_{j=1}g_jh_j \in
Q_{2k+1}(B^A,\beta^A)~=~A/\{c+c^2\,\vert\, c \in A\}~.$$
\hfill$\qed$}
\end{Example}

For $k=0$ we have~:

\begin{Proposition} Suppose that the involution on $A$ is even.
If $(K,\psi)$ is a nonsingular quadratic form over
$A$ and $L$ is a lagrangian of $(K,\psi+\psi^*)$ then $L$ is a
lagrangian of $(K,\psi)$, the Witt class is
$$(K,\psi)~=~0 \in L_0(A)~,$$
the algebraic normal invariant is
$$(\phi,\theta)~=~0 \in Q_1(\Ca(f),\gamma)~=~0$$
and the generalized Arf invariant is
$$(K,\psi;L)~=~ (g,\chi)_{\%}(\phi,\theta)~=~0 \in \widehat{L}^{4*+1}(A)~=~Q_1(B^A,\beta^A)~.$$
\end{Proposition}
\begin{proof} By hypothesis $\widehat{H}^1(\Z_2;A)=0$, and $L=A^{\ell}$,
so that by Proposition \ref{norminvform} (i)
$$Q_1(\Ca(f),\gamma)~=~\widehat{H}^0(\Z_2;S(L^*),-T)~=~
\bigoplus\limits_{\ell}\widehat{H}^1(\Z_2;A)~=~0~.$$
\end{proof}

For $k=1$ we have~:

\begin{Theorem} \label{genArfthm}
Let $A$ be an $r$-even ring with $A_2$-module basis
$\{x_1=1,x_2,\dots,x_r\} \subset \widehat{H}^0(\Z_2;A)$,
and let
$$X~=~\begin{pmatrix} x_1 & 0 & 0 & \dots & 0 \\
0 & x_2 & 0 & \dots & 0 \\
0 & 0 & x_3 & \dots & 0 \\
\vdots & \vdots & \vdots & \ddots & \vdots \\
0 & 0 & 0 & \dots & x_r \end{pmatrix} \in {\rm Sym}_r(A)$$
so that by Theorem \ref{equ.ba-ba0}
$$Q_3(B^A,\beta^A)~=~\dfrac{{\rm Sym}_r(A)}
{{\rm Quad}_r(A) + \{L-LXL\,\vert\,L \in {\rm Sym}_r(A)\}}~.$$
{\rm (i)} Given $M \in {\rm Sym}_r(A)$ define the nonsingular $(-1)$-quadratic
form over $A$
$$(K_M,\psi_M)~=~(A^r \oplus (A^r)^*,\begin{pmatrix} X & 1 \\ 0 & M \end{pmatrix})$$
such that $L_M=A^r \subset K_M$ is a lagrangian of $(K_M,\psi_M-\psi^*_M)$.
The function
$$Q_3(B^A,\beta^A) \to \widehat{L}^{4*+3}(A)~;~M \mapsto (K_M,\psi_M;L_M)$$
is an isomorphism, with inverse given by the generalized Arf invariant.\\
{\rm (ii)}
Let $(K,\psi)$ be a nonsingular $(-1)$-quadratic form over $A$ of the type
$$(K,\psi)~=~(L\oplus L^*,\begin{pmatrix} \mu & 1 \\ 0 & \nu \end{pmatrix})$$
with
$$\mu-\mu^*~=~0~:~L \to L^*~,~\nu-\nu^*~=~0~:~L^* \to L$$
and let $g:L \to A^r$, $h:L^* \to A^r$ be $A$-module morphisms such that
$$\mu~=~g^* X g \in \widehat{H}^0(\Z_2;S(L),T)~,~
\nu~=~h^* X h \in \widehat{H}^0(\Z_2;S(L^*),T)~.$$
The generalized Arf invariant of $(K,\psi;L)$ is
$$(K,\psi;L)~=~g\nu g^*~=~gh^*Xhg^* \in  Q_3(B^A,\beta^A)~.$$
If $L=A^{\ell}$ then
$$g~=~(g_{ij})~:~L~=~A^{\ell} \to A^r~,~h~=~(h_{ij})~:~L^*~=~A^{\ell} \to A^r$$
with the coefficients $g_{ij},h_{ij} \in A$ such that
$$\begin{array}{c}
\mu(e_j)(e_j)~=~\sum\limits^r_{i=1}(g_{ij})^2x_i~,~
\nu(e^*_j)(e^*_j)~=~\sum\limits^r_{i=1}(h_{ij})^2x_i \in \widehat{H}^0(\Z_2;A)\\[2ex]
(e_j=(0,\dots,0,1,0,\dots,0) \in L=A^{\ell}~,~
e^*_j=(0,\dots,0,1,0,\dots,0) \in L^*=A^{\ell})
\end{array}$$
and
$$(K,\psi;L)~=~gh^*Xhg^* ~=~\begin{pmatrix} c_1 & 0 & 0 & \dots & 0 \\
0 & c_2 & 0 & \dots & 0 \\
0 & 0 & c_3 & \dots & 0 \\
\vdots & \vdots & \vdots & \ddots & \vdots \\
0 & 0 & 0 & \dots & c_r \end{pmatrix}
\in  Q_3(B^A,\beta^A)$$
with
$$c_i~=~\sum\limits^r_{k=1}(\sum\limits^{\ell}_{j=1}g_{ij}h_{kj})^2x_k
\in \widehat{H}^0(\Z_2;A)~.$$
{\rm (iii)} For any $M=(m_{ij}) \in {\rm Sym}_r(A)$ let $h=(h_{ij}) \in
M_r(A)$ be such that
$$m_{jj}~=~\sum\limits^r_{i=1}(h_{ij})^2x_i \in \widehat{H}^0(\Z_2;A)~
(1 \leqslant j \leqslant r)~,$$
so that
$$M~=~\begin{pmatrix} m_{11} & 0 & 0 & \dots & 0 \\
0 & m_{22} & 0 & \dots & 0 \\
0 & 0 & m_{33} & \dots & 0 \\
\vdots & \vdots & \vdots & \ddots & \vdots \\
0 & 0 & 0 & \dots & m_{rr}\end{pmatrix}~=~h^*Xh
\in \widehat{H}^0(\Z_2;M_r(A),T)~=~\dfrac{{\rm Sym}_r(A)}
{{\rm Quad}_r(A)}$$
and the generalized Arf invariant of
the triple $(K_M,\psi_M;L_M)$ in {\rm (i)} is
$$(K_M,\psi_M;L_M)~=~h^*Xh~=~M \in Q_3(B^A,\beta^A)$$
{\rm (}with $g=(\delta_{ij})$ here{\rm )}.
\end{Theorem}
\begin{proof} (i) The isomorphism
$Q_3(B^A,\beta^A) \to \widehat{L}^3(A);M \mapsto (K_M,\psi_M;L_M)$
is given by Proposition \ref{cob}.\\
(ii)  As in Definition \ref{genArf} let $(\phi,\theta) \in Q_3(\Ca(f),\gamma)$
be the algebraic normal invariant of the 3-dimensional
(symmetric, quadratic) Poincar\'e pair $(f:C \to D,(\delta\phi,\psi))$
concentrated in degree 1, with
$$f~=~\begin{pmatrix} 1 & 0 \end{pmatrix}~:~
C_1~=~K^*~=~L^* \oplus L \to D_1~=~L^*~,~\delta\phi~=~0~.$$
The $A$-module morphism
$$\widehat{v}_2(\gamma)~:~H_2(\Ca(f))~=~H^1(D)~=~
L \to \widehat{H}^0(\Z_2;A)~;~y \mapsto \mu(y)(y)$$
is induced by the $A$-module chain map
$$g~:~\Ca(f)~\simeq~L_{*-2} \to B^A(1)$$
and
$$(g,0)~:~(\Ca(f),\gamma) \to (B^A(1),\beta^A(1))
\to (B^A,\beta^A)$$
is a classifying chain bundle map. The induced morphism
$$\begin{array}{l}
(g,0)_{\%}~:~Q_3(\Ca(f),\gamma)~=~
{\rm coker}(J_{\mu}:H^0(\Z_2;S(L^*),T) \to \widehat{H}^0(\Z_2;S(L^*),T))\\[1ex]
\to Q_3(B^A,\beta^A)~=~{\rm coker}(J_X:H^0(\Z_2;M_r(A),T) \to
\widehat{H}^0(\Z_2;M_r(A),T))~;~\sigma \mapsto g\sigma g^*
\end{array}$$
sends the algebraic normal invariant
$$(\phi,\theta)~=~\nu~=~h^* X h \in Q_3(\Ca(f),\gamma)$$
to the generalized Arf invariant
$$(g,0)_{\%}(\phi,\theta)~=~gh^*Xhg^* \in Q_3(B^A,\beta^A)~.$$
(iii) By construction.
\end{proof}

In particular, the generalized Arf invariant for $A=\Z_2$ is just the
classical Arf invariant.

\section{The Generalized Arf Invariant for Linking Forms}

An $\epsilon$-quadratic formation $(Q,\psi;F,G)$ over $A$ corresponds
to a $1$-dimensional $\epsilon$-quadratic Poincar\'e complex.
The 1-dimensional $\epsilon$-quadratic $L$-group $L_1(A,\epsilon)$ is the Witt group
of $\epsilon$-quadratic formations, or equivalently the cobordism
group of 1-dimensional $\epsilon$-quadratic Poincar\'e complexes over $A$.
We could define a generalized Arf invariant $\alpha \in Q_2(B^A,\beta^A,\epsilon)$
for any formation with a null-cobordism of the 1-dimensional
$\epsilon$-symmetric Poincar\'e complex, so that
$$\begin{array}{l}
(Q,\psi;F,G)~=~\partial(\alpha)
\in {\rm ker}(1+T_{\epsilon}:L_1(A,\epsilon) \to L^{4*+1}(A,\epsilon))\\[1ex]
\hskip150pt =~{\rm im}(\partial:Q_2(B^{A,\epsilon},\beta^{A,\epsilon},\epsilon) \to L_1(A,\epsilon))~.
\end{array}$$
However, we do not need quite such a generalized Arf invariant here.
For our application to $\unil$, it suffices to work with a
localization $S^{-1}A$ of $A$ and to only consider a formation $(Q,\psi;F,G)$ such that
$$F \cap G~=~\{0\}~,~S^{-1}(Q/(F+G))~=~0$$
which corresponds to a $(-\epsilon)$-quadratic linking form $(T,\lambda,\mu)$
over $(A,S)$ with
$$T~=~Q/(F+G)~,~\lambda~:~T \times T \to S^{-1}A/A~.$$
Given a lagrangian $U \subset T$ for the $(-\epsilon)$-symmetric linking form
$(T,\lambda)$ we define in this section a `linking Arf invariant'
$$(T,\lambda,\mu;U) \in Q_2(B^{A,\epsilon},\beta^{A,\epsilon},\epsilon)~=~
\widehat{L}^{4*+2}(A,\epsilon)$$
such that
$$\begin{array}{l}
(Q,\psi;F,G)~=~\partial(T,\lambda,\mu;U) \in
{\rm ker}(1+T_{\epsilon}:L_1(A,\epsilon) \to L^{4*+1}(A,\epsilon))\\[1ex]
\hphantom{(Q,\psi;F,G)~=~\partial(T,\lambda,\mu;U) \in}~=~
{\rm im}(\partial:Q_2(B^{A,\epsilon},\beta^{A,\epsilon}) \to L_1(A,\epsilon))~.
\end{array}$$

\subsection{Linking Forms and Formations}

Given a ring with involution $A$ and a multiplicative subset $S \subset
A$ of central non-zero divisors such that $\overline{S}=S$ let
$S^{-1}A$ be the localized ring with involution obtained from $A$ by
inverting $S$. We refer to \cite{ranicki2} for the localization exact
sequences in $\epsilon$-symmetric and $\epsilon$-quadratic algebraic $L$-theory
$$\begin{array}{l}
\dots \to L^n(A,\epsilon) \to L_I^n(S^{-1}A,\epsilon) \to
L^n(A,S,\epsilon) \to L^{n-1}(A,\epsilon) \to \dots~,\\[1ex]
\dots \to L_n(A,\epsilon) \to L^I_n(S^{-1}A,\epsilon) \to
L_n(A,S,\epsilon) \to L_{n-1}(A,\epsilon) \to \dots
\end{array}$$
with $I={\rm im}(\widetilde{K}_0(A) \to \widetilde{K}_0(S^{-1}A))$,
$L^n(A,S,\epsilon)$ the cobordism group of $(n-1)$-dimensional
$\epsilon$-symmatric Poincar\'e complexes $(C,\phi)$ over $A$ such that
$H_*(S^{-1}C)=0$, and similarly for $L_n(A,S,\epsilon)$.

An {\it $(A,S)$-module} is an $A$-module $T$ with a
1-dimensional f.g. projective $A$-module resolution
$$\xymatrix@C-10pt{0 \ar[r] & P \ar[r]^{d} & Q \ar[r] & T \ar[r] & 0}$$
such that $S^{-1}d:S^{-1}P \to S^{-1}Q$ is an $S^{-1}A$-module
isomorphism. In particular,
$$S^{-1}T~=~0~.$$
The {\it dual $(A,S)$-module} is defined by
$$\begin{array}{ll}
T\widehat{~}&=~{\rm Ext}^1_A(T,A)~=~{\rm Hom}_A(T,S^{-1}A/A)\\[1ex]
&=~{\rm coker}(d^*:Q^* \to P^*)
\end{array}$$
with
$$A \times T\widehat{~} \to T\widehat{~}~;~
(a,f) \mapsto (x \mapsto f(x)\overline{a})~.$$
For any $(A,S)$-modules $T$, $U$ there is defined a duality isomorphism
$${\rm Hom}_A(T,U) \to {\rm Hom}_A(U\widehat{~},T\widehat{~})~;~ f
\mapsto f\widehat{~}$$
with
$$f\widehat{~}~:~U\widehat{~} \to T\widehat{~}~;~ g \mapsto (x \mapsto
g(f(x)))~.$$
An element $\lambda \in {\rm Hom}_A(T,T\widehat{~})$
can be regarded as a sesquilinear linking pairing
$$\lambda ~:~T \times T \to S^{-1}A/A~;~(x,y) \mapsto \lambda(x,y)~=~
\lambda(x)(y)$$
with
$$\begin{array}{l}
\lambda(x,ay+bz)~=~a\lambda(x,y)+b\lambda(x,z)~,\\[1ex]
\lambda(ay+bz,x)~=~\lambda(y,x)\overline{a}+\lambda(z,x)\overline{b}~,\\[1ex]
\lambda\widehat{~}(x,y)~=~\overline{\lambda(y,x)} \in
S^{-1}A/A~~(a,b \in A,\, x,y,z \in T)~.
\end{array}$$

\begin{Definition} {\rm Let $\epsilon = \pm 1$. \\
(i) An {\it $\epsilon$-symmetric linking form over $(A,S)$} $(T,\lambda)$
is an $(A,S)$-module $T$ together with $\lambda \in {\rm Hom}_A(T,T\widehat{~})$
such that $\lambda\widehat{~}=\epsilon\lambda$, so that
$$\overline{\lambda(x,y)}~=~\epsilon \lambda(y,x) \in S^{-1}A/A~~(x,y \in T)~.$$
The linking form is {\it nonsingular} if $\lambda:T \to T\widehat{~}$ is
an isomorphism. A {\it lagrangian} for $(T,\lambda)$ is an $(A,S)$-submodule
$U \subset T$ such that the sequence
$$\xymatrix@C-10pt{0 \ar[r] & U \ar[r]^-{j} &
T \ar[r]^-{j\widehat{~}\lambda} & U\widehat{~} \ar[r] & 0}$$
is exact with $j \in {\rm Hom}_A(U,T)$ the inclusion.
Thus $\lambda$ restricts to 0 on $U$ and
$$U^{\perp}~=~\{x \in T \,\vert\, \lambda(x)(U)=\{0\}\subset S^{-1}A/A\}~=~U~.$$
(ii) A {\it {\rm (}nonsingular{\rm )} $\epsilon$-quadratic linking form
over $(A,S)$} $(T,\lambda,\mu)$ is a (nonsingular) $\epsilon$-symmetric
linking form $(T,\lambda)$ together with a function
$$\mu~:~T \to Q_{\epsilon}(A,S)~=~\dfrac{
\{b \in S^{-1}A \,\vert\, \epsilon\overline{b}=b\}}
{\{a+\epsilon\overline{a}\,\vert\,a\in A\}}$$
such that
$$\begin{array}{l}
\mu(ax)~=~a\mu(x)\overline{a}~,\\[1ex]
\mu(x+y)~=~\mu(x)+\mu(y)+\lambda(x,y)+\lambda(y,x) \in Q_{\epsilon}(A,S)~,\\[1ex]
\mu(x)~=~\lambda(x,x) \in {\rm im}(Q_{\epsilon}(A,S)\to S^{-1}A/A)~~
(x,y \in T,a \in A)~.
\end{array}$$
A {\it lagrangian} $U\subset T$ for $(T,\lambda,\mu)$ is a lagrangian for $(T,\lambda)$
such that $\mu\vert_U=0$.\hfill$\qed$}
\end{Definition}

We refer to Ranicki \cite[3.5]{ranicki2} for the development of the
theory of $\epsilon$-symmetric and $\epsilon$-quadratic linking formations
over $(A,S)$.

From now on, we shall only be concerned with $A,S$ which satisfy~:

\begin{Hypothesis} \label{hyp}
{\rm $A,S$ are such that
$$\widehat{H}^*(\Z_2;S^{-1}A)~=~0~.$$
\hfill$\qed$}
\end{Hypothesis}

\begin{Example} {\rm Hypothesis \ref{hyp} is satisfied if $1/2 \in S^{-1}A$,
e.g. if $A$ is even and
$$S~=~(2)^{\infty}~=~\{2^i\,\vert\,i \geqslant 0\} \subset A~,~
S^{-1}A~=~A[1/2]~.$$
\hfill$\qed$}
\end{Example}

\begin{Proposition} \label{Sbraid}
{\rm (i)} For $n=2$ $($resp. $1)$ the relative group $L^n(A,S,\epsilon)$ in the $\epsilon$-symmetric $L$-theory
localization exact sequence
$$\dots \to L^n(A,\epsilon) \to L_I^n(S^{-1}A,\epsilon) \to L^n(A,S,\epsilon) \to
L^{n-1}(A,\epsilon) \to \dots$$
is the Witt group of nonsingular $(-\epsilon)$-symmetric linking forms
$($resp. $\epsilon$-symmetric linking formations$)$ over $(A,S)$,
with $I={\rm im}(\widetilde{K}_0(A) \to \widetilde{K}_0(S^{-1}A))$.
The skew-suspension maps
$$\overline{S}~:~L^n(A,S,\epsilon) \to L^{n+2}(A,S,-\epsilon)~(n \geqslant 1)$$
are isomorphisms if and only if the skew-suspension maps
$$\overline{S}~:~L^n(A,\epsilon) \to L^{n+2}(A,-\epsilon)~(n \geqslant 0)$$
are isomorphisms.\\
{\rm (ii)} The relative group $L_n(A,S,\epsilon)$ for $n=2k$ $($resp. $2k+1)$
in the $\epsilon$-quadratic $L$-theory localization exact sequence
$$\dots \to L_n(A,\epsilon) \to L^I_n(S^{-1}A,\epsilon) \to
L_n(A,S,\epsilon) \to L_{n-1}(A,\epsilon) \to \dots$$
is the Witt group of nonsingular $(-1)^k\epsilon$-quadratic linking forms
$($resp. formations$)$ over $(A,S)$.\\
{\rm (iii)} The 4-periodic $\epsilon$-symmetric and
$\epsilon$-quadratic localization exact sequences interleave in a
commutative braid of exact sequences

$$\xymatrix@C-20pt{
Q_{n+1}(B^A,\beta^A,\epsilon)\ar[dr]^-{\partial^S}\ar@/^2pc/[rr]^-{\partial}&&
L_n(A,\epsilon)  \ar[dr]\ar@/^2pc/[rr]  &&L^I_n(S^{-1}A,\epsilon)  \\
& L_{n+1}(A,S,\epsilon)\ar[ur]\ar[dr] && \hbox{~~$L^{n+4*}(A,\epsilon)$~~} \ar[ur] \ar[dr]&&\\
L^I_{n+1}(S^{-1}A,\epsilon) \ar[ur]\ar@/_2pc/[rr]_-{}&&L^{n+4*+1}(A,S,\epsilon)
\ar[ur]\ar@/_2pc/[rr]_{}&&Q_n(B^A,\beta^A,\epsilon)}$$

\vskip3mm

\end{Proposition}
\begin{proof} (i)+(ii) See \cite[\S3]{ranicki2}.\\
(iii) For $A,S$ satisfying Hypothesis \ref{hyp}
the $\epsilon$-symmetrization maps for the $L$-groups of $S^{-1}A$ are
isomorphisms
$$1+T_{\epsilon}~:~L^I_n(S^{-1}A,\epsilon) \xymatrix{\ar[r]^-{\cong}&}
L_I^n(S^{-1}A,\epsilon)~.$$
\end{proof}

\begin{Definition} {\rm
{\rm (i)}  An {\it $\epsilon$-quadratic $S$-formation}
$(Q,\psi;F,G)$ over $A$ is an $\epsilon$-quadratic formation such that
$$S^{-1}F \oplus S^{-1}G~=~S^{-1}Q~,$$
or equivalently such that $Q/(F +G)$ is an $(A,S)$-module.\\
{\rm (iii)} A {\it stable isomorphism} of $\epsilon$-quadratic $S$-formations
over $A$
$$[f]~:~(Q_1,\psi_1;F_1,G_1) \to (Q_2,\psi_2;F_2,G_2)$$
is an isomorphism of the type
$$f~:~(Q_1,\psi_1;F_1,G_1)\oplus (N_1,\nu_1;H_1,K_1)
 \to (Q_2,\psi_2;F_2,G_2)\oplus  (N_2,\nu_2;H_2,K_2) $$
with $N_1=H_1 \oplus K_1$, $N_2=H_2 \oplus K_2$. \hfill$\qed$}
\end{Definition}

\begin{Proposition} \label{lform}
{\rm (i)} A $(-\epsilon)$-quadratic $S$-formation $(Q,\psi;F,G)$ over $A$
determines a nonsingular $\epsilon$-quadratic linking form $(T,\lambda,\mu)$
over $(A,S)$, with
$$\begin{array}{l}
T~=~Q/(F+G)~,\\[1ex]
\lambda~:~T \times T \to S^{-1}A/A~;~(x,y) \mapsto (\psi-\epsilon \psi^*)(x)(z)/s\\[1ex]
\mu~:~T  \to Q_{\epsilon}(A,S)~;~y \mapsto (\psi-\epsilon \psi^*)(x)(z)/s - \psi(y)(y)\\[1ex]
(x,y \in Q~,~z \in G~,~s \in S~,~sy-z \in F)~.
\end{array}$$
{\rm (ii)} The isomorphism classes of nonsingular $\epsilon$-quadratic
linking forms over $A$ are in one-one correspondence with the stable
isomorphism classes of $(-\epsilon)$-quadratic $S$-formations over $A$.
\end{Proposition}
\begin{proof} See Proposition 3.4.3 of \cite{ranicki2}.
\end{proof}

For any  $S^{-1}A$-contractible f.g. projective
$A$-module chain complexes concentrated in degrees $k,k+1$
$$\begin{array}{l}
C~:~\dots \to 0 \to C_{k+1} \to C_k \to 0 \to \dots~,\\[1ex]
D~:~\dots \to 0 \to D_{k+1} \to D_k \to 0 \to \dots
\end{array}$$
there are natural identifications
$$\begin{array}{l}
H^{k+1}(C)~=~H_k(C)\widehat{~}~,~H_k(C)~=~H^{k+1}(C)\widehat{~}~,\\[2ex]
H^{k+1}(D)~=~H_k(D)\widehat{~}~,~H_k(D)~=~H^{k+1}(D)\widehat{~}~,\\[2ex]
H_0({\rm Hom}_A(C,D))~=~{\rm Hom}_A(H_k(C),H_k(D))~=~
{\rm Tor}_1^A(H^{k+1}(C),H_k(D))~,\\[2ex]
H_1({\rm Hom}_A(C,D))~=~H^{k+1}(C)\otimes_A H_k(D)~=~
{\rm Ext}^1_A(H_k(C),H_k(D))~,\\[2ex]
H_{2k}(C\otimes_AD)~=~H_k(C)\otimes_AH_k(D)~=~
{\rm Ext}_A^1(H^{k+1}(C),H_k(D))~,\\[2ex]
H_{2k+1}(C\otimes_AD)~=~{\rm Hom}_A(H^{k+1}(C),H_k(D))~=~
{\rm Tor}_1^A(H_k(C),H_k(D))~.
\end{array}$$
In particular, an element $\lambda \in H_{2k+1}(C\otimes_AD)$
is a sesquilinear linking pairing
$$\lambda~:~H^{k+1}(C) \times H^{k+1}(D) \to S^{-1}A/A~.$$
An element $\phi \in H_{2k}(C\otimes_AD)$ is a chain
homotopy class of chain maps $\phi:C^{2k-*} \to D$,
classifying the extension
$$0 \to H_k(D) \to H_k(\phi) \to H^{k+1}(C) \to 0~.$$

\begin{Proposition} \label{link}
Given an $(A,S)$-module $T$ let
$$B~:~\dots \to 0 \to B_{k+1} \xymatrix{\ar[r]^-{d}&} B_k \to 0 \to \dots$$
be a f.g. projective $A$-module chain complex
concentrated in degrees $k,k+1$ such that $H^{k+1}(B)=T$, $H^k(B)=0$,
so that $H_k(B)=T\widehat{~}$, $H_{k+1}(B)=0$.
The $Q$-groups in the exact sequence
$$\xymatrix{Q^{2k+2}(B)=0 \ar[r]&
\widehat{Q}^{2k+2}(B) \ar[r]^-{H}&} Q_{2k+1}(B)
\xymatrix{\ar[r]^-{1+T}&} Q^{2k+1}(B)
\xymatrix{\ar[r]^-{J}&} \widehat{Q}^{2k+1}(B)$$
have the following interpretation in terms of $T$.\\
{\rm (i)} The symmetric $Q$-group
$$Q^{2k+1}(B)~=~H^0(\Z_2;{\rm Hom}_A(T,T\widehat{~}),(-1)^{k+1})$$
is the additive group of $(-1)^{k+1}$-symmetric linking pairings $\lambda$ on $T$,
with $\phi \in Q^{2k+1}(B)$ corresponding to
$$\lambda ~:~T \times T \to
S^{-1}A/A~;~(x,y) \mapsto \phi_0(d^*)^{-1}(x)(y) ~~(x,y \in B^{k+1})~.$$
{\rm (ii)} The quadratic $Q$-group
$$\begin{array}{l}
Q_{2k+1}(B)~=\\[1ex]
~\displaystyle{\frac{\{
(\psi_0,\psi_1) \in {\rm Hom}_A(B^k,B_{k+1}) \oplus S(B^k)
\,\vert\, d\psi_0=\psi_1+(-1)^{k+1}\psi_1^* \in S(B^k)\}}
{\{((\chi_0+(-1)^{k+1}\chi_0^*)d^*,d\chi_0d^*+\chi_1+(-1)^k\chi_1^*)
\,\vert\,(\chi_0,\chi_1) \in S(B^{k+1}) \oplus S(B^k)\}}}
\end{array}$$
is the additive group of $(-1)^{k+1}$-quadratic linking structures
$(\lambda,\mu)$ on $T$. The element
$\psi=(\psi_0,\psi_1)\in Q_{2k+1}(B)$ corresponds to
$$\begin{array}{l}
\lambda ~:~T \times T \to S^{-1}A/A~;~
(x,y) \mapsto \psi_0(d^*)^{-1}(x)(y) ~~(x,y \in B^{k+1})~,\\[1ex]
\mu~:~T \to Q_{(-1)^{k+1}}(A,S)~;~ x \mapsto \psi_0(d^*)^{-1}(x)(x)~.
\end{array}$$
{\rm (iii)} The hyperquadratic $Q$-groups of $B$
$$\widehat{Q}^n(B)~=~H_n(\widehat{d}^{\%}:\widehat{W}^{\%}B_{k+1} \to
\widehat{W}^{\%}B_k)$$
are such that
$$\begin{array}{l}
\widehat{Q}^{2k}(B)~=~
\dfrac{\{(\delta,\chi) \in S(B^{k+1}) \oplus S(B^k) \,\vert\,
\delta^*=(-1)^{k+1}\delta,d \delta d^*=\chi+(-1)^{k+1}\chi^*\}}
{\{(\mu+(-1)^{k+1}\mu^*,d\mu d^*+\nu+(-1)^k\nu^*)\,\vert\,
(\mu,\nu) \in S(B^{k+1}) \oplus S(B^k) \}}~,\\[3ex]
\widehat{Q}^{2k+1}(B)~=~
\dfrac{\{(\delta,\chi) \in S(B^{k+1}) \oplus S(B^k) \,\vert\,
\delta^*=(-1)^k\delta,d \delta d^*=\chi+(-1)^k\chi^*\}}
{\{(\mu+(-1)^k\mu^*,d\mu d^*+\nu+(-1)^{k+1}\nu^*)\,\vert\,
(\mu,\nu) \in S(B^{k+1}) \oplus S(B^k) \}}~,
\end{array}$$
with universal coefficient exact sequences
$$\begin{array}{l}
0 \to T\widehat{~}\otimes_A \widehat{H}^k(\Z_2;A)
\to \widehat{Q}^{2k}(B) \xymatrix{\ar[r]^{\widehat{v}_{k+1}}&}
{\rm Hom}_A(T,\widehat{H}^{k+1}(\Z_2;A)) \to 0~,\\[2ex]
0 \to T\widehat{~}\otimes_A \widehat{H}^{k+1}(\Z_2;A)
\to \widehat{Q}^{2k+1}(B)
\xymatrix{\ar[r]^{\widehat{v}_k}&}{\rm Hom}_A(T,\widehat{H}^k(\Z_2;A)) \to 0~.
\end{array}$$
\hfill$\qed$
\end{Proposition}

Let $f:C \to D$ be a chain map of $S^{-1}A$-contractible $A$-module
chain complexes concentrated in degrees $k,k+1$, inducing the $A$-module
morphism
$$f^*~=~j~:~U~=~H^{k+1}(D) \to T~=~H^{k+1}(C)~.$$
By Proposition \ref{link} (i) a $(2k+1)$-dimensional symmetric
Poincar\'e complex $(C,\phi)$ is essentially the same as a nonsingular
$(-1)^{k+1}$-symmetric linking form $(T,\lambda)$, and a
$(2k+2)$-dimensional symmetric Poincar\'e pair $(f:C \to D,(\delta
\phi,\phi))$ is essentially the same as a lagrangian $U$ for
$(T,\lambda)$, with $j=f^*:U \to T$ the inclusion.  Similarly, a
$(2k+1)$-dimensional quadratic Poincar\'e complex $(C,\psi)$ is
essentially the same as a nonsingular $(-1)^{k+1}$-quadratic linking
form $(T,\lambda,\mu)$, and a $(2k+2)$-dimensional quadratic Poincar\'e
pair $(f:C \to D,(\delta \psi,\psi))$ is essentially the same as a
lagrangian $U\subset T$ for $(T,\lambda,\mu)$.  A $(2k+2)$-dimensional
(symmetric, quadratic) Poincar\'e pair $(f:C \to D,(\delta\phi,\psi))$
is a nonsingular $(-1)^{k+1}$-quadratic form $(T,\lambda,\mu)$ together
with a lagrangian $U\subset T$ for the nonsingular
$(-1)^{k+1}$-symmetric linking form $(T,\lambda)$.

\begin{Proposition} \label{linktw}
Let $U$ be an $(A,S)$-module together with an $A$-module morphism
$\mu_1:U \to \widehat{H}^{k+1}(\Z_2;A)$, defining a
$(-1)^{k+1}$-quadratic linking form $(U,\lambda_1,\mu_1)$ over $(A,S)$
with $\lambda_1=0$.\\
{\rm (i)} There exists a map of chain bundles
$(d,\chi):(B_{k+2},0) \to (B_{k+1},\delta)$
concentrated in degree $k+1$ such that the cone chain bundle
$(B,\beta)=\Ca(d,\chi)$ has
$$\begin{array}{l}
H_{k+1}(B)~=~U~,~H^{k+2}(B)~=~U\widehat{~}~,~
H_{k+2}(B)~=~H^{k+1}(B)~=~0~,\\[1ex]
\beta~=~[\delta]~=~\mu_1 \in
\widehat{Q}^0(B^{0-*})~=~{\rm Hom}_A(U,\widehat{H}^{k+1}(\Z_2;A))~.
\end{array}$$
{\rm (ii)} The $(2k+2)$-dimensional twisted quadratic $Q$-group of $(B,\beta)$
as in {\rm (i)}
$$\begin{array}{l}
Q_{2k+2}(B,\beta)~=\\[2ex]
\hskip20pt
\dfrac{\{(\phi,\theta) \in S(B^{k+1}) \oplus S(B^{k+1}) \,\vert\,
\phi^*=(-1)^{k+1}\phi,\phi-\phi \delta\phi^*=\theta+(-1)^{k+1}\theta^*\}}
{\{(d,\chi)_{\%}(\nu)+(0,\eta+(-1)^k\eta^*)\,\vert\,
\nu \in S(B^{k+2}), \eta \in S(B^{k+1})\}}\\[3ex]
((d,\chi)_{\%}(\nu)=(d(\nu+(-1)^{k+1}\nu^*)d^*,d\nu d^*-d(\nu+(-1)^{k+1}\nu^*)\chi
(\nu^*+(-1)^{k+1}\nu)d^*))
\end{array}$$
is the additive group of isomorphism classes of
extensions of $U$ to a nonsingular $(-1)^{k+1}$-quadratic linking form
$(T,\lambda,\mu)$ over $(A,S)$ such that $U\subset T$ is a lagrangian of the
$(-1)^{k+1}$-symmetric linking form $(T,\lambda)$ and
$$\beta~ =~ \mu \vert_U~:~H_{k+1}(B)~=~U \to \widehat{H}^{k+1}(\Z_2;A)~=~
{\rm ker}(Q_{(-1)^{k+1}}(A,S) \to S^{-1}A/A)~.$$
{\rm (iii)} An element $(\phi,\theta) \in Q_{2k+2}(B,\beta)$ is the algebraic
normal invariant $(\ref{norminv})$ of the $(2k+2)$-dimensional
$($symmetric, quadratic$)$ Poincar\'e pair $(f:C \to D,(\delta\phi,\psi) \in
Q^{2k+2}_{2k+2}(f))$ with
$$\begin{array}{l}
d_C~=~\begin{pmatrix} d & \phi \\ 0 & d^* \end{pmatrix}~:~
C_{k+1}~=~B_{k+2} \oplus B^{k+1}  \to C_k~=~B_{k+1} \oplus B^{k+2}~,\\[2ex]
f~=~{\rm projection}~:~C \to D~=~B^{2k+2-*}
\end{array}$$
constructed as in Proposition  {\rm \ref{norm} (ii)},
corresponding to the quadruple $(T,\lambda,\mu;U)$ given by
$$j~=~f^*~:~U~=~H^{k+1}(D)~=~H_{k+1}(B) \to T~=~H^{k+1}(C)~.$$
The $A$-module extension
$$0 \to U \to T \to U\widehat{~} \to 0$$
is classified by
$$[\phi] \in H_{2k+2}(B\otimes_AB)~=~U\otimes_AU~=~{\rm Ext}^1_A(U\widehat{~},U)~.$$
{\rm (iv)} The $(-1)^{k+1}$-quadratic linking form $(T,\lambda,\mu)$
in {\rm (iii)} corresponds to
the $(-1)^k$-quadratic $S$-formation $(Q,\psi;F,G)$ with
$$\begin{array}{l}
(Q,\psi)~=~H_{(-1)^k}(F)~,~F~=~B_{k+2} \oplus B^{k+1}~,\\[1ex]
G~=~{\rm im}\bigg( \begin{pmatrix} 1 & 0 \\ -\delta d & 1 -\delta \phi \\
0 & (-1)^{k+1}d ^* \\ d & \phi \end{pmatrix}:B_{k+2} \oplus B^{k+1}
\to B_{k+2} \oplus B^{k+1} \oplus B^{k+2} \oplus B_{k+1}\bigg)
\subset F \oplus F^*
\end{array}$$
such that
$$F \cap G~=~\{0\}~,~Q/(F+G)~=~H^{k+1}(C)~=~T~.$$
The inclusion $U \to T$ is resolved by
$$\xymatrix@C+26pt@R+10pt{
0 \ar[r] & B_{k+2} \ar[r]_-{\displaystyle{d}} \ar[d]^-{\begin{pmatrix} 1 \\ 0 \end{pmatrix}}
& B_{k+1} \ar[r] \ar[d]^-{\begin{pmatrix} 0 \\ 1 \end{pmatrix}} &
U \ar[r] \ar[d] & 0 \\
0 \ar[r] & B_{k+2}\oplus B^{k+1} \ar[r]^-{
\begin{pmatrix} 0 & (-)^{k+1}d^* \\
d & \phi \end{pmatrix} } & B^{k+2} \oplus B_{k+1} \ar[r] & T
\ar[r] & 0}$$
{\rm (v)} If the involution on $A$ is even and $k=-1$ then
$$Q_0(B,\beta)~=~\dfrac{\{\phi \in {\rm Sym}(B^0) \,\vert\,
\phi-\phi \delta\phi\in {\rm Quad}(B^0)\}}
{\{d \sigma  d^*\,\vert\,\sigma \in {\rm Quad}(B^1)\}}~.$$
An extension of $U={\rm coker}(d:B_1 \to B_0)$ to a nonsingular
quadratic linking form $(T,\lambda,\mu)$ over $(A,S)$ with
$\mu\vert_U=\mu_1$  and $U\subset T$ a lagrangian of $(T,\lambda)$
is classified by $\phi \in Q_0(B,\beta)$ such that
$\lambda:T \to T\widehat{~}$ is resolved by
$$\xymatrix@C+26pt@R+10pt{
0 \ar[r] & B_1 \oplus B^0
\ar[d]_-{\begin{pmatrix} 1 & 0 \\ -\delta d & 1 - \delta \phi \end{pmatrix} }
\ar[r]^-{\begin{pmatrix} 0 & d^* \\ d & \phi \end{pmatrix} }
& B^1 \oplus B_0 \ar[r]
\ar[d]^-{\begin{pmatrix} 1 & -d^*\delta \\ 0 & 1 -  \phi\delta \end{pmatrix} } &
T \ar[r] \ar[d]^-{\displaystyle{\lambda}} & 0 \\
0 \ar[r] & B_1\oplus B^0 \ar[r]^-{\begin{pmatrix} 0 & d^* \\
d & \phi \end{pmatrix} } & B^1 \oplus B_0 \ar[r] & T\widehat{~}\ar[r] & 0}$$
and
$$\begin{array}{l}
T~=~{\rm coker}(\begin{pmatrix} 0 & d^* \\ d & \phi \end{pmatrix}:
B_1 \oplus B^0 \to B^1 \oplus B_0)~,\\[1ex]
\lambda~:~T \times T \to S^{-1}A/A~;\\[1ex]
((x_1,x_0),(y_1,y_0)) \mapsto -d^{-1}\phi (d^*)^{-1}(x_1)(y_1)
+ d^{-1}(x_1)(y_0) + (d^*)^{-1}(x_0)(y_1)~,\\[1ex]
\mu~:~T \to Q_{+1}(A,S)~;\\[1ex]
(x_1,x_0) \mapsto -d^{-1}\phi (d^*)^{-1}(x_1)(x_1)
+ d^{-1}(x_1)(x_0) + (d^*)^{-1}(x_0)(x_1) - \delta(x_0)(x_0)~,\\[1ex]
(x_0,y_0 \in B_0~,~x_1,y_1 \in B^1)~.
\end{array}$$
\hfill$\qed$
\end{Proposition}

\subsection{The Linking Arf Invariant}

\begin{Definition} \label{linkArf} {\rm The {\it linking Arf invariant}
of a nonsingular $(-1)^{k+1}$-quadratic linking form $(T,\lambda,\mu)$
over $(A,S)$ together with a lagrangian $U\subset T$ for $(T,\lambda)$
is the image
$$(T,\lambda,\mu;U)~=~(g,\chi)_{\%}(\phi,\theta)\in
\widehat{L}^{4*+2k+2}(A)~=~Q_{2k+2}(B^A,\beta^A)$$
of the algebraic normal invariant $(\phi,\theta) \in Q_{2k+2}(\Ca(f),\gamma)$
$($\ref{norminv}$)$ of the corresponding
$(2k+2)$-dimensional (symmetric, quadratic) Poincar\'e pair
$(f:C \to D,(\delta\phi,\psi)\in Q^{2k+2}_{2k+2}(f))$ concentrated
in degrees $k,k+1$ with
$$f^*~=~j~:~H^{k+1}(D)~=~U \to H^{k+1}(C)~=~T~,$$
and $(g,\chi)_{\%}$ induced by the classifying
chain bundle map $(g,\chi):(\Ca(f),\gamma) \to (B^A,\beta^A)$.\\
\hbox{~~\hskip100pt} \hfill$\qed$}
\end{Definition}

The chain bundle $(\Ca(f),\gamma)$ in \ref{linkArf} is (up to
equivalence) of the type $(B,\beta)$ considered in Proposition
\ref{linktw} (i)\,: the algebraic normal invariant $(\phi,\theta)\in
Q_{2k+2}(B,\beta)$ classifies the extension of $(U,\beta)$ to a
lagrangian of a $(-1)^{k+1}$-symmetric linking form $(T,\lambda)$ with
a $(-1)^{k+1}$-quadratic function $\mu$ on $T$ such that
$\mu\vert_U=\beta$.  The linking Arf invariant $(T,\lambda,\mu;U) \in
Q_{2k+2}(B^A,\beta^A)$ gives the Witt class of $(T,\lambda,\mu;U)$.
The boundary map
$$\partial~:~Q_{2k+2}(B^A,\beta^A) \to L_{2k+1}(A)~;~
(T,\lambda,\mu;U) \mapsto (Q,\psi;F,G)$$
sends the linking Arf invariant to
the Witt class of the $(-1)^k$-quadratic formation  $(Q,\psi;F,G)$
constructed in \ref{linktw} (iv).

\begin{Theorem} \label{linkArfthm}
Let $A$ be an $r$-even ring with  $A_2$-module basis
$\{x_1=1,x_2,\dots,x_r\} \subset \widehat{H}^0(\Z_2;A)$, and let
$$X~=~\begin{pmatrix} x_1 & 0 & 0 & \dots & 0 \\
0 & x_2 & 0 & \dots & 0 \\
0 & 0 & x_3 & \dots & 0 \\
\vdots & \vdots & \vdots & \ddots & \vdots \\
0 & 0 & 0 & \dots & x_r \end{pmatrix} \in {\rm Sym}_r(A)~,$$
so that by Theorem \ref{equ.ba-ba0}
$$Q_{2k}(B^A,\beta^A)~=~
\begin{cases}
\dfrac{\{M \in {\rm Sym}_r(A)\,\vert\, M-MXM \in {\rm Quad}_r(A)\}}
{4{\rm Quad}_r(A) + \{2(N+N^t)-N^tXN\,\vert\,N \in M_r(A)\}}&\hbox{if $k=0$}\\[1ex]
0&\hbox{if $k=1$}~.
\end{cases}$$
{\rm (i)} Let
$$S~=~(2)^{\infty} \subset A~,$$
so that
$$S^{-1}A~=~A[1/2]$$
and $\widehat{H}^0(\Z_2;A)$ is an $(A,S)$-module. The
hyperquadratic $L$-group
$\widehat{L}^0(A)$ fits into the exact sequence
$$\dots \to L^1(A,S) \to \widehat{L}^0(A) \to L_0(A,S) \to
L^0(A,S) \to \dots~.$$
The linking Arf invariant of a nonsingular quadratic
linking form $(T,\lambda,\mu)$ over $(A,S)$ with a lagrangian
$U \subset T$ for $(T,\lambda)$ is the Witt class
$$(T,\lambda,\mu;U) \in Q_0(B^A,\beta^A)~=~\widehat{L}^{4*}(A)~.$$
{\rm (ii)} Given $M \in {\rm Sym}_r(A)$ such that $M-MXM \in {\rm Quad}_r(A)$
let $(T_M,\lambda_M,\mu_M)$ be the nonsingular quadratic linking form
over $(A,S)$ corresponding to the $(-1)$-quadratic $S$-formation over $A$
(\ref{lform})
$$(Q_M,\psi_M;F_M,G_M)~=~ (H_-(A^{2r});A^{2r},
{\rm im}\bigg(\begin{pmatrix}
\begin{pmatrix} I & 0 \\ -2X & I-XM \end{pmatrix} \\[3ex]
\begin{pmatrix} 0 & 2I \\ 2I & M \end{pmatrix}
\end{pmatrix}:A^{2r} \to A^{2r} \oplus A^{2r}\bigg))$$
and let
$$U_M~=~(A_2)^r \subset T_M~=~Q_M/(F_M+G_M)~=~{\rm coker}(G_M \to F^*_M)$$
be the lagrangian for the nonsingular symmetric linking form
$(T_M,\lambda_M)$ over $(A,S)$ with the inclusion $U_M \to T_M$
resolved by
$$\xymatrix@C+23pt@R+10pt{
0 \ar[r] & A^r \ar[r]^-{\displaystyle{2I}}
\ar[d]_-{\begin{pmatrix} 1 \\ 0 \end{pmatrix}}
& A^r\ar[r] \ar[d]^-{\begin{pmatrix} 0 \\ 1 \end{pmatrix}} &
U_M \ar[r] \ar[d] & 0 \\
0 \ar[r] & A^r \oplus A^r \ar[r]^-{
\begin{pmatrix} 0 & 2I \\
2I & M \end{pmatrix} } & A^r  \oplus A^r \ar[r] & T_M
\ar[r] & 0}$$
The function
$$Q_0(B^A,\beta^A) \to \widehat{L}^{4*}(A)~;~M \mapsto (T_M,\lambda_M,\mu_M;U_M)$$
is an isomorphism, with inverse given by the linking Arf invariant.\\
{\rm (iii)} Let $(T,\lambda,\mu)$ be a nonsingular quadratic linking
form over $(A,S)$ together with a lagrangian $U \subset T$ for $(T,\lambda)$.
For any f.g. projective $A$-module resolution of $U$
$$0 \to B_1 \xymatrix{\ar[r]^-{d}&} B_0 \to U \to 0$$
let
$$\delta \in {\rm Sym}(B_0)~,~\phi \in {\rm Sym}(B^0)~,~
\beta ~=~[\delta]~=~\mu\vert_U\in \widehat{Q}^0(B^{0-*})~=~
{\rm Hom}_A(U,\widehat{H}^0(\Z_2;A))$$
be as in
Proposition {\rm \ref{linktw}~(i),(v)}, so that
$$d^*\delta d \in {\rm Quad}(B_1)~,~\phi - \phi \delta \phi \in
{\rm Quad}(B^0)$$
and
$$\phi \in Q_0(B,\beta)~=~
\dfrac{{\rm ker}(J_{\delta}:{\rm Sym}(B^0) \to {\rm Sym}(B^0)/{\rm Quad}(B^0))}
{{\rm im}((d^*)^{\%}:{\rm Quad}(B^1) \to {\rm Sym}(B^0))}$$
classifies $(T,\lambda,\mu;U)$. Lift $\beta:U \to \widehat{H}^0(\Z_2;A)$
to an $A$-module morphism $g:B_0 \to A^r$ such that
$$gd(B_1) \subseteq 2A^r~,~\delta~=~g^*X g \in \widehat{H}^0(\Z_2;S(B^0),T)~=~
{\rm Sym}(B^0)/{\rm Quad}(B^0)~.$$
The linking Arf invariant is
$$(T,\lambda,\mu;U)~=~g\phi g^* \in Q_0(B^A,\beta^A)~.$$
{\rm (iv)} For any $M=(m_{ij}) \in {\rm Sym}_r(A)$ with $m_{ij} \in 2A$
$$M-MXM~=~2(M/2-2(M/2)X(M/2)) \in {\rm Quad}_r(A)$$
and so $M$ represents an element $M \in Q_0(B^A,\beta^A)$.
The invertible matrix
$$\begin{pmatrix} -M/2 & I \\ I & 0 \end{pmatrix} \in M_{2r}(A)$$
is such that
$$\begin{array}{l}
\begin{pmatrix} -M/2 & I \\ I & 0 \end{pmatrix}
\begin{pmatrix} 0 & 2I \\ 2I & M \end{pmatrix} ~=~
\begin{pmatrix} 2I & 0 \\ 0 & 2I \end{pmatrix}~,\\[3ex]
\begin{pmatrix} I & 0 \\ -2X & I-XM \end{pmatrix}
\begin{pmatrix} -M/2 & I \\ I & 0 \end{pmatrix}~=~
\begin{pmatrix} -M/2 & I \\ I & -2X \end{pmatrix}
\end{array}$$
so that $(Q_M,\psi_M;F_M,G_M)$
is isomorphic to the $(-1)$-quadratic $S$-formation
$$(Q'_M,\psi'_M;F'_M,G'_M)~=~(H_-(A^{2r});A^{2r},
{\rm im}\bigg(\begin{pmatrix}
\begin{pmatrix} -M/2 & I \\ I & -2X \end{pmatrix} \\[3ex]
\begin{pmatrix} 2I& 0 \\ 0&2I \end{pmatrix}
\end{pmatrix}:A^{2r} \to A^{2r} \oplus A^{2r}\bigg))~,$$
corresponding to the nonsingular quadratic linking form over $(A,S)$
$$(T'_M,\lambda'_M,\mu'_M)~=~((A_2)^r \oplus (A_2)^r,
\begin{pmatrix} -M/4& I/2 \\ I/2&0 \end{pmatrix},
\begin{pmatrix} -M/4 \\ -X\end{pmatrix})$$
with $2T'_M=0$, and $U'_M=0 \oplus (A_2)^r\subset T'_M$ a
lagrangian for the symmetric linking form $(T'_M,\lambda'_M)$.
The linking Arf invariant of $(T'_M,\lambda'_M,\mu'_M;U'_M)$ is
$$(T'_M,\lambda'_M,\mu'_M;U'_M)~=~M \in Q_0(B^A,\beta^A)~.$$
\end{Theorem}
\begin{proof} (i)
$\widehat{H}^0(\Z_2;A)$ has an $S^{-1}A$-contractible f.g. free $A$-module resolution
$$\xymatrix{0 \ar[r] & A^r \ar[r]^-2 & A^r \ar[r]^-x & \widehat{H}^0(\Z_2;A)
\ar[r] & 0~.}$$
The exact sequence for $\widehat{L}^0(A)$ is given by
the exact sequence of Proposition \ref{Sbraid} (iii)
$$\dots \to L^{4*+1}(A,S) \to Q_0(B^A,\beta^A) \to L_0(A,S) \to L^{4*}(A,S) \to \dots$$
and the isomorphism $Q_0(B^A,\beta^A) \cong \widehat{L}^{4*}(A)$.\\
(ii) The isomorphism
$$Q_0(B^A,\beta^A) \to \widehat{L}^{4*}(A)~;~M \mapsto (T_M,\lambda_M,\mu_M;U_M)$$
is given by Proposition \ref{cob}.\\
(iii) Combine (ii) and Proposition \ref{linktw}.\\
(iv) By construction.
\end{proof}

\section{Application to $\unil$}

\subsection{Background}

The topological context for the unitary nilpotent $L$-groups
$\unil_\ast$ is the following. Let $N^n$ be a closed connected manifold
together with a decomposition into $n$-dimensional connected submanifolds
$N_-, N_+ \subset N$ such that
$$ N~=~N_- \cup N_+ $$
and
$$ N_\cap~=~N_- \cap N_+~=~\partial N_-~=~\partial N_+ \subset N $$
is a connected $(n-1)$-manifold with
$\pi_1 (N_\cap) \rightarrow \pi_1 (N_\pm)$ injective. Then
$$ \pi_1 (N)~=~\pi_1 (N_-) \ast_{\pi_1 (N_\cap)} \pi_1 (N_+) $$
with $\pi_1 (N_\pm) \rightarrow \pi_1 (N)$ injective. Let
$M$ be an $n$-manifold. A homotopy equivalence $f:M\rightarrow N$ is
called {\it splittable along $N_\cap$} if it is homotopic to a map
$f'$, transverse regular to $N_\cap$ (whence $f'^{-1} (N_\cap)$ is an
$(n-1)$-dimensional submanifold of $M$), and whose restriction
$f'^{-1} (N_\cap)\rightarrow N_\cap$, and a fortiori also
$f'^{-1} (N_\pm)\rightarrow N_\pm$, is a homotopy equivalence.

We ask the following question~:~ given a simple homotopy equivalence
$f:M\rightarrow N,$ when is $M$ $h$-cobordant to a manifold $M'$ such
that the induced homotopy equivalence $f':M' \rightarrow N$ is
splittable along $N_\cap$? The answer is given by Cappell
\cite{cappell3}, \cite{cappell4}~:~ the problem has a positive solution
if and only if a Whitehead torsion obstruction
$$\overline{\Phi}(\tau (f))\in \widehat{H}^n (\Z_2; \ker (\widetilde{K}_0 (A)
\rightarrow \widetilde{K}_0 (B_+) \oplus
\widetilde{K}_0 (B_-)))$$
(which is 0 if $f$ is simple) and an algebraic $L$-theory obstruction
$$ \chi^h (f) \in \unil_{n+1} (A; \Na_-, \Na_+) $$
vanish, where
$$A~=~\Z [\pi_1(N_\cap)]~,~B_{\pm}~=~\Z[\pi_1(N_{\pm})]~,~
\Na_\pm~=~B_{\pm}-A~.$$

The groups $\unil_\ast(A; \Na_-, \Na_+)$ are $4$-periodic and 2-primary,
and  vanish if the inclusions
$\pi_1 (N_\cap) \hookrightarrow \pi_1 (N_\pm)$ are square root closed.
The groups $\unil_\ast(\Z; \Z, \Z)$ arising from the expression of
the infinite dihedral group as a free product
$$D_{\infty}~=~\Z_2*\Z_2$$
are of particular interest. Cappell \cite{cappell1} showed that
$$\unil_{4k+2}(\Z;\Z,\Z)~=~
\unil_{4k+2} (\Z; \Z [\Z_2 - \{1\}], \Z [\Z_2 - \{1\}])$$
contains $(\Z_2)^\infty,$ and deduced that there is a manifold
homotopy equivalent to the connected sum
$\real \proj^{4k+1} \# \real \proj^{4k+1}$ which does not have a compatible
connected sum decomposition. With
$$B~=~\Z [\pi_1 (N)]~=~B_1*_AB_2$$
the map
$$ \unil_{n+1} (A;\Na_-, \Na_+) \longrightarrow L_{n+1} (B) $$
given by sending the splitting obstruction $\chi^h (f)$ to
the surgery obstruction of an $(n+1)$-dimensional normal map between
$f$ and a split homotopy equivalence, is a split monomorphism, and
$$L_{n+1}(B)~=~L^K_{n+1}(A \to B_+ \cup B_-) \oplus \unil_{n+1}(A; \Na_-, \Na_+)$$
with $K= \ker (\widetilde{K}_0 (A)
\rightarrow \widetilde{K}_0 (B_+)\oplus \widetilde{K}_0 (B_-))$.
Farrell \cite{farrell1} established a factorization of this map as
$$  \unil_{n+1} (A;\Na_-, \Na_+) \longrightarrow
    \unil_{n+1} (B; B, B) \longrightarrow L_{n+1} (B)~. $$
Thus the groups $\unil_n (A; A, A)$ for any ring $A$ with involution
acquire special importance, and we shall use the usual abbreviation
$$ \unil_n (A)~=~\unil_n (A; A, A)~.  $$
Cappell \cite{cappell1}, \cite{cappell2}, \cite{cappell3} proved that
$\unil_{4k} (\Z)=0$ and that $\unil_{4k+2} (\Z)$ is infinitely
generated.  Farrell \cite{farrell1} showed that for any ring $A$,
$4\unil_\ast (A)=0$.  Connolly and Ko\'zniewski \cite{connkoz} obtained
$\unil_{4k+2} (\Z)=\bigoplus_1^\infty \Z_2.$

For any ring with involution $A$ let $NL_\ast$ denote the
$L$-theoretic analogues of the nilpotent $K$-groups
$$NK_\ast (A)~=~\ker (K_\ast (A[x])\rightarrow K_\ast (A))~,$$
that is
$$NL_\ast (A)~=~\ker (L_\ast (A[x])\rightarrow L_\ast (A))$$
where $A[x]\rightarrow A$ is the augmentation map $x\mapsto 0$.
Ranicki \cite[7.6]{ranicki2} used the geometric interpretation of
$\unil_*(A)$ to identify $NL_*(A)=\unil_*(A)$ in the case when
$A=\Z[\pi]$ is the integral group ring of a finitely presented
group $\pi$. The following was obtained by pure algebra~:

\begin{Proposition} \label{equ.connran}
{\rm (Connolly and Ranicki \cite{connran})}
For any ring with involution $A$
$$\unil_\ast (A)~\cong~NL_\ast (A)~.$$
\hfill$\qed$
\end{Proposition}

It was further shown in \cite{connran} that $\unil_1 (\Z)=0$ and
$\unil_3 (\Z)$ was computed up to extensions,
thus showing it to be infinitely generated.

Connolly and Davis \cite{conndav} related $\unil_3 (\Z)$ to
quadratic linking forms over $\Z [x]$ and computed the Grothendieck
group of the latter. By Proposition \ref{equ.connran}
$$ \unil_3 (\Z)~\cong~\ker (L_3 (\Z [x])\rightarrow L_3 (\Z))
   ~=~L_3 (\Z [x])~, $$
using the classical fact $L_3 (\Z)=0.$ From a diagram chase one gets
$$ L_3 (\Z [x])~\cong~\ker (L_0 (\Z [x],(2)^{\infty}) \to  L_0 (\Z,(2)^{\infty}))~. $$
By definition, $L_0(\Z [x], (2)^{\infty})$ is the Witt group
of nonsingular quadratic linking forms $(T,\lambda,\mu)$
over $(\Z [x],(2)^{\infty})$, with $2^nT=0$ for some $n \geqslant 1$.
Let $\La (\Z [x], 2)$ be a similar Witt group, the difference being that
the underlying module $T$ is required to satisfy $2T=0.$
The main results of \cite{conndav} are
$$ L_0(\Z [x], (2)^{\infty})~\cong~\La (\Z [x], 2) $$
and
$$ \La (\Z [x], 2)~ \cong~
  \frac{x\Z_4[x]}{ \{ 2(p^2 -p)\,\vert\, p\in x\Z_4 [x] \} }   \oplus \Z_2[x]~. $$

By definition, a ring $A$ is {\it $0$-dimensional} if it is hereditary and
noetherian, or equivalently if every submodule of a f.g. projective
$A$-module is f.g. projective. In particular, a Dedekind ring $A$
is 0-dimensional. The symmetric and hyperquadratic $L$-groups of a
0-dimensional $A$ are 4-periodic
$$L^n(A)~=~L^{n+4}(A)~,~\widehat{L}^n(A)~=~\widehat{L}^{n+4}(A)~.$$

\begin{Proposition} \label{prop.uq8}
{\rm (Connolly and Ranicki \cite{connran})}
For any $0$-dimensional ring $A$ with involution
$$Q_{n+1} (B^{A[x]}, \beta^{A[x]})~=~
Q_{n+1} (B^A, \beta^A)\oplus \unil_n (A) ~~(n \in \Z)~.$$
\end{Proposition}
\begin{proof} For any ring with involution $A$ the inclusion
$A \to A[x]$ and the augmentation $A[x] \to A;x \mapsto 0$
determine a functorial splitting of the exact sequence
$$\dots \to L_n(A[x]) \to L^n(A[x]) \to \widehat{L}^n(A[x])
\to L_{n-1}(A[x]) \to \dots$$
as a direct sum of the exact sequences
$$\begin{array}{l}
\dots \to L_n(A) \to L^n(A) \to \widehat{L}^n(A)
\to L_{n-1}(A) \to \dots~,\\[1ex]
\dots \to NL_n(A) \to NL^n(A) \to N\widehat{L}^n(A)
\to NL_{n-1}(A) \to \dots~.
\end{array}$$
with $\widehat{L}^{n+4*}(A)=Q_n(B^A,\beta^A)$. It is proved in \cite{connran}
that for a 0-dimensional $A$
$$L^n(A[x])~=~L^n(A)~,~NL^n(A)~=~0~,~
N\widehat{L}^{n+1}(A)~=~NL_n(A)~=~{\rm UNil}_n(A)~.$$
\end{proof}

\begin{Example} {\rm Proposition \ref{prop.uq8} applies to $A=\Z$, so that
$$Q_{n+1} (B^{\Z[x]}, \beta^{\Z[x]})~=~
Q_{n+1} (B^\Z, \beta^\Z)\oplus \unil_n (\Z)$$
with $Q_*(B^\Z, \beta^\Z)=\widehat{L}^*(\Z)$ as
given by Example \ref{expl.QBZ}.}\hfill$\qed$
\end{Example}

\subsection{The Computation of $Q_*(B^{A[x]},\beta^{A[x]})$ for 1-even $A$}

We shall now compute the groups
$$\widehat{L}^n(A[x])~=~Q_n(B^{A[x]},\beta^{A[x]})~~(n(\bmod\,4))$$
for a 1-even ring $A$. The special case $A=\Z$ computes
$$\widehat{L}^n(\Z[x])~=~ Q_n(B^{\Z[x]},\beta^{\Z[x]})~=~
\widehat{L}^n(\Z) \oplus {\rm UNil}_{n-1}(\Z)~.$$

\begin{Proposition} \label{qba}
The universal chain bundle over $A[x]$ is given by
$$(B^{A[x]},\beta^{A[x]})~=~
\bigoplus\limits^{\infty}_{i=-\infty}(C(X),\gamma(X))_{*+2i}$$
with $(C(X),\gamma(X))$ the chain bundle over $A[x]$ given by the
construction of {\rm (}\ref{cx}{\rm )} for
$$X~=~\begin{pmatrix} 1 & 0 \\ 0 & x \end{pmatrix} \in {\rm Sym}_2(A[x])~.$$
The twisted quadratic $Q$-groups of $(B^{A[x]},\beta^{A[x]})$ are
$$\begin{array}{l}
Q_n(B^{A[x]},\beta^{A[x]})~=\\[2ex]
\begin{cases}
Q_0(C(X),\gamma(X))~=~
\dfrac{\{M \in {\rm Sym}_2(A[x])\,\vert\, M-MXM \in {\rm Quad}_2(A[x])\}}
{4{\rm Quad}_2(A[x]) + \{2(N+N^t)-N^tXN\,\vert\,N \in M_2(A[x])\}}&
\hbox{\it if $n=0$}\\[3ex]
{\rm im}(N_{\gamma(X)}:Q_1(C(X),\gamma(X)) \to Q^1(C(X)))~=~{\rm ker}(J_{\gamma(X)}:
Q^1(C(X)) \to \widehat{Q}^1(C(X)))\\[1ex]
=~\dfrac{\{N \in M_2(A[x])\,\vert\,N+N^t \in 2{\rm Sym}_2(A[x]),
\dfrac{1}{2}(N+N^t)-N^tXN \in {\rm Quad}_2(A[x])\}}{2M_2(A[x])}&\text{if $n=1$}\\[2ex]
0&\hbox{\it if $n=2$}\\[2ex]
Q_{-1}(C(X),\gamma(X))~=~\dfrac{{\rm Sym}_2(A[x])}
{{\rm Quad}_2(A[x]) + \{L-LXL\,\vert\,L \in {\rm Sym}_2(A[x])\}}&
\hbox{\it if $n=3$}~.
\end{cases}
\end{array}$$
\end{Proposition}
\begin{proof}  A special case of Theorem \ref{equ.ba-ba0},
noting that by Proposition \ref{2even} $A[x]$ is 2-even, with
$\{1,x\}$ an $A_2[x]$-module basis for $\widehat{H}^0(\Z_2;A[x])$.
\end{proof}

Our strategy for computing $Q_*(B^{A[x]},\beta^{A[x]})$ will
be to first compute $Q_*(C(1),\gamma(1))$, $Q_*(C(x),\gamma(x))$
and then to compute $Q_*(C(X),\gamma(X))$ for
$$(C(X),\gamma(X))~=~(C(1),\gamma(1)) \oplus (C(x),\gamma(x))$$
using the exact sequence given by Proposition \ref{seq} (ii)
$$\begin{array}{l}
\dots \to H_{n+1}(C(1)\otimes_{A[x]}C(x))\xymatrix{\ar[r]^-{\partial}&}
Q_n(C(1),\gamma(1)) \oplus Q_n(C(x),\gamma(x))\\[1ex]
\hphantom{\dots \to H_{n+1}(C(1)\otimes_{A[x]}C(x))}
\to Q_n(C(X),\gamma(X)) \to H_n(C(1)\otimes_{A[x]}C(x)) \to \dots~.
\end{array}$$
The connecting maps $\partial$ have components
$$\begin{array}{l}
\partial(1)~:~H_{n+1}(C(1)\otimes_{A[x]}C(x)) \to
\widehat{Q}^{n+1}(C(1)) \to Q_n(C(1),\gamma(1))~;\\[1ex]
\hphantom{\partial(1)~:~}
(f(1):C(x)^{n+1-*} \to C(1)) \mapsto (0,\widehat{f(1)}^{\%}(S^{n+1}\gamma(x)))~,\\[1ex]
\partial(x)~:~H_{n+1}(C(1)\otimes_{A[x]}C(x)) \to
\widehat{Q}^{n+1}(C(x)) \to Q_n(C(x),\gamma(x))~;\\[1ex]
\hphantom{\partial(1)~:~}
(f(x):C(1)^{n+1-*} \to C(x)) \mapsto (0,\widehat{f(x)}^{\%}(S^{n+1}\gamma(1)))~.
\end{array}$$

\begin {Proposition} \label{Qbaba}
{\rm (i)} The twisted quadratic $Q$-groups
$$Q_n(C(1),\gamma(1))~=\\[2ex]
\begin{cases}
\dfrac{A[x]}{2A[x] + \{a-a^2\,\vert\,a \in A[x]\}}
&\text{if $n=-1$}\\[2ex]
\dfrac{\{a \in A[x]\,\vert\, a-a^2 \in 2A[x]\}}
{8A[x] + \{4b-4b^2\,\vert\,b \in A[x]\}}
&\text{if $n=0$}\\[2ex]
\dfrac{\{a \in A[x]\,\vert\,a-a^2 \in 2A[x]\}}{2A[x]}&\text{if $n=1$}
\end{cases}$$
{\rm (}as given by Theorem  \ref{expl.Q0fchi}{\rm )} are such that
$$Q_n(C(1),\gamma(1))~\cong~\begin{cases}
A_2[x]&\text{if $n=-1$}\\
A_8\oplus A_4[x] \oplus A_2[x]&\text{if $n=0$}\\
A_2&\text{if $n=1$}
\end{cases}$$
with isomorphisms
$$\begin{array}{l}
f_{-1}(1)~:~Q_{-1}(C(1),\gamma(1)) \to A_2[x]~;~
\sum\limits^{\infty}_{i=0}a_ix^i
\mapsto a_0+\sum\limits^{\infty}_{i=0}(\sum\limits^{\infty}_{j=0}a_{(2i+1)2^j})x^{i+1}~,\\[1ex]
f_0(1)~:~Q_0(C(1),\gamma(1)) \to A_8 \oplus A_4[x] \oplus A_2[x]~;\\[1ex]
\hphantom{Q_0(C(1),\gamma(1))} \sum\limits^{\infty}_{i=0}a_ix^i \mapsto
(a_0,\sum\limits^{\infty}_{i=0}(\sum\limits^{\infty}_{j=0}a_{(2i+1)2^j}/2)x^i,
\sum\limits^{\infty}_{k=0}(a_{2k+2}/2)x^k)\\[1ex]
f_1(1)~:~Q_1(C(1),\gamma(1)) \to A_2~;~a~=~\sum\limits^{\infty}_{i=0}a_ix^i \mapsto a_0~.
\end{array}$$
The connecting map components $\partial(1)$ are given by
$$\begin{array}{l}
\partial(1)~:~H_1(C(1)\otimes_{A[x]}C(x))~=~A_2[x] \to Q_0(C(1),\gamma(1))~;~
c \mapsto (0,2c,0)~,\\[1ex]
\partial(1)~:~H_0(C(1)\otimes_{A[x]}C(x))~=~A_2[x] \to Q_{-1}(C(1),\gamma(1));~
c \mapsto cx.
\end{array}$$
{\rm (ii)} The twisted quadratic $Q$-groups
$$Q_n(C(x),\gamma(x))~=\\[2ex]
\begin{cases}
\dfrac{A[x]}{2A[x] + \{a-a^2x\,\vert\,a \in A[x]\}}
&\text{if $n=-1$}\\[2ex]
\dfrac{\{a \in A[x]\,\vert\, a-a^2x \in 2A[x]\}}
{8A[x] + \{4b-4b^2x\,\vert\,b \in A[x]\}}
&\text{if $n=0$}\\[2ex]
\dfrac{\{a \in A[x]\,\vert\,a-a^2x \in 2A[x]\}}{2A[x]}&\text{if $n=1$}
\end{cases}$$
{\rm (}as given by Theorem  \ref{expl.Q0fchi}{\rm )} are such that
$$Q_n(C(x),\gamma(x))~\cong~\begin{cases}
A_2[x]&\text{if $n=-1$}\\
A_4[x] \oplus A_2[x]&\text{if $n=0$}\\
0&\text{if $n=1$}
\end{cases}$$
with isomorphisms
$$\begin{array}{l}
f_{-1}(x)~:~Q_{-1}(C(x),\gamma(x)) \to A_2[x]~;~a~=~\sum\limits^{\infty}_{i=0}a_ix^i
\mapsto \sum\limits^{\infty}_{i=0}(\sum\limits^{\infty}_{j=0}
a_{(2i+1)2^j-1})x^i\\[1ex]
f_0(x)~:~Q_0(C(x),\gamma(x)) \to A_4[x] \oplus A_2[x]~;\\[1ex]
\hphantom{Q_0(C(x),\gamma(x))}
\sum\limits^{\infty}_{i=0}a_ix^i \mapsto
(\sum\limits^{\infty}_{i=0}(\sum\limits^{\infty}_{j=0}a_{(2i+1)2^j-1}/2)x^i,
\sum\limits^{\infty}_{k=0}(a_{2k+1}/2)x^k)~.
\end{array}$$
The connecting map components $\partial(x)$ are given by
$$\begin{array}{l}
\partial(x)~:~H_1(C(1)\otimes_{A[x]}C(x))~=~A_2[x] \to Q_0(C(x),\gamma(x))~
;~ c \mapsto (2c,0)~,\\[1ex]
\partial(x)~:~H_0(C(1)\otimes_{A[x]}C(x))~=~A_2[x] \to Q_{-1}(C(x),\gamma(x))~
;~ c \mapsto c~.
\end{array}$$
\end{Proposition}
\begin{proof} (i) We start with $Q_1(C(1),\gamma(1))$.
A polynomial $a(x)=\sum\limits^{\infty}_{i=0}a_ix^i\in A[x]$
is such that $a(x)-a(x)^2 \in 2A[x]$ if and only if
$$a_{2i+1}~,~a_{2i+2}-(a_{i+1})^2 \in 2A~~(i \geqslant 0)~,$$
if and only if $a_k \in 2A$ for all $k \geqslant 1$, so that
$f_1(1)$ is an isomorphism.\\
\indent Next, we consider $Q_{-1}(C(1),\gamma(1))$.
A polynomial $a(x)=\sum\limits^{\infty}_{i=0}a_ix^i \in A[x]$ is such that
$$a(x) \in 2A[x]+\{b(x)-b(x)^2\,\vert\,b(x) \in A[x]\}$$
if and only if there exist $b_1,b_2,\dots \in A$ such that
$$a_0~=~0~,~a_1~=~b_1~,~a_2~=~b_2-b_1~,~a_3~=~b_3~,~a_4~=~b_4-b_2~,~\dots
\in A_2~,$$
if and only if
$$a_0~=~\sum\limits^{\infty}_{j=0}a_{(2i+1)2^j} ~=~0 \in A_2~~(i \geqslant 0)$$
(with $b_{(2i+1)2^j}=\sum\limits^j_{k=0}a_{(2i+1)2^k} \in A_2$ for any
$i,j \geqslant 0$).
Thus $f_{-1}(1)$ is well-defined and injective.  The morphism
$f_{-1}(1)$ is surjective, since
$$\sum\limits^{\infty}_{i=0}c_ix^i~=~
f_{-1}(1)(c_0+\sum\limits^{\infty}_{i=0}c_{i+1}x^{2i+1})
\in A_2[x] ~~(c_i \in A)~.$$
The map $\widehat{Q}^1(C(1))\to Q_0(C(1),\gamma(1))$ is given by
$$\begin{array}{l}
\widehat{Q}^1(C(1))~=~A_2[x] \to
Q_0(C(1),\gamma(1))~=~A_8 \oplus A_4[x] \oplus
A_2[x]~;\\[1ex]
\hphantom{\widehat{Q}^1(C(1))~=~A_2[x]}
a=\sum\limits^{\infty}_{i=0}a_ix^i \mapsto (4a_0,
\sum\limits^{\infty}_{i=0}(\sum\limits^{\infty}_{j=0}2a_{(2i+1)2^j})x^i,0)~.
\end{array}$$
If $a=c^2x$ for $c=\sum\limits^{\infty}_{i=0}c_ix^i \in A_2[x]$ then
$$(4a_0,\sum\limits^{\infty}_{i=0}(\sum\limits^{\infty}_{j=0}2a_{(2i+1)2^j})x^i)
~=~ (0,2c) \in A_8 \oplus A_4[x]~,$$
so that the composite
$$\partial(1)~:~
H_1(C(1)\otimes_{A[x]}C(x))=A_2[x] \to \widehat{Q}^1(C(1))\to Q_0(C(1),\gamma(1))
=A_8 \oplus A_4[x] \oplus \Z_2[x]$$
is given by $c \mapsto (0,2c,0)$.

\indent Next, we consider $Q_0(C(1),\gamma(1))$.
A polynomial $a(x) \in A \oplus 2xA[x]$ is such that
$$a(x) \in 8 A[x] + \{4(b(x)-b(x)^2)\,\vert\,b(x) \in A[x]\}$$
if and only if there exist $b_1,b_2,\dots \in A$ such that
$$a_0~=~0~,~a_1~=~4b_1~,~a_2~=~4(b_2-b_1)~,~a_3~=~4b_3~,~a_4~=~4(b_4-b_2)~,~
\dots \in A_8~,$$
if and only if
$$\begin{array}{l}
a_1 ~=~ a_2 ~=~ a_3 ~=~ a_4 ~=~ \dots ~=~0 \in A_4~,\\[1ex]
a_0~=~\sum\limits^{\infty}_{j=0}a_{(2i+1)2^j}~=~0 \in A_8~~(i \geqslant 0)~.
\end{array}$$
Thus $f_0(1)$ is well-defined and injective. The morphism $f_0(1)$ is surjective, since
$$\begin{array}{ll}
(a_0,\sum\limits^{\infty}_{i=0}b_ix^i,\sum\limits^{\infty}_{i=0}c_ix^i)&=~
f_0(1)(a_0+2\sum\limits^{\infty}_{i=0}b_ix^{2i+1}+
2\sum\limits^{\infty}_{i=0}c_ix^{2i+2})\\[1ex]
&\in A_8 \oplus A_4[x] \oplus A_2[x]~~(a_0,b_i,c_i \in A)~.
\end{array}$$
The map $\widehat{Q}^0(C(1))\to Q_{-1}(C(1),\gamma(1))$ is given by
$$\begin{array}{l}
\widehat{Q}^0(C(1))~=~A_2[x] \to  Q_{-1}(C(1),\gamma(1))~=~A_2[x]~;\\[1ex]
\hphantom{\widehat{Q}^0(C(1))~=~A_2[x]}
a=\sum\limits^{\infty}_{i=0}a_ix^i \mapsto
a_0+\sum\limits^{\infty}_{i=0}(\sum\limits^{\infty}_{j=0}a_{(2i+1)2^j})x^{i+1}~.$$
\end{array}$$
If $a=c^2x$ for $c=\sum\limits^{\infty}_{i=0}c_ix^i \in A_2[x]$ then
$$a_0+\sum\limits^{\infty}_{i=0}(\sum\limits^{\infty}_{j=0}a_{(2i+1)2^j})x^{i+1}
~=~ cx \in A_2[x]~,$$
so that the composite
$$\partial(1)~:~H_0(C(1)\otimes_{A[x]}C(x))=A_2[x] \to \widehat{Q}^0(C(1))\to Q_{-1}(C(1),\gamma(1))
=A_2[x]$$
is given by $c \mapsto cx$.\\
(ii)  We start with $Q_1(C(x),\gamma(x))$.
For any polynomial $a=\sum\limits^{\infty}_{i=0}a_ix^i \in A[x]$
$$a-a^2x~=~\sum\limits^{\infty}_{i=0}a_ix^i - \sum\limits^{\infty}_{i=0}a_ix^{2i+1} \in A_2[x]~.$$
Now $a-a^2x \in 2A[x]$ if and only if the coefficients $a_0,a_1,\dots \in A$ are
such that
$$a_0 ~=~ a_1-a_0 ~=~ a_2 ~=~ a_3-a_1 ~=~ \dots ~=~ 0 \in A_2~,$$
if and only if
$$a_0 ~=~ a_1 ~=~ a_2 ~=~ a_3 ~=~ \dots ~=~ 0 \in A_2~.$$
It follows that $Q_1(C(x),\gamma(x))=0$.

Next, $Q_{-1}(C(x),\gamma(x))$.
A polynomial $a=\sum\limits^{\infty}_{i=0}a_ix^i \in A[x]$ is such that
$$a \in 2A[x]+\{b-b^2x\,\vert\,v \in A[x]\}$$
if and only if there exist $b_0,b_1,\dots \in A$ such that
$$a_0~=~b_0~~,~~a_1~=~b_1-b_0~~,~~a_2~=~b_2~~,~~
a_3~=~b_3-b_1~~,~~a_4~=~b_4~,~\dots~ \in A_2~,$$
if and only if
$$\sum\limits^{\infty}_{j=0}a_{(2i+1)2^j-1} ~=~0 \in A_2~~(i \geqslant 0)~.$$
Thus $f_{-1}(x)$ is well-defined and injective.
The morphism $f_{-1}(x)$ is surjective, since
$$\sum\limits^{\infty}_{i=0}c_ix^i~=~
f_{-1}(x)(\sum\limits^{\infty}_{i=0}c_ix^{2i})
\in A_2[x]~~(c_i \in A)~.$$
The map $\widehat{Q}^0(C(x))\to Q_{-1}(C(x),\gamma(x))$ is given by
$$\begin{array}{l}
\widehat{Q}^0(C(x))~=~A_2[x] \to  Q_{-1}(C(x),\gamma(x))~=~A_2[x]~;\\[1ex]
\hphantom{\widehat{Q}^0(C(x))~=~A_2[x]}
b=\sum\limits^{\infty}_{i=0}b_ix^i \mapsto
\sum\limits^{\infty}_{i=0}(\sum\limits^{\infty}_{j=0}b_{(2i+1)2^j-1})x^i~.$$
\end{array}$$
If $b=c^2$ for $c=\sum\limits^{\infty}_{i=0}c_ix^i \in A_2[x]$ then
$$\sum\limits^{\infty}_{i=0}(\sum\limits^{\infty}_{j=0}b_{(2i+1)2^j-1})x^i
~=~ c \in A_2[x]~,$$
so that the composite
$$\partial(x)~:~H_0(C(1)\otimes_{A[x]}C(x))=A_2[x] \to
\widehat{Q}^0(C(x))\to Q_{-1}(C(x),\gamma(x))=A_2[x]$$
is just the identity $c \mapsto c$.

Next, $Q_0(C(x),\gamma(x))$. For any $a \in A[x]$
$$a \in 8A[x]+ \{4(b-b^2x)\,\vert\,b \in A[x]\}$$
if and only there exist $b_0,b_1,\dots \in A$ such that
$$a_0~=~4b_0~,~a_1~=~4(b_1-b_0)~,~a_2~=~4b_2~,~a_3~=~
4(b_3-b_1)~,~\dots~ \in A_8,$$
if and only if
$$\begin{array}{l}
a_0 ~=~ a_1 ~=~ a_2 ~=~ a_3 ~=~ \dots ~=~ 0 \in A_4~,\\[2ex]
\sum\limits^{\infty}_{j=0}a_{(2i+1)2^j-1}~=~0 \in A_8~~(i \geqslant 0)~.
\end{array}$$
Thus $f_0(x)$ is well-defined and injective.
The morphism $f_0(x)$ is surjective, since
$$(\sum\limits^{\infty}_{i=0}c_ix^i,\sum\limits^{\infty}_{i=0}d_ix^i)~=~
f_0(x)(\sum\limits^{\infty}_{i=0}2c_ix^{2i}+
\sum\limits^{\infty}_{i=0}2d_ix^{2i+1}) \in A_4[x] \oplus A_2[x]~~
(c_i,d_i \in A)~.$$
The map $\widehat{Q}^1(C(x))\allowbreak \to Q_0(C(x),\gamma(x))$ is given by
$$\begin{array}{l}
\widehat{Q}^1(C(x))~=~A_2[x] \to
Q_0(C(x),\gamma(x))~=~A_4[x] \oplus A_2[x]~;\\[1ex]
\hphantom{\widehat{Q}^1(C(x))~=~A_2[x]}
b=\sum\limits^{\infty}_{i=0}b_ix^i \mapsto (
\sum\limits^{\infty}_{i=0}(\sum\limits^{\infty}_{j=0}2b_{(2i+1)2^j-1})x^i,0)~.
\end{array}$$
If $b=c^2$ for $c=\sum\limits^{\infty}_{i=0}c_ix^i \in A_2[x]$ then
$$\sum\limits^{\infty}_{i=0}(\sum\limits^{\infty}_{j=0}2a_{(2i+1)2^j})x^i)
~=~2c \in A_4[x]~,$$
so that the composite
$$\partial(x)~:~H_1(C(1)\otimes_{A[x]}C(x))=A_2[x] \to
\widehat{Q}^1(C(x))\to Q_0(C(x),\gamma(x))=A_4[x] \oplus A_2[x]$$
is given by $c \mapsto (2c,0)$.
\end{proof}

We can now prove Theorem \ref{main2}~:

\begin{Theorem} \label{thm.q0bbmat}
The hyperquadratic $L$-groups of $A[x]$ for a 1-even $A$ are given by
$$\widehat{L}^n(A[x])~=~Q_n(B^{A[x]},\beta^{A[x]})~=~\begin{cases}
A_8 \oplus A_4[x] \oplus A_2[x]^3&\text{if $n\equiv 0(\bmod\, 4)$}\\[1ex]
A_2 &\text{if $n\equiv 1(\bmod\, 4)$}\\[1ex]
0&\text{if $n\equiv 2(\bmod\, 4)$}\\[1ex]
A_2[x]&\text{if $n\equiv 3(\bmod\, 4)$}~.
\end{cases}$$
{\rm (i)} For $n=0$
$$Q_0(B^{A[x]},\beta^{A[x]})~=~
\dfrac{\{M \in {\rm Sym}_2(A[x])\,\vert\, M-MXM \in {\rm Quad}_2(A[x])\}}
{4{\rm Quad}_2(A[x]) + \{2(N+N^t)-4N^tXN\,\vert\,N \in M_2(A[x])\}}~.$$
An element $M \in Q_0(B^{A[x]},\beta^{A[x]})$ is represented by
a matrix
$$M~=~\begin{pmatrix} a & b \\ b & c \end{pmatrix}
\in {\rm Sym}_2(A[x])~~(a=\sum\limits^{\infty}_{i=0}a_ix^i,
c=\sum\limits^{\infty}_{i=0}c_ix^i \in A[x])$$
with $a-a_0,b,c\in 2A[x]$.  The isomorphism
$$Q_0(B^{A[x]},\beta^{A[x]}) \xymatrix{\ar[r]^{\cong}&} \widehat{L}^0(A[x])~=~
\widehat{L}^1(A[x],(2)^{\infty})~;~
M \mapsto (T_M,\lambda_M,\mu_M;U_M)$$
sends $M$ to the Witt class of the nonsingular quadratic linking form
$(T_M,\lambda_M,\mu_M)$ over $(A[x],(2)^{\infty})$ with a lagrangian
$U_M \subset T_M$ for $(T_M,\lambda_M)$ corresponding to the
$(-1)$-quadratic $(2)^{\infty}$-formation over $A[x]$
$$\partial(M)~=~(H_-(A[x]^4);A[x]^4,{\rm im}\bigg(\begin{pmatrix}
\begin{pmatrix} I & 0 \\ -2X & I-XM \end{pmatrix} \\[3ex]
\begin{pmatrix} 0 & 2I \\ 2I & M \end{pmatrix}
\end{pmatrix}:A[x]^4 \to A[x]^4 \oplus A[x]^4\bigg))$$
{\rm (\ref{linkArfthm})}, with
$$\partial~:~Q_0(B^{A[x]},\beta^{A[x]})~=~\widehat{L}^0(A[x])
\to L_{-1}(A[x])~;~M \mapsto \partial(M)~.$$
The inverse isomorphism is defined by the
linking Arf invariant {\rm (\ref{linkArf})}.\\
\indent Writing
$$2\Delta~:~A_2[x] \to A_4[x] \oplus A_4[x]~;~d \mapsto (2d,2d)$$
there are defined isomorphisms
$$\begin{array}{l}
Q_0(B^{A[x]},\beta^{A[x]}) \xymatrix{\ar[r]^-{\cong}&}
A_8 \oplus {\rm coker}(2\Delta) \oplus A_2[x] \oplus A_2[x]~;\\[2ex]
M~=~\begin{pmatrix} a & b \\ b & c \end{pmatrix}~=~
\begin{pmatrix} a & 0 \\ 0 & c' \end{pmatrix}~(c'=c-b^2)\\[2ex]
\mapsto
\big(a_0,\big[\sum\limits^{\infty}_{i=0}(\sum\limits^{\infty}_{j=0}a_{(2i+1)2^j}/2)x^i,
\sum\limits^{\infty}_{i=0}(\sum\limits^{\infty}_{j=0}c'_{(2i+1)2^j-1}/2)x^i\big],
\sum\limits^{\infty}_{k=0}(a_{2k+2}/2)x^k,
\sum\limits^{\infty}_{k=0}(c'_{2k+1}/2)x^k\big)~,\\[2ex]
{\rm coker}(2\Delta) \xymatrix{\ar[r]^-{\cong}&}
A_4[x] \oplus A_2[x]~;~[d,e] \mapsto (d-e,d)~.
\end{array}$$
In particular $M \in Q_0(B^{A[x]},\beta^{A[x]})$ can be represented by
a diagonal matrix $\begin{pmatrix} a & 0 \\ 0 & c' \end{pmatrix}$.
{\rm (ii)} For $n=1$
$$\begin{array}{l}
Q_1(B^{A[x]},\beta^{A[x]})\\[2ex]
=~\dfrac{\{N \in M_2(A[x])\,\vert\,N+N^t
\in 2{\rm Sym}_2(A[x]),
\dfrac{1}{2}(N+N^t)-N^tXN \in {\rm Quad}_2(A[x])\}}{2M_2(A[x])}
\end{array}$$
and there is defined an isomorphism
$$Q_1(B^{A[x]},\beta^{A[x]}) \xymatrix{\ar[r]^-{\cong}&}
Q_1(B^A,\beta^A)~=~A_2~;~
N~=~\begin{pmatrix} a & b \\ c & d \end{pmatrix} \mapsto a_0~,$$
with
$$\begin{array}{l}
\partial~:~Q_1(B^{A[x]},\beta^{A[x]})~=~\widehat{L}^1(A[x])~=~A_2 \to
L_0(A[x])~;\\[1ex]
a_0 \mapsto A[x]\otimes_A (A\oplus A,
\begin{pmatrix} a_0(a_0-1)/2 & 1-2a_0 \\ 0 & -2 \end{pmatrix})~.
\end{array}$$
{\rm (iii)} For $n=2$
$$Q_2(B^{A[x]},\beta^{A[x]})~=~0~.$$
{\rm (iv)} For $n=3$
$$Q_3(B^{A[x]},\beta^{A[x]})~=~\dfrac{{\rm Sym}_2(A[x])}
{{\rm Quad}_2(A[x]) + \{M-MXM\,\vert\,M \in {\rm Sym}_2(A[x])\}}~.$$
There is defined an isomorphism
$$\begin{array}{l}
Q_3(B^{A[x]},\beta^{A[x]}) \xymatrix{\ar[r]^-{\cong}&} A_2[x]~;\\[2ex]
M~=~\begin{pmatrix} a & b \\ b & c \end{pmatrix}~=~
\begin{pmatrix} a' & 0 \\ 0 & c' \end{pmatrix}
\mapsto d_0+\sum\limits^{\infty}_{i=0}
(\sum\limits^{\infty}_{j=0}d_{(2i+1)2^j})x^{i+1}\\[2ex]
(a'=a-b^2x~,~c'=c-b^2 \in A[x]~,~d=a'+c'x=a+cx \in A_2[x])~.
\end{array}$$
The isomorphism
$$Q_3(B^{A[x]},\beta^{A[x]}) \xymatrix{\ar[r]^{\cong}&} \widehat{L}^0(A[x])~=~
\widehat{L}^3(A[x])~;~
M~=~\begin{pmatrix} a & b \\ b & c \end{pmatrix} \mapsto (K_M,\psi_M;L_M)$$
sends $M$ to the Witt class of the nonsingular $(-1)$-quadratic form over $A[x]$
$$(K_M,\psi_M)~=~
(A[x]^2 \oplus A[x]^2,\begin{pmatrix} X & 1 \\ 0 & M \end{pmatrix})$$
with a lagrangian $L_M=A[x]^2\oplus 0 \subset K_M$
for $(K_M,\psi_M-\psi_M^*)$ {\rm (\ref{genArfthm})}, and
$$\partial~:~Q_3(B^{A[x]},\beta^{A[x]})~=~\widehat{L}^3(A[x])
\to L_2(A[x])~;~M \mapsto (K_M,\psi_M)~.$$
In particular $M \in Q_3(B^{A[x]},\beta^{A[x]})$ can be represented by
a diagonal matrix $\begin{pmatrix} a & 0 \\ 0 & c' \end{pmatrix}$.
The inverse isomorphism is defined by the generalized Arf invariant
{\rm (\ref{genArf})}.
\end{Theorem}
\begin{proof} Proposition \ref{qba} expresses
$Q_n(B^{A[x]},\beta^{A[x]})$
in terms of $2 \times 2$ matrices.
We deal with the four cases separately.\\
(i) Let $n=0$. Proposition \ref{Qbaba} gives an exact sequence
$$0 \to H_1(C(1)\otimes_{A[x]}C(x)) \xymatrix@C-5pt{\ar[r]^-{\partial}&} Q_0(C(1),\gamma(1))
\oplus Q_0(C(x),\gamma(x)) \to Q_0(C(X),\gamma(X)) \to 0$$
with
$$\begin{array}{l}
H_1(C(1)\otimes_{A[x]}C(x))=A_2[x]  \to \\[1ex]
Q_0(C(1),\gamma(1))
\oplus Q_0(C(x),\gamma(x)) =(A_8\oplus A_4[x]\oplus A_2[x])
\oplus (A_4[x]\oplus A_2[x])~;\\[1ex]
\hskip150pt x \mapsto ((0,2c,0),(2c,0)
\end{array}$$
so that there is defined an isomorphism
$${\rm coker}(\partial) \xymatrix{\ar[r]^-{\cong}&} A_8 \oplus
{\rm coker}(2\Delta) \oplus A_2[x] \oplus A_2[x]~;~
(s,t,u,v,w) \mapsto (s,[t,v],u,w)~.$$
We shall define an isomorphism $Q_0(C(X),\gamma(X)) \cong {\rm coker}(\partial)$
by constructing a splitting map
$$Q_0(C(X),\gamma(X)) \to Q_0(C(1),\gamma(1))\oplus Q_0(C(x),\gamma(x))~.$$
An element in $Q_0 (C(X),\gamma(X))$ is represented by a symmetric matrix
$$ M~=~\begin{pmatrix} a & b \\ b & c \end{pmatrix} \in {\rm Sym}_2(A[x])$$
such that
$$ M - MXM ~=~\begin{pmatrix} a-a^2-b^2x & b-ab-bcx \\
b-ab-bcx & c-b^2-c^2x \end{pmatrix}
\in \qad_2 (A[x])~,$$
so that
$$a - a^2 - b^2 x ~,~c-b^2-c^2x \in 2A[x]~.$$
Given $a=\sum\limits^{\infty}_{i=0}a_ix^i \in A[x]$ let
$$d~=~{\rm max}\{i \geqslant 0 \,\vert\,a_i \notin 2A\}~(=0~{\rm if}~a\in 2A[x])$$
so that $a \in A_2[x]$ has degree $d \geqslant 0$,
$$(a_d)^2~=~a_d \neq 0 \in A_2$$
and $a-a^2 \in A_2[x]$ has degree $2d$. Thus
if $b \neq 0 \in A_2[x]$ the degree of $a-a^2=b^2x \in A_2[x]$ is
both even and odd, so $b \in 2A[x]$ and hence also $a-a^2,c-c^2x \in 2A[x]$.
It follows from $a(1-a)=0 \in A_2[x]$ that $a=0$ or $1 \in A_2[x]$,
so $a-a_0 \in 2A[x]$. Similarly, it follows from $c(1-cx)=0 \in A_2[x]$ that
$c=0 \in A_2[x]$, so $c \in 2A[x]$. The matrices defined by
$$N~=~\begin{pmatrix} 0 & -b/2 \\ 0 & 0 \end{pmatrix} \in M_2(A[x])~,~
M'~=~\begin{pmatrix} a & 0 \\ 0 & c-b^2 \end{pmatrix} \in {\rm Sym}_2(A[x])$$
are such that
$$M+2(N+N^t)-4N^tXN~=~M'\in {\rm Sym}_2(A[x])$$
and so $M=M' \in Q_0(C(X),\gamma(X))$.
The explicit splitting map is given by
$$Q_0(C(X),\gamma(X)) \to Q_0(C(1),\gamma(1))\oplus Q_0(C(x),\gamma(x))~;~
M=M'\mapsto  (a,c-b^2)~.$$
The isomorphism
$$Q_0(C(X),\gamma(X)) \xymatrix{\ar[r]^-{\cong}&} {\rm coker}(\partial)~;~
M \mapsto (a,c-b^2)$$
may now be composed with the isomorphisms given in the proof of Proposition
\ref{Qbaba} (i)
$$\begin{array}{l}
Q_0(C(1),\gamma(1)) \xymatrix{\ar[r]^-{\displaystyle{\cong}}&}
A_8 \oplus A_4[x] \oplus A_2[x]~;\\[1ex]
\hphantom{Q_0(C(1),\gamma(1))} \sum\limits^{\infty}_{i=0}d_ix^i \mapsto
(d_0,\sum\limits^{\infty}_{i=0}(\sum\limits^{\infty}_{j=0}d_{(2i+1)2^j}/2)x^i,
\sum\limits^{\infty}_{k=0}(d_{2k+2}/2)x^k)~,\\[1ex]
Q_0(C(x),\gamma(x))\xymatrix{\ar[r]^-{\displaystyle{\cong}}&} A_4[x] \oplus A_2[x]~;
\\[1ex]
\hphantom{Q_0(C(x),\gamma(x))}
\sum\limits^{\infty}_{i=0}e_ix^i \mapsto
(\sum\limits^{\infty}_{i=0}(\sum\limits^{\infty}_{j=0}e_{(2i+1)2^j-1}/2)x^i,
\sum\limits^{\infty}_{k=0}(e_{2k+1}/2)x^k)~.
\end{array}$$
(ii)  Let $n=1$.
If $N=\begin{pmatrix} a & b \\ c & d \end{pmatrix} \in M_2(A[x])$
represents an element $N \in Q_1(B^{A[x]},\beta^{A[x]})$
$$\begin{array}{l}
N+N^t~=~\begin{pmatrix} 2a & b+c \\ b+c & 2d \end{pmatrix} \in
2{\rm Sym}_2(A[x])~,\\[2ex]
\dfrac{1}{2}(N+N^t)-N^tXN ~=~\begin{pmatrix} a & (b+c)/2 \\ (b+c)/2 & d
\end{pmatrix} - \begin{pmatrix} a^2+c^2x & ab+cdx \\ ab+cdx & b^2+d^2x
\end{pmatrix}\\[2ex]
\hskip200pt \in {\rm Quad}_2(A[x])
\end{array}$$
then
$$b+c~,~a-a^2-c^2x~,~d-b^2-d^2x \in 2A[x]~.$$
If $d \notin 2A[x]$ then the degree of $d-d^2x=b^2 \in A_2[x]$ is
both even and odd, so that $d \in 2A[x]$ and hence $b,c \in 2A[x]$.
Thus $a-a^2 \in 2A[x]$ and so (as above) $a-a_0 \in 2A[x]$. It follows that
$$Q_1(B^{A[x]},\beta^{A[x]})~=~Q_1(B^{A},\beta^{A})~=~A_2~.$$
(iii) Let $n=2$. $Q_2(B^{A[x]},\beta^{A[x]})=0$ by \ref{Qbaba}.\\
(iv) Let $n=3$. Proposition \ref{Qbaba} gives an exact sequence
$$0 \to H_0(C(1)\otimes_{A[x]}C(x)) \xymatrix@C-5pt{\ar[r]^-{\partial}&}
Q_{-1}(C(1),\gamma(1))\oplus Q_{-1}(C(x),\gamma(x)) \to Q_3(C(X),\gamma(X)) \to 0$$
with
$$\begin{array}{l}
\partial~:~H_0(C(1)\otimes_{A[x]}C(x))=A_2[x] \to \\[1ex]
Q_{-1}(C(1),\gamma(1))
\oplus Q_{-1}(C(x),\gamma(x))=A_2[x] \oplus A_2[x]~;~c \mapsto (cx,c)~,
\end{array}$$
so that there is defined an isomorphism
$${\rm coker}(\partial) \xymatrix{\ar[r]^-{\cong}&} A_2[x]~;~
(a,b) \mapsto a+bx~.$$
We shall define an isomorphism $Q_3(C(X),\gamma(X)) \cong {\rm coker}(\partial)$
by constructing a splitting map
$$Q_3(C(X),\gamma(X)) \to Q_{-1}(C(1),\gamma(1))\oplus Q_{-1}(C(x),\gamma(x))~.$$
For any $M=\begin{pmatrix} a & b \\ b & c \end{pmatrix} \in {\rm Sym}_2(A[x])$
the matrices
$$L~=~\begin{pmatrix} 0 & -b \\ -b & 0 \end{pmatrix}~,~
M'~=~ \begin{pmatrix} a-b^2x & 0 \\ 0 & c-b^2 \end{pmatrix} \in {\rm Sym}_2(A[x])$$
are such that
$$M'~=~M + L - LXL \in {\rm Sym}_2(A[x])$$
so $M=M' \in Q_3(C(X),\gamma(X))$. The explicit splitting map is given by
$$Q_3(C(X),\gamma(X)) \to Q_{-1}(C(1),\gamma(1))\oplus Q_{-1}(C(x),\gamma(x))~;~
M=M'\mapsto  (a-b^2x,c-b^2)~.$$
The isomorphism
$$Q_3(C(X),\gamma(X)) \xymatrix{\ar[r]^-{\cong}&} Q_{-1}(C(1),\gamma(1))~;~
M \mapsto (a-b^2x)+(c-b^2)x~=~a+cx$$
may now be composed with the isomorphism given in the proof of Proposition
\ref{Qbaba} (ii)
$$Q_{-1}(C(1),\gamma(1))\xymatrix{\ar[r]^-{\cong}&} A_2[x]~;~
d=\sum\limits^{\infty}_{i=0}d_ix^i\mapsto d_0+\sum\limits^{\infty}_{i=0}
(\sum\limits^{\infty}_{j=0}d_{(2i+1)2^j})x^{i+1}~.$$
\end{proof}

\begin{Remark}
{\rm (i)
Substituting the computation of $Q_*(B^{\Z[x]},\beta^{\Z[x]})$ given by
Theorem \ref{thm.q0bbmat} in the formula
$$Q_{n+1}(B^{\Z[x]},\beta^{\Z[x]})~=~Q_{n+1}(B^{\Z},\beta^{\Z}) \oplus
{\rm UNil}_n(\Z)$$
recovers the computations
$${\rm UNil}_n(\Z)~=~NL_n(\Z)~=\begin{cases}
0&\text{if $n \equiv 0,1(\bmod\, 4)$}\\[1ex]
\Z_2[x]&\text{if $n \equiv 2(\bmod\, 4)$}\\[1ex]
\Z_4[x]\oplus \Z_2[x]^3&\text{if $n \equiv 3(\bmod\, 4)$}
\end{cases}$$
of Connolly and Ranicki \cite{connran} and Connolly and Davis \cite{conndav}.\\
(ii) The twisted quadratic $Q$-group
$$Q_0(B^{\Z[x]},\beta^{\Z[x]})~=~\Z_8 \oplus L_{-1}(\Z[x])~=~
\Z_8 \oplus {\rm UNil}_3(\Z)$$
fits into a commutative braid of exact sequences

$$\xymatrix@!C@C-80pt@R-10pt{
0 \ar[dr]\ar@/^2pc/[rr] && Q_0(B^{\Z[x]},\beta^{\Z[x]})
\ar[dr]\ar@/^2pc/[rr]&&L_{-1}(\Z[x]) \ar[dr]\ar@/^2pc/[rr]&&
L^{-1}(\Z[1/2][x])=0 \\
& L^0(\Z[x])=\Z  \ar[ur]\ar[dr]&& L_0(\Z[x],(2)^{\infty}) \ar[ur]\ar[dr]&&
L^{-1}(\Z[x])=0 \ar[ur]\ar[dr]& \\
L_0(\Z[x])=\Z\ar[ur]^-{\displaystyle{8}}\ar@/_2pc/[rr] &&
L^0(\Z[1/2][x])=\Z \oplus \Z_2 \ar[ur]\ar@/_2pc/[rr]  &&
\ar@/_2pc/[rr] \ar[ur] L^0(\Z[x],(2)^{\infty})=\Z_2 && 0}$$

\vskip4mm

\noindent with $L_0(\Z[x],(2)^{\infty})$ (resp.  $L^0(\Z[x],(2)^{\infty})$) the Witt group of
nonsingular quadratic (resp.  symmetric) linking forms over $(\Z[x],(2)^{\infty})$, and
$$L^0(\Z[x],(2)^{\infty}) \xymatrix{\ar[r]^-{\displaystyle{\cong}}&} \Z_2~;~
(T,\lambda) \mapsto n ~{\rm if}~\vert \Z\otimes_{\Z[x]}T \vert =2^n~.$$
The twisted quadratic $Q$-group $Q_0(B^{\Z[x]},\beta^{\Z[x]})$
is thus the Witt group of nonsingular
quadratic linking forms $(T,\lambda,\mu)$ over $(\Z[x],(2)^{\infty})$ with
$\vert \Z\otimes_{\Z[x]}T \vert = 4^m$ for some $m \geqslant 0$.
$Q_0(B^{\Z[x]},\beta^{\Z[x]})$
can also be regarded as the Witt group of nonsingular
quadratic linking forms $(T,\lambda,\mu)$ over $(\Z[x],(2)^{\infty})$ together with a
lagrangian $U \subset T$ for the symmetric linking form $(T,\lambda)$.  The
isomorphism class of any such quadruple $(T,\lambda,\mu;U)$ is an
element $\phi \in Q_0(B,\beta)$.  The chain bundle $\beta$ is
classified by a chain bundle map
$$(f,\chi)~:~(B,\beta) \to (B^{\Z[x]},\beta^{\Z[x]})$$
and the Witt class is given by the linking Arf invariant
$$(T,\lambda,\mu;U)~=~(f,\chi)_{\%}(\phi) \in
Q_0(B^{\Z[x]},\beta^{\Z[x]})~=~\Z_8 \oplus\Z_4[x] \oplus \Z_2[x]^3~.$$
\noindent (iii) Here is an explicit procedure obtaining the generalized
linking Arf invariant
$$(T,\lambda,\mu;U) \in
Q_0(B^{A[x]},\beta^{A[x]})~=~A_8 \oplus A_4[x] \oplus A_2[x]^3$$
for a nonsingular quadratic linking form $(T,\lambda,\mu)$ over
$(A[x],(2)^{\infty})$ together with a lagrangian $U \subset T$ for the
symmetric linking form $(T,\lambda)$ such that $[U]=0 \in \widetilde{K}_0(A[x])$,
for any 1-even ring $A$.\\
\indent Use a set of $A[x]$-module generators $\{g_1,g_2,\dots,g_{u}\} \subset U$
to obtain a f.g. free $A[x]$-module resolution
$$0 \to B_1  \xymatrix{\ar[r]^-{\displaystyle{d}}&} B_0=A[x]^{u}
\xymatrix@C+40pt{\ar[r]^-{\displaystyle{(g_1,g_2,\dots,g_{u})}}&} U \to 0~.$$
Let $(p_i,q_i) \in A_2[x] \oplus A_2[x]$ be the unique elements such that
$$\mu(g_i)~=~(p_i)^2+x(q_i)^2 \in \widehat{H}^0(\Z_2;A[x])~=~
A_2[x]~~(1 \leqslant i \leqslant u)~,$$
and use arbitrary lifts $(p_i,q_i) \in A[x] \oplus A[x]$ to define
$$\begin{array}{l}
b_i~=~(p_i)^2+x(q_i)^2 \in A[x]~,\\[1ex]
p~=~(p_1,p_2,\dots,p_{u})~,~q~=~(q_1,q_2,\dots,q_{u}) \in A[x]^{u}~.
\end{array}$$
The diagonal symmetric form on $B_0$
$$\beta~=~\begin{pmatrix} b_1 & 0 & \dots & 0 \\
0 & b_2 & \dots & 0 \\
\vdots & \vdots & \ddots & \vdots \\
0 & 0 & \dots & b_{u}\end{pmatrix} \in {\rm Sym}(B_0)$$
is such that
$$d^*\beta d \in {\rm Quad}(B_1) \subset {\rm Sym}(B_1)~,$$
and represents the chain bundle
$$\beta~=~\mu\vert_U \in \widehat{Q}^0(B^{-*})~=~
{\rm Hom}_A(U,\widehat{H}^0(\Z_2;A[x]))~.$$
The $A[x]$-module morphisms
$$\begin{array}{l}
f_0~=~\begin{pmatrix} p \\ q \end{pmatrix}~:~
B_0~=~A[x]^{u} \to B^{A[x]}_0~=~A[x] \oplus A[x]~;~
(a_1,a_2,\dots,a_{u}) \mapsto \sum\limits^{u}_{i=1}a_i(p_i,q_i)~,\\[2ex]
f_1~:~B_1~=~A[x]^{u} \to B^{A[x]}_1~=~A[x] \oplus A[x]~;~
a=(a_1,a_2,\dots,a_{u}) \mapsto \displaystyle{\frac{f_0d(a)}{2}}
\end{array}$$
define a chain bundle map
$$(f,0)~:~(B,\beta) \to (B^{A[x]},\beta^{A[x]})~,$$
with
$$\beta^{A[x]}_0~=~\begin{pmatrix} 1 & 0 \\ 0 & x \end{pmatrix}~:~
B^{A[x]}_0~=~A[x] \oplus A[x] \to (B_0^{A[x]})^*~=~A[x] \oplus A[x]~.$$
The $(2)^{\infty}$-torsion dual of $U$ has f.g. free $A[x]$-module resolution
$$0 \to B^0=A[x]^{u}  \xymatrix{\ar[r]^-{\displaystyle{d^*}}&} B^1
\to U\widehat{~} \to 0~.$$
Lift a set of $A[x]$-module generators
$\{h_1,h_2,\dots,h_{u}\} \subset U\widehat{~}$
to obtain a basis for $B^1$, and hence an identification $B^1=A[x]^{u}$.
Also, lift these generators to elements $\{h_1,h_2,\dots,h_{u}\} \subset T$,
so that $\{g_1,g_2,\dots,g_{u},h_1,h_2,\dots,h_{u}\} \subset T$ is a set
of $A[x]$-module generators such that
$$d^{-1}~=~(\lambda(g_i,h_j)) \in
\displaystyle{\frac{{\rm Hom}_{A[1/2][x]}(B_0[1/2],B_1[1/2])}{
{\rm Hom}_{A[x]}(B_0,B_1)}}~.$$
Lift the symmetric $u \times u$ matrix $(\lambda(h_i,h_j))$ with entries
in $A[1/2][x]/A[x]$ to a symmetric form on the f.g. free $A[1/2][x]$-module
$B^1[1/2]=A[1/2][x]^{u}$
$$\Lambda~=~(\lambda_{ij})\in {\rm Sym}(B^1[1/2])$$
such that $\lambda_{ii} \in A[1/2][x]$ has image $\mu(h_i) \in
A[1/2][x]/2A[x]$.  Let $\phi=(\phi_{ij})$ be the symmetric form on
$B^0=A[x]^{u}$ defined by
$$\phi~=~d\Lambda d^* \in {\rm Sym}(B^0) \subset {\rm Sym}(B^0[1/2])~. $$
Then $T$ has a f.g. free $A[x]$-module resolution
$$0 \to B_1\oplus B^0  \xymatrix@C+20pt
{\ar[r]^-{\displaystyle{\begin{pmatrix} 0 & d^* \\ d & \phi\end{pmatrix}}}&}
B^1 \oplus B_0 \xymatrix@C+75pt{\ar[r]^-{\displaystyle{(g_1,\dots,g_{u},h_1,\dots,h_{u})}}&}
T \to 0~,$$
and
$$\phi_{ii} - \sum\limits^{u}_{j=1}(\phi_{ij})^2b_j \in 2A[x]~.$$
The symmetric form on $(B^{A[x]}_0)^*=A[x] \oplus A[x]$ defined by
$$\begin{array}{c}
\begin{pmatrix} a & b \\ b & c \end{pmatrix}~=~
f_0 \phi f_0^* ~=~
\begin{pmatrix} \phi(p,p) & \phi(p,q) \\[1ex]
\phi(q,p) & \phi(q,q) \end{pmatrix} \in {\rm Sym}((B^{A[x]}_0)^*)\\[3ex]
(p=(p_1,p_2,\dots,p_{u}),q=(q_1,q_2,\dots,q_{u}) \in B^0=A[x]^{u})
\end{array}$$
is of the type considered in the proof of Theorem \ref{thm.q0bbmat} (i), with
$$a-a^2~=~b^2x~,~c-c^2x~=~b^2 \in A_2[x]~,~b \in 2A[x]~.$$
The Witt class is
$$\begin{array}{ll}
(T,\lambda,\mu;U)&=~(f,0)_{\%}(\phi)\\[1ex]
&=~\begin{pmatrix} a & b \\ b & c \end{pmatrix}~=~
\begin{pmatrix} a & 0 \\ 0 & c' \end{pmatrix}
\in Q_0(B^{A[x]},\beta^{A[x]})~~(c'=c-b^2)~,
\end{array}$$
with isomorphisms
$$\begin{array}{l}
Q_0(B^{A[x]},\beta^{A[x]})\xymatrix{\ar[r]^-{\displaystyle{\cong}}&}
A_8\oplus
{\rm coker}(2\Delta) \oplus A_2[x]\oplus A_2[x]~;\\[1ex]
\begin{pmatrix} a & 0 \\ 0 & c' \end{pmatrix}
\mapsto \big(a_0,\big[\sum\limits^{\infty}_{i=0}
(\sum\limits^{\infty}_{j=0}a_{(2i+1)2^j}/2)x^i,
\sum\limits^{\infty}_{i=0}(\sum\limits^{\infty}_{j=0}
c'_{(2i+1)2^j-1}/2)x^i\big],\\
\hskip172pt \sum\limits^{\infty}_{k=0}(a_{2k+2}/2)x^k,
\sum\limits^{\infty}_{k=0}(c'_{2k+1}/2)x^k\big)~,\\[2ex]
{\rm coker}(2\Delta) \xymatrix{\ar[r]^-{\displaystyle{\cong}}&} A_4[x] \oplus A_2[x]~;~
[m,n] \mapsto (m-n,m)~,
\end{array}$$
where
$$2\Delta~:~A_2[x] \to A_4[x] \oplus A_4[x] ~;~m \mapsto (2m,2m)$$
as in Theorem \ref{thm.q0bbmat}, and
$$Q_0(B^{A[x]},\beta^{A[x]})~=~A_8 \oplus
A_4[x] \oplus A_2[x]^3~.$$
For Dedekind $A$ the splitting formula of \cite{connran} gives
$${\rm UNil}_3(A)~\cong~
Q_0(B^{A[x]},\beta^{A[x]})/A_8~\cong~ A_4[x] \oplus A_2[x]^3~.$$
\hfill$\qed$}
\end{Remark}

\providecommand{\bysame}{\leavevmode\hbox to3em{\hrulefill}\thinspace}

\end{document}